%% file: fvdbltt.tex
\documentclass[dvipsnames, reqno]{amsart}
\usepackage{preamble_fvdbltt}
\RequirePackage[%
a4paper,margin=2cm,%
]{geometry}

\title[An Internal Logic of Virtual Double Categories ]
{An Internal Logic of Virtual Double Categories }
\author{Hayato Nasu}

\date{\today}

\begin{document}

\begin{abstract}
    We present a type theory called \ac{FVDblTT} designed specifically for formal category theory,
    which is a succinct reformulation of New and Licata's Virtual Equipment Type Theory (VETT).
    \ac{FVDblTT} formalizes reasoning on isomorphisms that are commonly employed in category theory.
    Virtual double categories are one of the most successful frameworks for developing formal category theory,
    and \ac{FVDblTT} has them as a theoretical foundation.
    We validate its worth as an internal language of virtual double categories 
    by providing a syntax-semantics duality between virtual double categories and specifications in \ac{FVDblTT}
    as a biadjunction.
\end{abstract}
    \maketitle

    \section{Introduction}
        \label{section:introfvdtt}
        \subfile{newsubfiles/introduction/intromain.tex}
    \subsection{Outline}
        \label{subsection:outline}
        \subfile{newsubfiles/introduction/introoutline.tex}
    \subsection{Acknowledgements}
        \label{subsection:acknowledgements}
        The author would like to thank his supervisor, Masahito Hasegawa, for his guidance and encouragement.
        Benedikt Ahrens, Nathanael Arkor, Keisuke Hoshino, Yuki Maehara, Hiroyuki Miyoshi, Paige North, and Yuta Yamamoto for their discussions, comments, and suggestions
        on the manuscript.
        The author would also like to thank the anonymous reviewers for their valuable feedback. 

    \section{Preliminaries on Virtual Double Categories}
        \label{section:vdc}
        \subfile{newsubfiles/prelim/vdc.tex}

    \section{Fibrational Virtual Double Type Theory}
        \label{sec:fvdtt}
        \subfile{newsubfiles/typetheory/typetheory.tex}
    \subsection{Syntax}
        \label{sec:syntax}
        \subfile{newsubfiles/typetheory/syntax.tex}

    \subsection{Semantics}
        \label{sec:semantics}
        \subfile{newsubfiles/typetheory/semantics.tex}
        \subfile{newsubfiles/typetheory/inddefsem.tex}
    \subsection{Protype isomorphisms}
        \label{sec:protypeiso}
        \subfile{newsubfiles/typetheory/proiso.tex}

    \section{Protype and type constructors for FVDblTT}
        \label{sec:additional}
    \subsection{Further structures in VDCs and the corresponding constructors}
        \label{subsec:additionaltype}
        \subfile{newsubfiles/constructor/tableofconstructors.tex}
        \subfile{newsubfiles/constructor/units.tex}
        \subfile{newsubfiles/constructor/others.tex}
        \subfile{newsubfiles/constructor/predlogic.tex}
    
        \section{Examples of calculus}
        \label{sec:examples}
        \subfile{newsubfiles/constructor/examples.tex}

    \section{A syntax-semantics adjunction for FVDblTT}
        \label{sec:synsemadj}
        \subfile{newsubfiles/adjunction/synsem.tex}
    \subsection{Syntactic presentation of virtual double categories}
        \label{sec:synfvd}
        \subfile{newsubfiles/adjunction/syntacticpres.tex}

    \subsection{Constructing the adjunction}
        \label{sec:new}
        \subfile{newsubfiles/adjunction/newproof.tex}

    \section{Future Work}
        \label{sec:relatedwork}
        \subfile{newsubfiles/discussion/conclusion.tex}

    \bibliographystyle{halpha}
    \bibliography{fvdbltt}

    \appendix
    \section{Cartesianness of Structured Virtual Double Categories}
        \label{sec:composition}
        \subfile{newsubfiles/appendix/cartesianstr.tex}

    \section{The derivation rules for the additional constructors}
        \label{sec:appendix1}
        \subfile{newsubfiles/appendix/cartesiansyn.tex}
    \vspace{1em}
\end{document}

%% file: newsubfiles/introduction/intromain.tex
Variants of category theory have been developed over the decades,
each with its own characteristics but sharing some basic concepts and principles.
For instance, monoidal category theory \cite{selingerSurveyGraphicalLanguages2011},
enriched category theories over monoidal categories \cite{kellyBasicConceptsEnriched2005},
internal category theories in toposes \cite{johnstoneSketchesElephantTopos2002a},
and fibered category theory \cite{streicherFiberedCategoriesJean2023}
all have well-developed theories and significant applications.
They often share several concepts, such as limits, representable functors, adjoints,
and fundamental results like the Yoneda lemma,
though there may be slight differences in their presentations.

\textit{Formal category theory} \cite{grayFormalCategoryTheory1974} is the abstract method that unifies
these various category theories.
As category theory offers us abstract results that can universally be applied to
mathematical structures,
formal category theory enables us to enjoy the universal results 
that hold for general category theories.
A comprehensive exposition of this field is given in 
\cite{libertiFORMALCATEGORYTHEORYa}.
The earliest attempt was to perform category theory in an arbitrary 2-category by pretending that it is the 2-category of categories \cite{grayFormalCategoryTheory1974}.
However, more than just 2-categories are needed to capture the big picture of category theory.
The core difficulty that one encounters in this approach
is that it does not embody the notion of presheaves,
or ``set-valued functors'' inside a 2-category.
Subsequently, many solutions have emerged to address this problem,
such as \textit{Yoneda structures} \cite{streetYonedaStructures2categories1978a}
and \textit{proarrow equipments} \cite{woodAbstractProarrows1982,woodProarrowsII1985}.

A recent and prominent approach to formal category theory is to use
\textit{\aclp{VDC}} or \textit{augmented virtual double categories} 
\cite{shulmanFramedBicategoriesMonoidal2009,koudenburgAugmentedVirtualDouble2020}.
General theories in (augmented) virtual double categories have recently been developed,
successful examples of which include
the Yoneda structures and total categories in augmented virtual double categories by Koudenburg
\cite{koudenburgAugmentedVirtualDouble2020,koudenburgFormalCategoryTheory2024}
and the theory of relative monads in virtual equipments by Arkor and McDermott \cite{arkorFormalTheoryRelative2024}.
The advantage of this framework is that it is built up with necessary components of category theory as 
primitive structures.
A virtual double category models the structure constituted by categories, functors, natural transformations,
and \textit{profunctors}, a common generalization of presheaves and copresheaves.
This allows us to capture far broader classes of category theories
since the virtual double category for a category can at least be defined
even when essential components, like presheaves or natural transformations, do not behave as well as
in the ordinary category theory.

In this paper, we provide a type theory called \textit{\textbf{\acl{FVDblTT}}} (\ac{FVDblTT}),
which is designed specifically for formal category theory and serves
as an internal language of virtual double categories.
It aims to function as a formal language to reason about category theory 
that can be applied to various category theories,
which may be used as the groundwork for computer-assisted proofs.
Arguing category theories is often divided into two parts:
one is a common argument independent of different category theories, which occasionally falls into \textit{abstract nonsense},
and the other is a specific discussion particular to a certain category theory.
What we can do with this type theory is to deal with massive proofs belonging to the former part in the formal language and 
make people focus on the latter part.
Our attempt is not the first in this direction,
as it follows New and Licata's Virtual Equipment Type Theory (VETT) 
\cite{newFormalLogicFormal2023}.
However, we design \ac{FVDblTT} based on the following desiderata that set it apart from the previous work:
\begin{enumerate}
    \item It admits a syntax-semantics duality between the category of virtual double categories (with suitable structures) and the category of syntactic presentations of them.
    \item It is built up from a plain core type theory but allows enhancement that is compatible with existing and future results in formal category theory.
    \item It allows reasoning with isomorphisms, a common practice in category theory.
\end{enumerate}
In order to explain how \ac{FVDblTT} achieves these goals,
we overview its syntax and semantics.

\subsection{Syntax and Semantics}
We start with reviewing virtual double categories.
While its name first appeared in the work of Cruttwell and Shulman \cite{cruttwellUnifiedFrameworkGeneralized2010},
the idea of virtual double categories has been studied in various forms in the past
under different names such as 
\textit{multicat\'egories} \cite{Bur71},
\textit{\textbf{fc}-multicategories} \cite{leinsterGeneralizedEnrichmentCategories2002,leinsterHigherOperadsHigher2004},
and \textit{lax double categories} \cite{dawsonPathsDoubleCategories2006}.
For these years, virtual double categories have gained the status of a guidepost for working out new category theories,
especially in the $\infty$-categorical setting \cite{gepnerEnrichedCategoriesNonsymmetric2015,riehlKanExtensionsCalculus2017,ruit2024formalcategorytheoryinftyequipments}.

A virtual double category has four kinds of data:
\textit{objects}, \textit{tight arrows}, \textit{loose arrows}, and \textit{virtual cells}.
The typical example is $\Prof$, which has categories, functors, \textit{profunctors}, and (generalized) natural transformations
as these data.
A profunctor from a category $\one{I}$ to a category $\one{J}$, written as $P(-,\bullet)\colon\one{I}\sto\one{J}$, is a functor from $\one{I}\op\times\one{J}$ to the category of sets $\Set$,
which is a common generalization of a presheaf on $\one{I}$ and a copresheaf on $\one{J}$.
One would expect these two kinds of arrows to have compositional structures,
and indeed, two profunctors $P(-,\bullet)\colon\one{I}\sto\one{J}$ and $Q(-,\bullet)\colon\one{J}\sto\one{K}$ can be composed by a certain kind of
colimits called coends in $\Set$.
However, the composition of profunctors may not always be defined within a general category theory,
for instance, an enriched category theory with the enriching base category lacking enough colimits.
Virtual cells are introduced to liberate loose arrows from their composition 
and yet to keep seizing their compositional behaviors.
As in \Cref{fig:cell2},
a virtual cell has two tight arrows, one loose arrow, and one sequence of 
loose arrows as its underlying data,
and in the case of $\Prof$, virtual cells are 
natural families with multiple inputs.
This pliability enables us to express category theoretic phenomena 
with a weaker assumption on the category theory one works with.
\begin{figure}[htbp]
    \centering
    \begin{minipage}[t]{0.6\columnwidth}
    \begin{itemize}
    \item A virtual cell in $\Prof$:
        \[
            \begin{tikzcd}[virtual]
                \one{I}_0
                \ar[d, "S"']
                \sar[r, "\alpha_1"]
                \ar[rrrd, phantom, "\mu"]
                & \one{I}_1
                \sar[r]
                & \cdots
                \sar[r, "\alpha_n"]
                & \one{I}_n
                \ar[d, "T"] \\
                \one{J}_0
                \sar[rrr, "\beta"']
                & & & \one{J}_1
            \end{tikzcd}
        \]
    \item A family of functions natural in $i_0,i_n$ and dinatural in $i_1,\dots,i_{n-1}$:
        \begin{equation*}
            \mu_{i_0,\dots,i_n}\colon \alpha_1(i_0,i_1)\times\dots\times\alpha_n(i_{n-1},i_n)
            \to\  \beta(S(i_0),T(i_n)) 
        \end{equation*}
    \item An interpretation of the proterm 
    {\footnotesize
    \begin{align*}
    \syn{x}_0:\syn{I}_0\smcl\dots\smcl\syn{x}_n:\syn{I}_n&\mid
    \syn{a}_1:\syn{\alpha}_1(\syn{x}_0\smcl\syn{x}_1)\smcl\dots\smcl
    \syn{a}_n:\syn{\alpha}_n(\syn{x}_{n-1}\smcl\syn{x}_n)\\
    &\vdash\syn{\mu}:\syn{\beta}(\syn{S}(\syn{x}_0),\syn{T}(\syn{x}_n)).
    \end{align*}
    }
    \end{itemize}
    \end{minipage}
    \caption{A virtual cell in $\Prof$ and a proterm that corresponds to it.}
    \label{fig:cell2}  
\end{figure}

Corresponding to these four kinds of entities, 
\ac{FVDblTT} has four kinds of core judgments:
\textit{types}, \textit{terms}, \textit{\textbf{pro}types}, and \textit{proterms} (\Cref{fig:judgment}).
In the semantics in the virtual double category $\Prof$,
types, terms, and protypes are interpreted as categories, functors, and \textbf{pro}functors,
while proterms are interpreted as virtual cells with the functors on both sides being identities.
We restrict the interpretation in this way in order to have the linearized presentation of virtual cells in the type theory.
This enables us to bypass diagrammatic presentations of virtual cells,
which often occupy considerable space in papers\footnote{
    This thesis is a good example of this.
}.
Nevertheless, it does not lose the expressive power because we assume the 
semantic stage to be a \emph{fibrational} virtual double category.
\begin{figure}[htbp]
    \centering
    \begin{gather*}
        \text{Type} \quad \syn{I}\ \textsf{type}\  ,\quad \\
        \text{Term} \quad \syn{\Gamma} \vdash \syn{s}:\syn{I}\  ,\quad \\
        \text{Protype} \quad \syn{\Gamma}\smcl\syn{\Delta} \vdash \syn{\alpha} \ \textsf{protype}\  ,\quad\\
        \text{Proterm} \quad \syn{\Gamma}_0\smcl\dots\smcl\syn{\Gamma}_n \mid \syn{a}_1:\syn{\alpha}_1\smcl\dots\smcl\syn{a}_n:\syn{\alpha}_n \vdash \syn{\mu}:\syn{\beta}\  ,\quad\\
        (\syn{\Gamma},\syn{\Delta},\dots \text{ are contexts like } \syn{x}_1:\syn{I}_1,\dots,\syn{x}_n:\syn{I}_n.)
    \end{gather*}
    \caption{Judgments of \ac{FVDblTT}.} 
    \label{fig:judgment}  
\end{figure}

Fibrationality is satisfied in most virtual double categories for our purposes
and is conceptually a natural assumption since it represents the possibility 
of instantiating functors $S$ and $T$ in a profunctor $\alpha(-,\bullet)$.
Furthermore, the fibrationality reflects how we practically reason about cells in the virtual double categories 
for formal category theory.
For instance, 
a virtual cell in $\Prof$ is defined as a natural family $\mu$, as in \Cref{fig:cell2},
and it only refers to the instantiated profunctor $\beta(S(-),T(\bullet))$.
Accordingly, we let the type theory describe a virtual cell as a proterm as in \Cref{fig:cell2}.
The fibrationality condition is defined in terms of universal property and assumed to hold in the semantics.
We will further assume \acp{VDC} to have suitable finite products
to interpret finite products in \ac{FVDblTT},
which alleviates the complexity of syntactical presentation.

A byproduct of this type theory is
its aspect as an all-encompassing language for predicate logic.
The double category $\Rel$ of sets, functions, relations as objects, tight arrows, and loose arrows
would also serve as the stage of the semantics of \ac{FVDblTT}.
In this approach, protypes symbolize relations (two-sided \textbf{pro}positions),
and proterms symbolize Horn formulas.
In other words, category theory based on categories, functors, and profunctors 
can be perceived as \textit{generalized logic}.
The unity of category theory and logic dates back to the work of Lawvere \cite{lawvereMetricSpacesGeneralized1973a},
in which he proposed that the theories of categories or metric spaces are generalized logic,
with the truth value sets being some closed monoidal categories.

The interpretation of \ac{FVDblTT} is summarized in \Cref{table:common}.
\begin{table*}[htbp]
    \setlength{\aboverulesep}{0pt}
    \setlength{\belowrulesep}{0pt}
    \footnotesize
    \rowcolors{2}{gray!25}{white}
    \resizebox{\textwidth}{!}{%
    \begin{tabular}{|c||c|c|}
        \toprule
        \rowcolor{gray!50}
        \textbf{Items in FVDblTT} & \textbf{Formal category theory} & \textbf{Predicate logic} \\
        \midrule\midrule
        Types $\syn{I}$ &  categories $\one{I}$ & sets $I$ \\
        Terms $\syn{x}:\syn{I}\vdash\syn{s}:\syn{J}$ & functors $S\colon\one{I}\to\one{J}$ & functions $s\colon I\to J$ \\
        Protypes $\syn{\alpha}(\syn{x}\smcl\syn{y})$ & profunctors $\alpha\colon\one{I}\sto\one{J}$ &  formulas $\alpha(x,y)\ (x\in I,y\in J)$ \\
        \begin{tabular}{@{}c@{}} Proterms  \\ $\syn{a}:\syn{\alpha}(\syn{x}\smcl\syn{y})\smcl\syn{b}:\syn{\beta}(\syn{y}\smcl\syn{z})$\\$\quad\vdash\syn{\mu}:\syn{\gamma}(\syn{x}\smcl\syn{z})$\end{tabular} & \begin{tabular}{@{}c@{}}natural transformations\\ $\mu_{x,y,z}\colon \alpha(x,y)\times \beta(y,z)\to \gamma(x,z)$\end{tabular} & \begin{tabular}{@{}c@{}} proofs of Horn clauses\\ $\alpha(x,y),\beta(y,z)\Rightarrow \gamma(x,z)$\end{tabular} \\
        Product types $\syn{I}\times\syn{J}$& product categories $\one{I}\times\one{J}$& product sets $I\times J$ \\
        Product protypes  $\syn{\alpha}\land\syn{\beta}$ & product profunctors $\alpha(x,y)\times \beta(x,y)$ & conjunctions $\alpha(x,y)\land\beta(x,y)$ \\
        \midrule
        path protype $\ide{\!}$ & hom profunctor $\one{I}(-,-)$ & equality relation $=_I$ \\
        composition protype $\odot$ & composition of profunctors by coend & composition of relations by $\exists$ \\
        \midrule
        \begin{tabular}{@{}c@{}} Protype Isomorphisms \\ $\syn{\Upsilon}:\syn{\alpha}\ccong\syn{\beta}$\end{tabular} &\begin{tabular}{@{}c@{}}natural isomorphisms\\ $\Upsilon_{x,y}\colon \alpha(x,y)\cong \beta(x,y)$\end{tabular}&\begin{tabular}{@{}c@{}} equivalence of formulas\\ $\alpha(x,y)\equiv\beta(x,y)$\end{tabular} \\
        \bottomrule
    \end{tabular} 
    }
    \vspace{0.5em}
    \caption{Interpretation of \ac{FVDblTT} in $\PROF$ and $\Rel{}$\quad
    (All rows except the last three are included in the core of \ac{FVDblTT}.)}
    \label{table:common}
\end{table*}

\subsection{Realizing the desiderata}

\myparagraph{(i) Syntax-semantics duality for \ac{VDC}}
Categorical structures have been studied as the stages for semantics.
Good examples include the Lawvere theories in categories with finite products \cite{lawvereFunctorialSemanticsAlgebraic1963},
simply typed lambda calculus in cartesian closed categories \cite{lambekIntroductionHigherOrder1986},
extensional Martin-L\"of type theory in locally cartesian closed categories 
\cite{seelyLocallyCartesianClosed1984},
and homotopy type theory in $\infty$-groupoids \cite{hofmannGroupoidInterpretationType1998,STREICHER201445}.
Thus, it has been discovered that there are dualities between syntax and categorical structures
\cite{jacobsCategoricalLogicType1999a,clairambaultBiequivalenceLocallyCartesian2014},
endorsing the principle that type theory corresponds to category theory.
It is worth noting that the above examples all started from the development of calculi,
and the corresponding categorical structures were determined.

We will define specifications for \ac{FVDblTT} and construct an adjunction
between the category of virtual double categories with some structures and the category of those specifications
whose counit is componentwise an equivalence,
which justifies the type theory as an internal language
and directly implies the soundness and completeness of the type theory. 
Here, we have proceeded in the reversed direction to the traditional developments:
knowing that virtual double categories are the
appropriate structures for formal category theory,
we extract a calculus from it.
This principle can be seen in \cite{ahrensBicategoricalTypeTheory2023}.

\myparagraph{(ii) Additional constructors}
Additional type and protype constructors are introduced to make \ac{FVDblTT}
expressive enough to describe sophisticated arguments in category theory.
For example, the hom-profunctor $\one{I}(-,\bullet)\colon\one{I}\sto\one{I}$
cannot be achieved in the core \ac{FVDblTT},
and we introduce \textit{path protype}
$\syn{x}:\syn{I}\smcl\syn{y}:\syn{I}\vdash\syn{x}\ide{\syn{I}}\syn{y} : \textsf{protype}$
as its counterpart.
Just as a variable $\syn{x}:\syn{I}$ serves as an object variable in $\one{I}$,
a provariable $\syn{a}:\syn{x}\ide{\syn{I}}\syn{y}$ serves as a morphism variable in $\one{I}$.
The introduction rule for this is similar to the path induction in homotopy type theory.
Using this constructor, one can formalize, for instance, the fully-faithfulness of a functor (\Cref{fig:fullyfaithful}),
as it is defined merely through the bahevior on the hom-sets.
Interpreting this in the virtual double categories of enriched categories,
one obtains the existing definition of fully-faithful enriched functors.
In addition, we introduce \textit{composition protype}, \textit{filler protype},
and \textit{comprehension type} in this paper,
by which one can formalize a myriad of concepts in category theory,
including (weighted) (co)limits,
pointwise Kan extensions, and the Grothendieck construction of (co)presheaves,
which is only possible with the protype constructors.

\begin{figure}[htbp]
    \begin{minipage}[t]{0.38\columnwidth}
    \captionsetup{width=\textwidth}
    \begin{framed}
    {\footnotesize
    A term $\syn{x}:\syn{I}\vdash \syn{s}(\syn{x}): \syn{J}$
    is fully faithfull if 
    the following proterm has an inverse.
    \begin{mathpar}
        \inferrule*
        {
        \inferrule*
        {
        \syn{y}:\syn{J}\mid\ \vdash \refl_{\syn{J}}:\syn{y}\ide{\syn{J}}\syn{y} 
        }
        {
        \syn{x}:\syn{I}\mid\ \vdash \refl_{\syn{J}}[\syn{s}(x)/\syn{y}]:\syn{s}(\syn{x})\ide{\syn{I}}\syn{s}(x)
        }}
        {
        {\begin{aligned}
        &\syn{x}:\syn{I}\smcl\syn{x'}:\syn{I}\mid\syn{a}:\syn{x}\ide{\syn{I}}\syn{x'}\\
        &\vdash
        \ideind{\syn{I}}\{\refl_{\syn{J}}[\syn{s}(x)/\syn{y}]\}:
        \syn{s}(\syn{x})\ide{\syn{J}}\syn{s}(\syn{x'})
        \end{aligned}}
        }
    \end{mathpar}
    Having an inverse is formulated using protype isomorphisms:
    it comes with the following protype isomorphism
    \[
        \syn{x}:\syn{I}\smcl\syn{x'}:\syn{I}\vdash
        \textsf{FF}:
        \syn{x}\ide{\syn{I}}\syn{x'}
        \ccong
        \syn{s}(\syn{x})\ide{\syn{J}}\syn{s}(\syn{x'})
    \]
    that satisfies the following equation (context is omitted).
    \[
    \textsf{FF}\{\syn{a}\}
    \equiv 
    \ideind{\syn{I}}\{\refl_{\syn{J}}[\syn{s}(\syn{x})/\syn{y}]\}
    \]
    }
    \end{framed}
    \caption{Fully faithfulness}
    \label{fig:fullyfaithful}
    \end{minipage}
    \hfill
    \begin{minipage}[t]{0.61\columnwidth}
    \begin{minipage}[t]{\columnwidth}
    \begin{framed}
    {\footnotesize
    A term $\syn{y}:\syn{J}\vdash \syn{l}(\syn{y}): \syn{K}$
    is a pointwise left Kan extension of $\syn{x}:\syn{I}\vdash \syn{s}(\syn{x}): \syn{K}$
    along $\syn{x}:\syn{I}\vdash \syn{t}(\syn{x}): \syn{J}$
    if it comes with the following protype isomorphism.
    \[
    \syn{y}:\syn{J}\smcl\syn{z}:\syn{K}\vdash 
    \textsf{Lan}:
    \syn{l}(\syn{y})\ide{\syn{K}}\syn{z}
    \ccong
    (\syn{t}(\syn{x})\ide{\syn{J}}\syn{y})\triangleright_{\syn{x}:\syn{I}}
    (\syn{s}(\syn{x})\ide{\syn{K}}\syn{z})
    \]}
    \end{framed}
    \caption{Pointwise Kan extensions}
    \label{fig:pointwisekan}
    \end{minipage}
    \begin{minipage}[b]{\columnwidth}
    \vspace{1.8em}
    \begin{framed}
    {\footnotesize
    A pointwise left Kan extension $\syn{l}(\syn{y})$
    of $\syn{s}(\syn{x})$
    along a fully faithful functor $\syn{t}(\syn{x})$
    admits an isomorphism $\syn{l}(\syn{t}(\syn{x}))\ccong\syn{s}(\syn{x})$.
    
    \textbf{Proof.} (Contexts are omitted.)
    \begin{align*}
    \syn{l}(\syn{t}(\syn{x}'))\ide{\syn{K}}\syn{z}
    &\ccong
    (\syn{t}(\syn{x})\ide{\syn{J}}(\syn{t}(\syn{x}')))\triangleright_{\syn{x}:\syn{I}}
    (\syn{s}(\syn{x})\ide{\syn{K}}\syn{z}) 
    &
    \text{($\textsf{Lan}$)} \\
    &\ccong
    (\syn{x}\ide{\syn{I}}\syn{x'})\triangleright_{\syn{x}:\syn{I}}
    (\syn{s}(\syn{x})\ide{\syn{K}}\syn{z})
    &
    \hspace{-3em}
    \text{($\textsf{FF}\sinv\triangleright_{\syn{x}:\syn{I}}(\syn{s}(\syn{x})\ide{\syn{K}}\syn{z})$)} \\
    &\ccong
    \syn{s}(\syn{x'})\ide{\syn{K}}\syn{z}
    &
    \hspace{-3em}
    \text{($\textsf{Yoneda}$ \Cref{example:ninja-yoneda})}
    \end{align*}
    }
    \end{framed}
    \caption{Pointwise Kan extensions along fully faithful functors}
    \label{fig:pointwisekan2}
    \end{minipage}
    \end{minipage}
\end{figure}

\myparagraph{(iii) Isomorphism reasoning}
We will enhance our type theory with \textit{protype isomorphisms},
a new kind of judgment for isomorphisms between protypes.
\[
    \text{Protype Isomorphism} \quad \syn{\Gamma}\smcl\syn{\Delta} \vdash \syn{\Upsilon}:\syn{\alpha}\ccong\syn{\beta}\ .
\]
They serve as a convenient gadget for up-to-isomorphism reasoning that is ubiquitous in category theory.
One often proves two things are isomorphic
by constructing some pieces of mutual inverses and then
combining them to form the intended isomorphism.
We bring this custom into the type theory as protype isomorphisms,
interpreted as isomorphisms between profunctors, \textit{i.e.}, an invertible natural transformation between profunctors.
For instance, pointwise Kan extensions are concisely defined using protype isomorphisms (\Cref{fig:pointwisekan}).
Protype isomorphisms capture isomorphisms between functors as well since isomorphisms between functors $F,G\colon\one{I}\to\one{J}$
correspond to natural isomorphisms between $\one{J}(F-,\bullet),\one{J}(G-,\bullet)\colon\one{I}\sto\one{J}$ according to the Yoneda lemma.
A formal proof of a well-known fact that a pointwise left Kan extension along 
a fully faithful functor admits an isomorphism to the original functor 
can be given by isomorphism reasoning (\Cref{fig:pointwisekan2}).
Although we do not present the proof that this isomorphism is achieved by the unit of the Kan extension here,
we can formalize it within the type theory since a protype isomorphism introduces a proterm that witnesses the isomorphism by the following rule:
\[
\inferrule*
{\syn{\Gamma}\smcl\syn{\Delta}\vdash \syn{\Upsilon}:\syn{\alpha}\ccong\syn{\beta}}
{\syn{\Gamma}\smcl\syn{\Delta}\mid \syn{a}:\syn{\alpha} \vdash \syn{\Upsilon}\{\syn{a}\}:\syn{\beta}}
\]

\subsection{Related Work}

The most closely related work to \ac{FVDblTT} is VETT by New and Licata \cite{newFormalLogicFormal2023}. 
Along with the desiderata, we compare the differences between the two type theories.
Regarding (i), their type theory is designed to have the adjunction between the category of hyperdoctrines of virtual equipments
and that of its syntax,
which originates from the polymorphic feature of VETT.
The type theory has different type-theoretic entities corresponding to the hierarchy of abstractness.
It has \textit{categories}, \textit{sets}, and meta-level entities called \textit{types}, all with equational theory.
The distinction between categories in VETT and types in \ac{FVDblTT} is that the former has the equational theory as elements
of a meta-level type ``Cat,'' while the latter does not.
Although this is advantageous when different layers of category theories are in question,
it possibly obfuscates the overall type theory as a language for formal category theory.
In contrast, \ac{FVDblTT} formalizes a single layer of category theory, namely one virtual double category,
and the type theory is designed correspondingly to its components.
It also gives rise to the syntax-semantics duality between the category of 
cartesian fibrational virtual double categories and the category of its syntax,
which substantiates the type theory as an internal language of those virtual double categories. 

Regarding (ii), VETT has more constructors for types and terms than \ac{FVDblTT} in its core.
On the other hand, we focus on minimal type theory to start with
and introduce additional constructors as needed.
This is because we aim to have a type theory that reflects results in formal category theory,
which is still under development.
For instance, when we introduce the path protype to \ac{FVDblTT},
it seems plausible that it is compatible with the default finite products in the type theory,
as in \Cref{sec:appendix1}, which is supported by a category-theoretic observation 
in \Cref{sec:composition} but cannot be found in VETT.

Regarding (iii), the capability to reason about isomorphisms is a novel feature of our type theory that is not found in VETT.
It facilitates reasoning in a category theory, as explained above.

There have been other attempts to obtain a formal language for category theory.
A calculus for profunctors is presented in \cite{loregianCoEndCalculus2021} 
on the semantical level, which is followed by its type-theoretic treatment in \cite{LarLorVel24}.
Its usage is quite similar to that of \ac{FVDblTT}, but they have different focuses.
Although the calculus is similar to \ac{FVDblTT} in that it deals with profunctors and some constructors for them,
the semantics uses ordinary categories, functors, and profunctors,
while general categorical structures as its semantic environment are not given,
still less its syntax-semantics duality.
On the other hand, the coend of an endoprofunctor $\alpha(-,\bullet)$,
which cannot be handled it using \ac{FVDblTT} at the moment,
is in the scope of their calculus.
It would be interesting to know the general categorical setting where the calculus 
can be interpreted,
and it is worth investigating whether the calculus can be integrated into \ac{FVDblTT}.

%% file: newsubfiles/introduction/introoutline.tex
\Cref{section:vdc} summarizes the terminology and notation used in this paper.
\Cref{sec:fvdtt} introduces the syntax and the equational theory of \ac{FVDblTT}
and its semantics in virtual double categories.
\Cref{sec:additional} explains the type theory's possible extensions with additional constructors
and how they work in the semantics with examples.
In \Cref{sec:synsemadj}, we present the main result of this paper, 
the biadjunction between the 2-category of virtual double categories and the 2-category of \ac{FVDblTT} specifications.
This result directly implies the soundness and completeness of the type theory.

%% file: newsubfiles/prelim/vdc.tex
In this section, we briefly recall the definition of a virtual double category
and introduce the notion of a cartesian fibrational virtual double category. 

\begin{definition}[{\cite[Definition 2.1]{cruttwellUnifiedFrameworkGeneralized2010}}]
    A \emph{\acl{VDC}} (\ac{VDC}) $\dbl{D}$ is a structure consisting of the following data.
    \begin{itemize}
        \item A category $\dbl{D}_{\textbf{t}}$. Its objects are simply called \emph{objects},
        and its arrows are called \emph{tight arrows}, which are depicted vertically in this paper.
        \item A class of \emph{loose arrows} $\dbl{D}(I,J)_0$ for each pair of objects $I,J \in \dbl{D}_{\textbf{t}}$. These
        arrows are depicted horizontally with slashes as $\alpha\colon I \sto J$.
        \item A class of \emph{(virtual) cells} 
        \begin{equation}
            \label{eq:cell1}
            \begin{tikzcd}[virtual]
                I_0
                \ar[d, "s"']
                \sar[r, "\alpha_1"]
                \ar[rrrd, phantom, "\mu"]
                & I_1
                \sar[r]
                & \cdots
                \sar[r, "\alpha_n"]
                & I_n
                \ar[d, "t"] \\
                J_0
                \sar[rrr, "\beta"']
                & & & J_1
            \end{tikzcd}
        \end{equation}
        for each dataset consisting of $n\geq 0$, objects $I_0, \ldots, I_n, J_0, J_1 \in \dbl{D}_{\textbf{t}}$, 
        tight arrows $s\colon I_0 \to J_0$ and $t\colon I_n \to J_1$, and loose arrows $\alpha_1, \ldots, \alpha_n, \beta$.
        We will write the finite sequence of loose arrows as $\ol\alpha = \alpha_1;\dots;\alpha_n$.
        When $s$ and $t$ are identities, we call the cell a \emph{globular cell} 
        and let $\mu\colon \ol\alpha \Rightarrow \beta$ denote the cell.
        The class of globular cells $\ol\alpha\Rightarrow\beta$ would
        be denoted by $\dbl{D}(\ol{I})(\ol\alpha,\beta)$ in which $\ol{I} = I_0;\dots;I_n$.
        \item A composition operation on cells that assigns to each dataset of cells
        \[
            \begin{tikzcd}[column sep=4em,virtual]
                I_{1,0}
                \ar[d, "s_0"']
                \sard[r, "\ol\alpha_1"]
                \ar[dr, phantom, "\mu_1"]
                & I_{1,m_1}
                \ar[d, "s_1"']
                \sard[r, "\ol\alpha_{2}"]
                \ar[dr, phantom, "\mu_2"]
                & I_{2,m_2}
                \ar[d, "s_2"']
                \sard[r]
                & \cdots
                \sard[r, "\ol\alpha_{n}"]
                \ar[dr, phantom, "\mu_n"]
                & I_{n,m_n}
                \ar[d, "s_n"] \\
                J_{0}
                \ar[d, "t_0"']
                \sar[r, "\beta_1"] 
                \ar[drrrr, phantom, "\nu"]
                & J_{1}
                \sar[r, "\beta_2"]
                & J_{2}
                \sar[r]
                & \cdots
                \sar[r, "\beta_n"]  
                & J_{n}
                \ar[d, "t_1"] \\
                K_{0}
                \sar[rrrr, "\gamma"']
                & & & & K_{1}
            \end{tikzcd}
        \]
        a cell 
        \[
            \begin{tikzcd}[column sep=4em,virtual]
                I_{1,0}
                \ar[d, "s_0"']
                \sard[r, "\ol\alpha_1"]
                \ar[ddrrrr, phantom, "{\nu\{\mu_1\smcl\dots\smcl\mu_n\}}"]
                & I_{1,m_1}
                \sard[r, "\ol\alpha_{2}"]
                & I_{2,m_2}
                \sard[r]
                & \cdots
                \sard[r, "\ol\alpha_{n}"]
                & I_{n,m_n}
                \ar[d, "s_n"] \\
                J_{0}
                \ar[d, "t_0"']
                &&&& J_{n}
                \ar[d, "t_1"] \\
                K_{0}
                \sar[rrrr, "\gamma"']
                & & & & K_{1}
            \end{tikzcd},
        \]
        where the dashed line represents finite sequences of loose arrows
        for which associativity axioms hold.
        We will write the finite sequence of cells as $\ol\mu = \mu_1;\dots;\mu_n$.
        \item An identity cell for each loose arrow $\alpha\colon I \sto J$
        \[
            \begin{tikzcd}[virtual]
                I
                \ar[d, "\id_I"']
                \sar[r, "\alpha"]
                \ar[dr, phantom, "\id_\alpha"]
                & J
                \ar[d, "\id_J"] \\
                I
                \sar[r, "\alpha"']
                & J
            \end{tikzcd},
        \]
        for which identity laws axioms hold.
        (Henceforth, we will just write $=$ for the identity tight arrows.)
    \end{itemize}
\end{definition}

We say two object $I,J$ in a virtual double category are isomorphic 
if they are isomorphic in the underlying tight category $\dbl{D}_{\textbf{t}}$, and write $I\cong J$.
For any objects $I,J$ in a virtual double category, we write $\dbl{D}(I,J)$ for the
category whose objects are loose arrows $\alpha\colon I\sto J$ and whose arrows are cells $\mu\colon\alpha\Rightarrow\beta$. 
A cell is called an \emph{isomorphism cell} if it is invertible in this category.
More generally, 
we say two loose arrows $\alpha,\beta$ are isomorphic if there exist two cells 
\[             
\begin{tikzcd}[virtual]
    I
    \ar[d, "s"']
    \sar[r, "\alpha"]
    \ar[dr, phantom, "\mu"]
    & J
    \ar[d, "t"] \\
    K
    \sar[r, "\beta"']
    & L
\end{tikzcd}
\quad\text{and}\quad
\begin{tikzcd}[virtual]
    K
    \ar[d, "s'"']
    \sar[r, "\beta"]
    \ar[dr, phantom, "\nu"]
    & L
    \ar[d, "t'"] \\
    I
    \sar[r, "\alpha"']
    & J
\end{tikzcd}
\]
such that $\mu\{\nu\} = \id_{\beta}$ and $\nu\{\mu\} = \id_{\alpha}$,
and call the cells $\mu$ and $\nu$ \emph{isomorphism cells}.
It is always the case that $I\cong K$ and $J\cong L$ through the tight arrows $s,t,s',t'$.

\begin{example}
\label{example:dblcat}
A double category is a virtual double category where every sequence of loose arrows is composable.
In this case, a cell \cref{eq:cell1} is a cell whose top loose arrow is the composite of the loose arrows $\alpha_1,\dots,\alpha_n$.
\end{example}

\begin{definition}[{\cite[Definition 3.1]{cruttwellUnifiedFrameworkGeneralized2010}}]
    A \emph{virtual double functor} $F\colon \dbl{D} \to \dbl{E}$ between virtual double categories $\dbl{D}$ and $\dbl{E}$ 
    consists of the following data and conditions:
    \begin{itemize}
        \item A functor $F_{\textbf{t}}\colon \dbl{D}_{\textbf{t}} \to \dbl{E}_{\textbf{t}}$.
        \item A family of functions $F_1\colon \dbl{D}(I,J)_0 \to \dbl{E}(F_{\textbf{t}}(I), F_{\textbf{t}}(J))_0$ for each pair of objects $I,J$ of $\dbl{D}$.
        \item A family of functions sending each cell $\mu$ of $\dbl{D}$
        on the left below
        to a cell $F_1(\mu)$ of $\dbl{E}$ on the right below:
        \begin{equation}
            \label{eq:vdfunc}
            \begin{tikzcd}[virtual]
                I_0
                \ar[d, "s_0"']
                \sar[r, "\alpha_1"]
                \ar[rrrd, phantom, "\mu"]
                & I_1
                \sar[r]
                & \cdots
                \sar[r, "\alpha_n"]
                & I_n
                \ar[d, "s_1"] \\
                J_0
                \sar[rrr, "\beta"']
                & & & J_1
            \end{tikzcd}
            \quad \mapsto \quad
            \begin{tikzcd}[virtual]
                F_{\textbf{t}}(I_0)
                \ar[d, "F_{\textbf{t}}(s_0)"']
                \sar[r, "F_1(\alpha_1)"]
                \ar[rrrd, phantom, "F_1(\mu)"]
                & F_{\textbf{t}}(I_1)
                \sar[r]
                & \cdots
                \sar[r, "F_1(\alpha_n)"]
                & F_{\textbf{t}}(I_n)
                \ar[d, "F_{\textbf{t}}(s_1)"] \\
                F_{\textbf{t}}(J_0)
                \sar[rrr, "F_1(\beta)"']
                & & & F_{\textbf{t}}(J_1)
            \end{tikzcd}
            .
        \end{equation}
        \item The identity cells are preserved.
        \item Composition of cells is preserved.
    \end{itemize}
    As usual, we will often omit the subscripts of the functor and functions $F_{\textbf{t}}$ and $F_1$.

    A \emph{vertical transformation} $\theta\colon F \to G$ between virtual double functors $F, G\colon \dbl{D} \to \dbl{E}$
    consists of the following data and conditions:
    \begin{itemize}
        \item A natural transformation $\theta_0\colon F_{\textbf{t}} \to G_{\textbf{t}}$.
        \item A cell $\theta_{1,\alpha}$ for each loose arrow $\alpha\colon I \sto J$ of $\dbl{D}$:
        \[
            \begin{tikzcd}[virtual]
                FI 
                \ar[d, "\theta_{0,I}"']
                \sar[r, "F\alpha"]
                \ar[dr, phantom, "\theta_{1,\alpha}"]
                & FJ
                \ar[d, "\theta_{0,J}"] \\
                GI
                \sar[r, "G\alpha"']
                & GJ
            \end{tikzcd}
        \]
        \item The naturality condition for cells:
        \[
            \begin{tikzcd}[virtual]
                FI_0
                \ar[d, "Fs_0"']
                \sard[r, "F\ol\alpha"]
                \ar[dr, phantom, "F\mu"]
                & FI_n
                \ar[d, "Fs_n"] \\
                FJ_0
                \ar[d, "\theta_{J_0}"']
                \sar[r, "F\beta"']
                \ar[dr, phantom, "\theta_\beta",yshift=-0.3em]
                & FJ_1
                \ar[d, "\theta_{J_1}"] \\
                GJ_0
                \sar[r, "G\beta"']
                & GJ_1
            \end{tikzcd}
            \quad = \quad
            \begin{tikzcd}[virtual]
                FI_0
                \ar[d, "\theta_{I_0}"']
                \sard[r, "F\ol\alpha"]
                \ar[dr, phantom, "\theta_{\ol\alpha}"]
                & FI_n
                \ar[d, "\theta_{I_n}"] \\
                GI_0
                \ar[d, "Gs_0"']
                \sard[r, "G\ol\alpha"']
                \ar[dr, phantom, "G\mu",yshift=-0.3em]
                & GI_n
                \ar[d, "Gs_n"] \\
                GJ_0
                \sar[r, "G\beta"']
                & GJ_1
            \end{tikzcd}.
        \]
    \end{itemize}
    $\VDbl$ is the 2-category of virtual double categories, virtual double functors, and vertical transformations.
\end{definition}

\begin{definition}[{\cite[Definition 7.1]{cruttwellUnifiedFrameworkGeneralized2010}}]
    \label{def:fibrational}
    Let $\dbl{D}$ be a virtual double category. 
    A \emph{restriction} of a loose arrow $\alpha\colon I \sto J$ 
    along a pair of tight arrows $s\colon I' \to I$ and $t\colon J' \to J$ is the loose arrow $\alpha[s\smcl t]\colon I' \sto J'$ 
    equipped with a cell 
    \[
        \begin{tikzcd}[virtual]
            I'
            \ar[d, "s"']
            \sar[r, "{\alpha[s\smcl t]}"]
            \ar[dr, phantom, "\restc"]
            & J'
            \ar[d, "t"] \\
            I
            \sar[r, "\alpha"']
            & J
        \end{tikzcd}
    \]
    with the following universal property: any cell $\mu$ of the form on the left below factors uniquely through the cell 
    $\restc$ as on the right below.
    \[
        \begin{tikzcd}[virtual]
            K 
            \ar[d, "u"']
            \sard[r, "\ol\beta"]
            \ar[ddr, phantom, "\mu"]
            & 
            L 
            \ar[d, "v"] \\
            I'
            \ar[d, "s"']
            &
            J'
            \ar[d, "t"] \\
            I
            \sar[r, "\alpha"']
            & J
        \end{tikzcd}
        \quad = \quad
        \begin{tikzcd}[virtual]
            K 
            \ar[d, "u"']
            \sard[r, "\ol\beta"]
            \ar[dr, phantom, "\widehat{\mu}"]
            &
            L
            \ar[d, "v"] \\
            I'
            \ar[d, "s"']
            \sar[r, "{\alpha[s\smcl t]}"']
            \ar[dr, phantom, "{\restc}", yshift=-0.3em]
            & J'
            \ar[d, "t"] \\
            I
            \sar[r, "\alpha"']
            & J
        \end{tikzcd}
    \]
    If the restrictions exist for all triples $(\alpha, s, t)$, then we say that $\dbl{D}$ is a \emph{\ac{FVDC}}\footnote{
    A more standard adjective to describe this property in the literature is \emph{fibrant} \cite{aleiferiCartesianDoubleCategories2018},
    but we prefer the term \emph{fibrational} because it has nothing to do with any model structure.}.

    A \emph{fibrational virtual double functor} $F\colon \dbl{D} \to \dbl{E}$ between fibrational virtual double categories $\dbl{D}$ and $\dbl{E}$
    is a virtual double functor that preserves restrictions.
    $\FibVDbl$ is the 2-category of fibrational virtual double categories, fibrational virtual double functors, and vertical transformations.
\end{definition}

\begin{lemma}
    \label{lemma:fibvdblequiv}
    A virtual double functor $F\colon \dbl{D} \to \dbl{E}$ is an equivalence in $\VDbl$ if and only if
    \begin{enumerate}
        \labeleditem  the functor $F_{\textbf{t}}\colon \dbl{D}_{\textbf{t}} \to \dbl{E}_{\textbf{t}}$ for $F$ is an equivalence of categories, \label{lem:fibvdblequiv1}
        \labeleditem \label{lem:fibvdblequiv2} for any loose arrow $\alpha\colon I \sto J$ in $\dbl{E}$, there exists
        a loose arrow $\beta\colon I' \sto J'$ in $\dbl{D}$ and an isomorphism cell $\mu$ as below:
        \[
            \begin{tikzcd}[virtual]
                FI'
                \ar[d, "\rotatebox{90}{$\cong$}"']
                \sar[r, "F\beta"]
                \ar[dr, phantom, "\mu\ \rotatebox{90}{$\cong$}"]    
                & FJ'
                \ar[d, "\rotatebox{90}{$\cong$}"] \\
                I
                \sar[r, "\alpha"']
                & J
            \end{tikzcd},
            \quad\text{and}
        \]
        \labeleditem \label{lem:fibvdblequiv3} for any quadruple $(s,t,\ol\alpha,\beta)$,
        the function $F$ on the cells \cref{eq:vdfunc} is a bijection.
    \end{enumerate}
    A fibrational virtual double functor $F\colon \dbl{D} \to \dbl{E}$ is an equivalence in $\FibVDbl$ if and only if
    \Cref{lem:fibvdblequiv1},\Cref{lem:fibvdblequiv2}, and the special case of \Cref{lem:fibvdblequiv3} where $s$ and $t$ are identities are satisfied.
\end{lemma}

\begin{proof}
    If we are given an inverse $G$ of $F$, then $G_{\textbf{t}}$ is the inverse of $F_{\textbf{t}}$,
    and the isomorphism $FG\Rightarrow\Idf$ gives the isomorphism cells $\mu$ above.
    The inverse of functions $F$ in \Cref{eq:vdfunc} is given by sending a cell $\nu$ on the right to $G_1(\nu)$
    and composing with the isomorphism cells obtained from the isomorphism $GF\Rightarrow\Idf$.

    Conversely, given the conditions, we can construct an inverse $G$ of $F$.
    The vertical part of $G$ is given by an inverse of $F_{\textbf{t}}$.
    Then, for each loose arrow $\alpha\colon I \sto J$ in $\dbl{E}$, 
    we can show that a loose arrow $\beta\colon GI \sto GJ$ in $\dbl{D}$ is isomorphic to $\alpha$ by the second condition.
    The bijection in \Cref{lem:fibvdblequiv3} determines how to send a cell in $\dbl{E}$ to a cell in $\dbl{D}$.
    The functoriality of $G$ follows from the one-to-one correspondence between cells in $\dbl{D}$ and $\dbl{E}$ in \Cref{lem:fibvdblequiv3}.

    For the fibrational case, we only need to check that \Cref{lem:fibvdblequiv3}, in general, follows from the cases where $s$ and $t$ are identities.
    However, it is straightforward by the universal property of the restrictions.
    It follows that the inverse is fibrational from the fact that any equivalence preserves restrictions.
\end{proof}

For later use, we define restrictions of cells along a sequence of tight arrows.
\begin{definition}
    \label{def:restrictioncell}
    Let $\dbl{D}$ be an \ac{FVDC}.
    Given a globular cell $\mu$ as in \Cref{eq:cell1} with $s$ and $t$ identities
    and a sequence of tight arrows $f_i\colon K_i \to I_i$ for $0\leq i\leq n$,
    we define the \emph{restriction} of $\mu$ along the sequence $\ol{f} = f_0\smcl\dots\smcl f_n$ as the globular cell 
    $\mu[\ol{f}]$ in the diagram below 
    defined as the unique cell that makes the following equation hold.
    \[
        \begin{tikzcd}[column sep=8ex,virtual]
            K_0
            \sar[r, "{\alpha_1[f_0\smcl f_1]}"{yshift=1ex}]
            \ar[d, "f_0"']
            \ar[phantom,rd, "\restc" description]
            & K_1
            \sar[r, phantom, "\cdots"]
            \ar[d, "f_1"']
            \ar[phantom,rd, "\cdots" description]
            & 
            K_{n-1}
            \sar[r, "{\alpha_n[f_{n-1}\smcl f_n]}"{yshift=1ex}]
            \ar[d, "f_{n-1}"']
            \ar[phantom,rd, "\restc" description]
            & K_n
            \ar[d, "f_n"] 
            \\
            I_0
            \sar[r, "\alpha_1"']
            \ar[d, equal]
            \ar[phantom,rrrd, "\mu" description]
            & 
            I_1
            \sar[r, phantom, "\cdots"]
            &
            I_{n-1}
            \sar[r, "\alpha_n"']
            &
            I_n
            \ar[d, equal]
            \\
            I_0
            \sar[rrr, "\beta"']
            &&&
            I_n
        \end{tikzcd}
        =
        \begin{tikzcd}[column sep=6ex,virtual]
            K_0
            \sar[r, "{\alpha_1[f_0\smcl f_1]}"{yshift=1ex}]
            \ar[d,equal]
            \ar[drrr, phantom, "{\mu[\ol{f}]}" description]
            & K_1
            \sar[r, phantom, "\cdots"]
            &
            K_{n-1}
            \sar[r, "{\alpha_n[f_{n-1}\smcl f_n]}"{yshift=1ex}]
            &
            K_n
            \ar[d, equal]
            \\
            K_0
            \ar[d, "f_0"']
            \sar[rrr, "{\beta[f_0\smcl f_n]}"']
            \ar[drrr, phantom, "\restc"{yshift=-1ex,description}]
            &&&
            K_n
            \ar[d, "f_n"]
            \\
            I_0
            \sar[rrr, "\beta"']
            &&&
            I_n
        \end{tikzcd}
    \]
\end{definition}

Next, we introduce the notion of a cartesian fibrational virtual double category.

\begin{definition}
    A \emph{cartesian object} in a 2-category $\bi{B}$ with finite products $\monunit,\otimes$ is an object $x$ of $\bi{B}$ such that
    the canonical 1-cells $!\colon x\to \monunit$ and $\Delta\colon x\to x\otimes x$ have right adjoints $1\colon\monunit\to x$ and $\times\colon x\otimes x\to x$, respectively.
    A \emph{cartesian 1-cell} (or \emph{cartesian arrow}) in $\bi{B}$ is a 1-cell $f\colon x\to y$ between cartesian objects $x$ and $y$ of $\bi{B}$ such that
    the canonical 2-cells obtained by the mate construction $\times\circ(f\otimes f)\Rightarrow f\circ \times$ and $f\circ 1\Rightarrow 1$ are invertible. 

    For a 2-category $\bi{B}$ with finite products, 
    we write $\bi{B}_{\bi{cart}}$ for the 2-category of cartesian objects, cartesian 1-cells, and arbitrary 2-cells in $\bi{B}$.
\end{definition}

\begin{lemma}
    \label{lemma:Cartesianequiv}
    Let $\bi{B}$ be a 2-category with finite products. 
    A 1-cell $f\colon x\to y$ in $\bi{B}_{\bi{cart}}$ is an equivalence in $\bi{B}_{\bi{cart}}$ 
    if and only if the underlying 1-cell of $f$ is an equivalence in $\bi{B}$.
\end{lemma}
\begin{proof}
    The only if part is clear since we have the forgetful 2-functor $\bi{B}_{\bi{cart}}\to\bi{B}$.
    For the if part, 
    take the right adjoint $g$ of the underlying 1-cell of $f$ as its inverse.
    Taking the right adjoint of both sides of the isomorphism 2-cells $!\circ f\cong !$ and $(f\otimes f)\circ\Delta\cong\Delta\cong f$,
    we obtain the isomorphism 2-cells $1\fatsemi g\cong 1$ and $\times\circ (g\otimes g)\cong g\circ\times$.
    This shows that $g$ gives a cartesian morphism from $y$ to $x$,
    and $g$ is indeed the inverse of $f$ in $\bi{B}_{\bi{cart}}$.
\end{proof}

\begin{proposition}
    \label{prop:FibVDblCart}
    An \ac{FVDC} $\dbl{D}$ is cartesian if and only if the following conditions are satisfied:
    \begin{enumerate}
        \item $\dbl{D}_{\textbf{t}}$ has finite products;
        \item $\dbl{D}$ locally has finite products, that is, for each $I,J\in\dbl{D}_{\textbf{t}}$,
        \begin{enumerate}
        \item for any loose arrows
        $\alpha,\beta\colon I\sto J$ in $\dbl{D}$,
        there exists a loose arrow $\alpha\land\beta\colon I\sto J$ and cells 
        $\pi_0\colon\alpha\land\beta\Rightarrow\alpha$, $\pi_1\colon\alpha\land\beta\Rightarrow\beta$
        such that for any finite sequence of loose arrows $\ol\gamma$ where $\gamma_i\colon I_{i-1}\sto I_i$ for $1\leq i\leq n$ where $I_0=I$ and $I_n=J$,
        the function 
        \[ 
            \dbl{D}(\ol{I})(\ol\gamma,\alpha\land\beta) \to
            \dbl{D}(\ol{I})(\ol\gamma,\alpha)\times\dbl{D}(\ol{I})(\ol\gamma,\beta)\quad;
            \quad \mu\mapsto (\pi_0\circ\mu,\pi_1\circ\mu)
        \]
        is a bijection, and
        \item there exists a loose arrow $\top\colon I\sto J$ such that $\dbl{D}(\ol{I})(\ol\gamma,\top)_0$ is
        a singleton for any finite sequence of loose arrows $\ol\gamma$;
        \end{enumerate}
        \item the local finite products are preserved by restrictions.
    \end{enumerate} 
    A morphism between cartesian \acp{FVDC} is a cartesian morphism if and only if the underlying tight functor preserves finite products
    and the morphism preserves local finite products.
\end{proposition}

\begin{proof}[Proof sketch]
The proof is similar to that of \cite[Prop 4.12]{aleiferiCartesianDoubleCategories2018}.
First, suppose that $\dbl{D}$ is cartesian.
Let $\Delta_I\colon I\to I\times I$ be the diagonal of $I$ and $!_I\colon I\to 1$ be the unique arrow to the terminal object.
If $\dbl{D}$ is cartesian, then $\alpha\land\beta$ and $\top$ in $\dbl{D}(I,J)$ are given by $(\alpha\times\beta)[\Delta_I\smcl\Delta_J]$
and $\Id_1(!_I,!_J)$, which brings the finite products in $\dbl{D}(I,J)$.
The local finite products are preserved by restrictions since, by the universal property of the restrictions,
we have 
\[
    (\alpha\times\beta)[\Delta_I\smcl\Delta_J][s\smcl t] \cong (\alpha\times\beta)[(s\times s)\smcl(t\times t)][\Delta_{I'}\smcl\Delta_{J'}] 
    \cong (\alpha[s\smcl t]\times\beta[s\smcl t])[\Delta_{I'}\smcl\Delta_{J'}],
\]
and similarly for $\top$.
Conversely, if $\dbl{D}$ locally has finite products, then assigning 
\[
    \alpha\times\beta\coloneqq\alpha[\pi_I\smcl\pi_J]\land\beta[\pi_K\smcl\pi_L]\colon I\times K\sto J\times L 
\]
to each pair $\alpha\colon I\sto J,\,\beta\colon K\sto L$ and a cell $\mu\times\nu$ naturally obtained from the universal property of the restrictions 
induces the functor $\times\colon\dbl{D}\times\dbl{D}\to\dbl{D}$ right adjoint to the diagonal functor,
and the functor $1\colon \dbl{1}\to\dbl{D}$ obtained by the terminal object in $\dbl{D}_{\textbf{t}}$ is the right adjoint of $!$.
The second statement follows from the construction of the equivalence above.
\end{proof}

\begin{remark}
    The third condition in \cref{prop:FibVDblCart} is necessary for \ac{FVDC} but not for equipments 
    as in \cite{aleiferiCartesianDoubleCategories2018}
    since the latter has extensions.
\end{remark}

\begin{example}
    \label{example:prof}
    One of the motivations for our type theory is to formalize category theory in formal language.
    The following examples of virtual double categories
    will provide a multitude of category theories that can be formalized in our type theory.
    \begin{enumerate}
        \item The double category $\Prof$ consists of small categories, functors, and profunctors as objects, tight arrows, and loose arrows, respectively.
        When we consider not necessarily small categories, however, 
        we do not have a composition of profunctors in general. 
        Nevertheless, we can still define a virtual double category $\PROF$ of categories, functors, and profunctors.
        It is a \ac{CFVDC}.
        \item Similarly, we can define \acp{CFVDC} $\one{V}\mhyphen\Prof$ and 
        $\one{V}\mhyphen\PROF$ of $\one{V}$-enriched categories, functors, and profunctors,
        where $\one{V}$ is a symmetric monoidal category.
        \item We can also define virtual double categories $\Prof(\one{S})$ of
        internal categories, functors, and profunctors in categories $\one{S}$ with finite limits.
    \end{enumerate}
\end{example}

\begin{example}
        \label{example:rel}
        Predicate logic deals with functions and relations between sets.
        We can see these two as tight and loose arrows, respectively,
        although we limit ourselves to the case where relations are binary.
        \begin{enumerate}
            \item We have a double category called $\Rel{}$ consisting of sets, functions, and relations \cite{lambertDoubleCategoriesRelations2022},
                where there is at most one cell for each frame and a cell \cref{eq:cell1} exists whenever 
                \[ 
                    \alpha_1(x_0,x_1)\land\dots\land\alpha_n(x_{n-1},x_n)\Rightarrow\beta(s_0(x_0),s_1(x_n))
                    \quad (\forall x_i\in I_i, \quad 0\le i\le n).
                \]
                This can be generalized to a regular category and further to a category with finite limits equipped with a stable factorization system 
                as in \cite{hoshinoDoubleCategoriesRelations2023}.
            \item Given a category $\one{S}$ with finite limits, a span in $\one{S}$ is a pair of arrows $f\colon Z\to X$ and $g\colon Z\to Y$.
                We can define a double category called $\Span(\one{S})$ consisting of objects, arrows, and spans in $\one{S}$ \cite{aleiferiCartesianDoubleCategories2018}.
                It can be seen as a proof-relevant version of $\Rel{}$.
        \end{enumerate}
\end{example}

%% file: newsubfiles/typetheory/typetheory.tex
This section will present a type theory that serves as an internal language for \acp{CFVDC}, which we call \emph{\acl{FVDblTT}} (\ac{FVDblTT}).
We will first introduce the type theory \ac{FVDblTT} in \cref{sec:syntax}.
Then, we will present how to interpret the type theory in a \ac{CFVDC} in \cref{sec:semantics}.

%% file: newsubfiles/typetheory/syntax.tex
The syntax of \ac{FVDblTT} is given by the following grammar.

\begin{figure}[h]
    {
    \small
    \begin{align*}
        \text{Type} &\quad \syn{I}\ \textsf{type} \\
        \text{Context} &\quad \syn{\Gamma}\ \textsf{ctx} \\
        \text{Term} &\quad \syn{\Gamma} \vdash \syn{s}:\syn{I}  \\
        \text{Term Substitution} &\quad \syn{\Gamma} \vdash \syn{S}\,/\,\syn{\Delta}  \\
        \text{Protype} &\quad \syn{\Gamma}\smcl\syn{\Delta} \vdash \syn{\alpha} \ \textsf{protype}  \\
        \text{Procontext} &\quad \syn{\Gamma}_0\smcl\dots\smcl\syn{\Gamma}_n \mid \syn{A} \ \textsf{proctx} \\
        \text{Proterm} &\quad \syn{\Gamma}_0\smcl\dots\smcl\syn{\Gamma}_n \mid \syn{A} \vdash \syn{\mu}:\syn{\beta}  \\
        \text{Term Equality} &\quad \syn{\Gamma} \vdash \syn{s} \equiv \syn{t}:\syn{I} &\\
        \text{Protype Equality} &\quad \syn{\Gamma}\smcl\syn{\Delta} \vdash \syn{\alpha} \equiv \syn{\beta} &\\
        \text{Proterm Equality} &\quad \syn{\Gamma}_0\smcl\dots\smcl\syn{\Gamma}_n \mid 
        \syn{A} \vdash \syn{\mu} \equiv \syn{\nu}:\syn{\beta} 
    \end{align*}}
\caption{Judgments in \ac{FVDblTT}}
\label{fig:judgments}
\end{figure}

Types, contexts, terms, and term substitutions are the same 
as those in the algebraic theory as in \cite{Crole_1994,jacobsCategoricalLogicType1999a}.
This fragment of the type theory serves as the theory of categories and functors.
As usual, substitution of terms for variables in terms is defined by induction on the structure of terms.

Protypes and proterms are particular to this type theory and 
encode the loose arrows and cells in an \ac{CFVDC}.
The prefix \emph{pro-} stands for ``pro''positions and ``pro''functors.
A protype $\syn{\alpha}$ depends on two contexts, $\syn{\Gamma}$ and $\syn{\Delta}$,
which will be interpreted as the source and the target of a loose arrow representing the protype.
We call the pair $(\syn{\Gamma},\syn{\Delta})$ the \emph{two-sided context} of the protype
and write $\syn{\Gamma}\smcl\syn{\Delta}$ for it.
In the type theory,
we distinguish semicolons ``$\smcl$'' from the ordinary concatenating symbol commas ``$,$'' by restricting using the former to concatenate items in the horizontal direction 
in a diagram in a \ac{VDC}.
Since the source and the target of a loose arrow can not be exchanged in any sense in a general \ac{VDC},
we need to respect the order when we use the semicolons.
Accordingly, a procontext $\syn{a}_1:\syn{\alpha}_1\smcl\dots\smcl\syn{a}_n:\syn{\alpha}_n$
with \textit{provariables} $\syn{a}_i$'s, which is formally defined as a finite sequence of protypes, is only well-typed
so that the second (target) context of a protype is the first (source) context of the subsequent protype,
and hence a procontext depends on a sequence of contexts.
As a particular case, we have the empty procontext $\cdot$ depending on a single context $\syn{\Gamma}$.
Another item is proterms $\syn{A}\vdash\syn{\mu}:\syn{\beta}$
where $\syn{A}$ is a procontext and $\syn{\beta}$ is a protype,
which are interpreted as globular cells in a \ac{VDC} whose domains and codomains are the interpretation of $\syn{A}$ and $\syn{\beta}$, respectively.

The type theory also has the equality judgments for terms, protypes, and proterms.
We incorporate the ordinary algebraic theory of terms with the equality judgments for terms,
and we also have the equality judgments for proterms to capture the equality of cells in a \ac{VDC}.
The rules for equality judgments, or the equational theory of the type theory,
are based on the basic axioms of reflexivity, symmetry, transitivity, and replacement
with respect to the substitution we will define later.
The equational theory for protypes is designed only to reflect the equational theory of terms
by the replacement rule for substitution of terms,
and we do not have any other rules for protypes
except for the basic axioms.
This is because, in formal category theory, we are mainly interested in isomorphisms of loose arrows,
which we will incorporate in the type theory as the protype isomorphisms later.

\begin{figure}[htbp]
    \begin{mathpar}
        \inferrule*
        {\syn{I}\ \ \textsf{type}\\ \syn{J} \ \textsf{type}}
        {\syn{I}\times \syn{J} \  \textsf{type}}
        \and        
        \inferrule*
        { }
        {\syn{1} \ \textsf{type}}
        \and 
        \inferrule*
        { }
        {\cdot \ \textsf{ctx}}
        \and
        \inferrule*
        {\syn{\Gamma} \ \textsf{ctx} \\ \syn{I} \ \textsf{type}}
        {\syn{\Gamma},\syn{x}:\syn{I} \ \textsf{ctx}}
        \and
        \inferrule*
        {\ }
        {\syn{\Gamma},\syn{x}:\syn{I},\syn{\Delta}\vdash \syn{x}:\syn{I}}
        \and
        \inferrule*
        {\syn{I}\ \textsf{type} \\ \syn{J}\ \textsf{type} \\ \syn{\Gamma}\vdash \syn{s}:\syn{I} \\ \syn{\Gamma}\vdash \syn{t}:\syn{J}}
        {\syn{\Gamma}\vdash \langle \syn{s},\syn{t} \rangle:\syn{I}\times \syn{J}}
        \and
        \inferrule*
        {\syn{\Gamma}\vdash \syn{t}:\syn{I}\times \syn{J}}
        {\syn{\Gamma}\vdash \pr{0}(\syn{t}):\syn{I}}
        \and
        \inferrule*
        {\syn{\Gamma}\vdash \syn{t}:\syn{I}\times \syn{J}}
        {\syn{\Gamma}\vdash \pr{1}(\syn{t}):\syn{J}}
        \and
        \inferrule*
        { }
        {\syn{\Gamma} \vdash \langle \ \rangle:\syn{1}}
        \and 
        \inferrule*
        { }
        {\syn{\Gamma} \cdot \vdash /\cdot}
        \and
        \inferrule*
        {\syn{\Gamma} \vdash \syn{S}\,/\,\syn{\Delta} \\ \syn{\Gamma} \vdash \syn{s}:\syn{I}}
        {\syn{\Gamma} \vdash \syn{S},\syn{s}\,/\,\syn{\Delta},\syn{x}:\syn{I}}
        \and
        \inferrule*
        {\syn{\Gamma} \vdash \syn{s}:\syn{I}\\ \syn{\Gamma} \vdash \syn{t}:\syn{J}}
        {\syn{\Gamma} \vdash \pr{0}(\langle \syn{s},\syn{t} \rangle)\equiv \syn{s} : \syn{I}}
        \and 
        \inferrule*
        {\syn{\Gamma} \vdash \syn{s}:\syn{I}\\ \syn{\Gamma} \vdash \syn{t}:\syn{J}}
        {\syn{\Gamma} \vdash \pr{1}(\langle \syn{s},\syn{t} \rangle)\equiv \syn{t} : \syn{J}}
        \and
        \inferrule*
        {\syn{\Gamma} \vdash \syn{s}:\syn{I}\times\syn{J} }
        {\syn{\Gamma} \vdash \langle \pr{0}(\syn{s}),\pr{1}(\syn{s}) \rangle \equiv \syn{s} : \syn{I}\times\syn{J}}
        \and
        \inferrule*
        {\syn{\Gamma} \vdash \syn{s}:\syn{1}}
        {\syn{\Gamma} \vdash \syn{s}\equiv \langle \ \rangle : \syn{1}}    
    \end{mathpar}
    \caption{The rules for types, contexts, and terms}
    \label{fig:typingterms}
\end{figure}
\begin{figure}[htbp]
    \begin{mathpar}
            \inferrule*
            {\syn{\Gamma}\smcl\syn{\Delta} \vdash \syn{\alpha} \ \textsf{protype} \\ \syn{\Gamma}\smcl\syn{\Delta} \vdash \syn{\beta} \ \textsf{protype}}
            {\syn{\Gamma}\smcl\syn{\Delta} \vdash \syn{\alpha}\land \syn{\beta} \ \textsf{protype}}
            \and
            \inferrule*
            { }
            {\syn{\Gamma}\smcl\syn{\Delta}  \vdash \top \ \textsf{protype}}
            \and 
            \inferrule*
            { }
            {\syn{\Gamma} \mid \cdot \ \textsf{proctx}}
            \and
            \inferrule*
            {\syn{\Gamma}_0\smcl\dots\smcl\syn{\Gamma}_{n}\mid \syn{A} \ \textsf{proctx} \\ \syn{\Gamma}_n\smcl\syn{\Delta} \vdash \syn{\alpha} \ \textsf{protype}}
            {\syn{\Gamma}_0\smcl\dots\smcl\syn{\Gamma}_n\smcl\syn{\Delta} \mid \syn{A},\syn{a}:\syn{\alpha} \ \textsf{proctx}}
            \and
            \inferrule*
            {\syn{\Gamma}\smcl\syn{\Delta}\vdash \syn{\alpha} \ \textsf{protype}\\
            \syn{\Gamma}'\vdash \syn{S}_0\equiv\syn{S}_1\,/\,\syn{\Gamma} \\ \syn{\Delta}'\vdash\syn{T}_0\equiv\syn{T}_1\,/\,\syn{\Delta}
            }
            {\syn{\Gamma}'\smcl\syn{\Delta}' \vdash \syn{\alpha}[\syn{S}_0/\syn{\Gamma}\smcl\syn{T}_0/\syn{\Delta}] \equiv 
            \syn{\alpha}[\syn{S}_1/\syn{\Gamma}\smcl\syn{T}_1/\syn{\Delta}]}
            \and
        \inferrule*
        {\syn{\Gamma}\smcl\syn{\Delta}\vdash \syn{\alpha} \ \textsf{protype}}
        {\syn{\Gamma}\smcl\syn{\Delta}\mid \syn{a}:\syn{\alpha} \vdash \syn{a}:\syn{\alpha}}
        \and
        \inferrule*
        {\ol{\syn{\Gamma}_i}\mid \syn{a}_{i,1}:\syn{\alpha}_{i,1}\smcl\dots\smcl\syn{a}_{i,n_i}:\syn{\alpha}_{i,n_i} \vdash \syn{\mu}_i:\syn{\beta}_i \ (i=1,\dots,m) \\
        \wt{\syn{\Gamma}} \mid\syn{b}_1:\syn{\beta}_1\smcl\dots\smcl\syn{b}_n:\syn{\beta}_n \vdash\syn{\nu}: \syn{\gamma}}
        {\syn{\ol{\Gamma}}\mid \syn{a}_{1,1}:\syn{\alpha}_{1,1}\smcl\dots\smcl\syn{a}_{m,n_m}:\syn{\alpha}_{m,n_m} 
        \vdash \syn{\nu}\{\syn{\mu}_1/\syn{b}_1:\syn{\beta}_1\smcl\dots\smcl\syn{\mu}_m/\syn{b}_m:\syn{\beta}_m\} : \syn{\gamma}}
        \and
        \inferrule*
        {\syn{\ol{\Gamma}}\mid \syn{A} \vdash \syn{\mu}:\syn{\alpha} \\ \syn{\ol{\Gamma}}\mid \syn{A} \vdash \syn{\nu}:\syn{\beta}}
        {\syn{\ol{\Gamma}}\mid \syn{A} \vdash \langle \syn{\mu},\syn{\nu} \rangle:\syn{\alpha}\land\syn{\beta}}
        \and
        \inferrule*
        {\syn{\ol{\Gamma}}\mid \syn{A} \vdash \syn{\mu}:\syn{\alpha}\land\syn{\beta}}
        {\syn{\ol{\Gamma}}\mid \syn{A} \vdash \syn{\pi}_0\{\syn{\mu}\}:\syn{\alpha}}
        \and
        \inferrule*
        {\syn{\ol{\Gamma}}\mid \syn{A} \vdash \syn{\mu}:\syn{\alpha}\land\syn{\beta}}
        {\syn{\ol{\Gamma}}\mid \syn{A} \vdash \syn{\pi}_1\{\syn{\mu}\}:\syn{\beta}}
        \and
        \inferrule*
        { }
        {\syn{\ol{\Gamma}}\mid \syn{A} \vdash \langle \ \rangle:\top}
        \and 
        \inferrule*
        {\ol{\syn{\Gamma}}\mid \syn{A}\vdash \syn{\mu}:\syn{\alpha}\\
        \ol{\syn{\Gamma}}\mid \syn{A}\vdash \syn{\nu}:\syn{\beta}}
        {\ol{\syn{\Gamma}}\mid \syn{A}\vdash \syn{\pi}_0(\langle \syn{\mu},\syn{\nu} \rangle)\equiv \syn{\mu}:\syn{\alpha}}
        \and
        \inferrule*
        {\ol{\syn{\Gamma}}\mid \syn{A}\vdash \syn{\mu}:\syn{\alpha}\\
        \ol{\syn{\Gamma}}\mid \syn{A}\vdash \syn{\nu}:\syn{\beta}}
        {\ol{\syn{\Gamma}}\mid \syn{A}\vdash \syn{\pi}_1(\langle \syn{\mu},\syn{\nu} \rangle)\equiv \syn{\nu}:\syn{\beta}}
        \and
        \inferrule*
        {\ol{\syn{\Gamma}}\mid \syn{A}\vdash \syn{\mu}:\syn{\alpha}\land\syn{\beta}}
        {\ol{\syn{\Gamma}}\mid \syn{A}\vdash \langle \syn{\pi}_0(\syn{\mu}),\syn{\pi}_1(\syn{\mu}) \rangle\equiv \syn{\mu}:\syn{\alpha}\land\syn{\beta}}
        \and 
        \inferrule*
        {\ol{\syn{\Gamma}}\mid \syn{A}\vdash \syn{\mu}:\top}
        {\ol{\syn{\Gamma}}\mid \syn{A}\vdash \syn{\mu}\equiv \langle \ \rangle:\top}
        \and
        \inferrule*
        {\ol{\syn{\Gamma}}\mid \syn{a}_1:\syn{\alpha}_1\smcl\dots\smcl\syn{a}_n:\syn{\alpha}_n\vdash \syn{\mu}:\syn{\beta} \\ 
        \syn{\Gamma}_0\smcl\syn{\Gamma}_1\vdash \syn{\alpha}_1\equiv\syn{\alpha}_1' \\ 
        \dots \\
        \syn{\Gamma}_{n-1}\smcl\syn{\Gamma}_n\vdash \syn{\alpha}_n\equiv\syn{\alpha}_n'\\
        \syn{\Gamma}_0\smcl\syn{\Gamma}_n\vdash \syn{\beta}\equiv\syn{\beta}'}
        {\syn{\Gamma}_0\smcl\syn{\Gamma}_n\mid \syn{a}_1:\syn{\alpha}_1'\smcl\dots\smcl\syn{a}_n:\syn{\alpha}_n'\vdash \syn{\mu}:\syn{\beta}'} 
    \end{mathpar}
    \caption{The rules for protypes, procontexts, and proterms}
    \label{fig:proterms}
\end{figure}

\myparagraph{Signatures.}
In algebraic theories, one often starts with a signature that specifies the sorts and operations of the theory. 
We present the notion of a signature for \ac{FVDblTT} as follows.
\begin{definition}
    A \emph{signature} $\Sigma$ for \ac{FVDblTT} is a quadruple $(T_{\Sigma},F_{\Sigma},P_{\Sigma},C_{\Sigma})$ where
    \begin{itemize}
        \item $T_{\Sigma}$ is a class of \emph{category symbols},
        \item $F_{\Sigma}(\syn{\sigma},\syn{\tau})$ is a family of classes of \emph{functor symbols} for any $\syn{\sigma},\syn{\tau}\in T_{\Sigma}$, 
        \item $P_{\Sigma}(\syn{\sigma}\smcl\syn{\tau})$ is a family of classes of \emph{profunctor symbols} for any $\syn{\sigma},\syn{\tau}\in T_{\Sigma}$,
        \item $C_{\Sigma}(\syn{\rho}_1\smcl\dots\smcl\syn{\rho}_n\mid\syn{\omega})$ is a family of classes of \emph{transformation symbols} for any $\syn{\sigma}_0,\dots,\syn{\sigma}_n\in T_{\Sigma}$, $\syn{\rho}_i\in P_{\Sigma}(\syn{\sigma}_{i-1}\smcl\syn{\sigma}_i)$ for $i=1,\dots,n$, and $\syn{\omega}\in P_{\Sigma}(\syn{\sigma}_0\smcl\syn{\sigma}_n)$
        where $n\geq 0$. 
    \end{itemize}

    For simplicity, in the last item, we omit the dependency of the class of transformation symbols on $\syn{\sigma}_i$'s.
    Henceforth, $\syn{f}\colon \syn{\sigma}\syto\syn{\tau}$ denotes a functor symbol $\syn{f}\in F_{\Sigma}(\syn{\sigma},\syn{\tau})$,
    $\syn{\rho}\colon \syn{\sigma}\systo\syn{\tau}$ denotes a profunctor symbol $\syn{\rho}\in P_{\Sigma}(\syn{\sigma}\smcl\syn{\tau})$,
    and $\syn{\kappa}\colon \syn{\rho}_1\smcl\dots\smcl\syn{\rho}_n\sydto\syn{\omega}$ denotes a transformation symbol $\syn{\kappa}\in C_{\Sigma}(\syn{\rho}_1\smcl\dots\smcl\syn{\rho}_n\mid\syn{\omega})$. 
    
    A \emph{morphism of signatures} $\Phi\colon\Sigma\to\Sigma'$ is 
    a family of functions sending the symbols of each kind in $\Sigma$ to symbols of the same kind in $\Sigma'$
    so that a symbol dependent on another kind of symbol is sent to a symbol dependent on the
    image of the former symbol.
    For instance, $\syn{\rho}\colon \syn{\sigma}\systo\syn{\tau}$ is sent to a profunctor symbol 
    of the form $\Phi(\syn{\rho})\colon \Phi(\syn{\sigma})\systo\Phi(\syn{\tau})$
    where the assignment of category symbols has already been determined.
\end{definition}

A typical example of a signature is the signature defined by a \ac{CFVDC} $\dbl{D}$.

\begin{definition}
    The \emph{associated signature} of a \ac{CFVDC} $\dbl{D}$ is
    the signature $\Sigma_{\dbl{D}}$ defined by
    \begin{itemize}
        \item $T_{\dbl{D}}$ is the set of objects of $\dbl{D}$, 
        where we write $\is{I}$ for $I\in \dbl{D}$ as a category symbol,
        \item $F_{\dbl{D}}(\is{I},\is{J})$ is the set of tight arrows 
        $I\to J$ in $\dbl{D}$, 
        where we write $\is{f}$ for $f\in F_{\dbl{D}}(\is{I},\is{J})$ as a functor symbol,
        \item $P_{\dbl{D}}(\is{I}\smcl\is{J})$ is the set of loose arrows
        $\alpha\colon I\sto J$ in $\dbl{D}$, 
        where we write $\is{\alpha}$ for $\alpha\in P_{\dbl{D}}(\is{I}\smcl\is{J})$ as a profunctor symbol,
        \item $C_{\dbl{D}}(\is{\alpha_1}\smcl\dots\smcl\is{\alpha_n}\mid\is{\beta})$ is the set of cells
        $\mu\colon \alpha_1;\dots;\alpha_n\Rightarrow\beta$
        in $\dbl{D}$, where we write $\is{\mu}$ for $\mu\in C_{\dbl{D}}(\is{\alpha_1}\smcl\dots\smcl\is{\alpha_n}\mid\is{\beta})$ as a transformation symbol.
    \end{itemize}
\end{definition}

A signature $\Sigma$ is what we start derivations with in the type theory.
In terms of formal category theory,
it signifies what one regard as categories, functors, profunctors, and 
natural transformations.
The rules for the signature are given as follows.
\begin{figure}[htbp]
    \begin{mathparpagebreakable}
        \inferrule*
        {\syn{\sigma}\in T_{\Sigma}}
        {\syn{\sigma}\ \textsf{type}}
        \and
        \inferrule*
        {\syn{f}\in F_{\Sigma}(\syn{\sigma},\syn{\tau})\\
        \syn{\Gamma}\vdash \syn{s}:\syn{\sigma}}
        {\syn{\Gamma}\vdash \syn{f}(\syn{s}):\syn{\tau}}
        \and
        \inferrule*
        {\syn{\rho}\in P_{\Sigma}(\syn{\sigma},\syn{\tau})\\
        \syn{\Gamma}\vdash \syn{s}:\syn{\sigma}\\
        \syn{\Delta}\vdash \syn{t}:\syn{\tau}}
        {\syn{\Gamma}\smcl\syn{\Delta}\vdash \syn{\rho}(\syn{s}\smcl\syn{t}):\textsf{protype}}
        \and
        \inferrule*
        {\syn{\kappa}\in C_{\Sigma}(\syn{\rho}_1\smcl\dots\smcl\syn{\rho}_n\mid\syn{\omega})\\
        \syn{\Gamma}_i\vdash \syn{s}_i:\syn{\sigma}_i\quad (i=0,\dots,n)\\
        \syn{\Gamma}_{i-1}\smcl\syn{\Gamma}_i\mid\syn{A}_i\vdash \syn{\mu}_i:\syn{\rho}_i(\syn{s}_{i-1}\smcl\syn{s}_i)\quad (i=1,\dots,n)}
        {\syn{\Gamma}_0\smcl\dots\smcl\syn{\Gamma}_n\mid\syn{A}_1\smcl\dots\smcl\syn{A}_n\vdash \syn{\kappa}(\syn{s}_0\smcl\dots\smcl\syn{s}_n)\{\syn{\mu_1}\smcl\dots\smcl\syn{\mu_n}\}}
    \end{mathparpagebreakable}
    \caption{The rules for the signature}
    \label{fig:signature}
\end{figure}

\myparagraph{Substitution.}
The substitution of terms for variables in terms, protypes, and proterms is defined
inductively as follows.
{\small
        \begin{align*}
            \syn{x}_i[\syn{S}/\syn{\Delta}]
            &\equiv \syn{s}_i \\
            & \quad (i=1,\dots,n,\ \syn{S}=(\syn{s}_1,\dots,\syn{s}_n))\\
            \syn{f}(\syn{s}_1,\dots,\syn{s}_n)[\syn{S}/\syn{\Delta}]
            &\equiv \syn{f}(\syn{s}_1[\syn{S}/\syn{\Delta}],\dots,\syn{s}_n[\syn{S}/\syn{\Delta}])\\
            \langle\syn{s},\syn{t}\rangle[\syn{S}/\syn{\Delta}]
            &\equiv \langle\syn{s}[\syn{S}/\syn{\Delta}],\syn{t}[{\syn{S}/\syn{\Delta}}]\rangle\\
            (\pr{i}(\syn{t}))[\syn{S}/\syn{\Delta}]
            &\equiv\pr{i}(\syn{t}[\syn{S}/\syn{\Delta}]) \\
            \langle\ \rangle[\syn{S}/\syn{\Delta}]
            &\equiv \langle\ \rangle\\
            (\syn{\rho}(\syn{s}\smcl\syn{t}))[\syn{S}/\syn{\Delta}\smcl\syn{T}/\syn{\Theta}]
            &\equiv \syn{\rho}(\syn{s}[\syn{S}/\syn{\Delta}]\smcl\syn{t}[{\syn{T}/\syn{\Theta}}])\\
            (\syn{\alpha}\land\syn{\beta})[\syn{S}\,/\,\syn{\Delta}\smcl\syn{T}\,/\,\syn{\Theta}]
            &\equiv \syn{\alpha}[\syn{S}/\syn{\Delta}\smcl\syn{T}/\syn{\Theta}]\land\syn{\beta}[\syn{S}/\syn{\Delta}\smcl\syn{T}/\syn{\Theta}] \\ 
            \top[\syn{S}\,/\,\syn{\Delta}\smcl\syn{T}\,/\,\syn{\Theta}]
            &\equiv \top\\
            \syn{a}[{\syn{S}\,/\,\syn{\Delta}\smcl\syn{T}\,/\,\syn{\Theta}}]
            &\equiv \syn{a}\\
            \left(\syn{\kappa}(\ol{\syn{s}_{\ul{i}}})\{\ol{\syn{\mu}_{\ul{i}}}\}\right)
            [\ol{\syn{S}_{\ul{i},\ul{j}}}/\ol{\syn{\Delta}_{\ul{i},\ul{j}}}] 
            &\equiv \syn{\kappa}\left(\ol{\syn{s}_{\ul{i}}[\syn{S}_{\ul{i},n_{\ul{i}}}/\syn{\Delta}_{\ul{i},n_{\ul{i}}}]}\right)         
            \{\ol{\syn{\mu}_{\ul{i}}[\ol{\syn{S}_{\ul{i},\ul{j}}}/\ol{\syn{\Delta}_{\ul{i},\ul{j}}}]} \} \\
            \langle\syn{\mu},\syn{\mu}'\rangle[\ol{\syn{S}_{\ul{i}}}/\ol{\syn{\Delta}_{\ul{i}}}] 
            &\equiv \langle\syn{\mu}[\ol{\syn{S}_{\ul{i}}}/\ol{\syn{\Delta}_{\ul{i}}}],\syn{\mu}'[\ol{\syn{S}_{\ul{i}}}/\ol{\syn{\Delta}_{\ul{i}}}]\rangle\\
            (\syn{\pi}_i\{\syn{\mu}\})[\ol{\syn{S}_{\ul{i}}}/\ol{\syn{\Delta}_{\ul{i}}}]
            &\equiv \syn{\pi}_i\{\syn{\mu}[\ol{\syn{S}_{\ul{i}}}/\ol{\syn{\Delta}_{\ul{i}}}]\}\\
            \langle\ \rangle[\ol{\syn{S}_{\ul{i}}}/\ol{\syn{\Delta}_{\ul{i}}}]
            &\equiv \langle\ \rangle\\
        \end{align*}
}

Since the type theory has a different layer consisting of protypes and proterms,
we need to define substitution for them as well,
which we call \emph{prosubstitution} and symbolize by $\psbsm{\cdot}$
to distinguish it from the usual substitution.
The prosubstitution is defined inductively as follows.
{\small 
\begin{align*}
    \syn{a}\psbsm{\syn{\mu}/\syn{a}}
    &\equiv \syn{\mu}\\
    \left(\syn{\kappa}(\ol{\syn{s}_{\ul{i}}})\{\ol{\syn{\mu}_{\ul{i}}}\}\right)
    \psb{\ol{\syn{\nu}_{\ul{i},\ul{j}}}/\ol{\syn{b}_{\ul{i},\ul{j}}}}
    &\equiv
    \syn{\kappa}(\ol{\syn{s}_{\ul{i}}})
    \left\{\ol{\syn{\mu}_{\ul{i}}\psb{\ol{\syn{\nu}_{i,\ul{j}}}/\ol{\syn{b}_{i,\ul{j}}}}}\right\}\\
    \langle\syn{\mu},\syn{\mu}'\rangle\psb{\ol{\syn{\nu}_{\ul{i}}}/\ol{\syn{b}_{\ul{i}}}} 
    &\equiv \langle\syn{\mu}\psb{\ol{\syn{\nu}_{\ul{i}}}/\ol{\syn{b}_{\ul{i}}}},\syn{\mu}'\psb{\ol{\syn{\nu}_{\ul{i}}}/\ol{\syn{b}_{\ul{i}}}}\rangle\\
    (\syn{\pi}_i\{\syn{\mu}\})\psb{\ol{\syn{\nu}_{\ul{i}}}/\ol{\syn{b}_{\ul{i}}}}
    &\equiv \syn{\pi}_i\{\syn{\mu}\psb{\ol{\syn{\nu}_{\ul{i}}}/\ol{\syn{b}_{\ul{i}}}}\}\\
    \langle\ \rangle\psb{\ol{\syn{\nu}_{\ul{i}}}/\ol{\syn{b}_{\ul{i}}}}
    &\equiv \langle\ \rangle\\
\end{align*}
}

In the above, we use overline notation to denote the concatenation of
terms, protypes, or proterms by $\smcl$,
and we use underlined notation to specify the range of indices
traversing the concatenation.
For example, we write
$\syn{\kappa}(\syn{s}_0\smcl\dots\smcl\syn{s}_n)\{\syn{\mu_1}\smcl\dots\smcl\syn{\mu_n}\}$ as $\syn{\kappa}(\ol{\syn{s}_{\ul{i}}})\{\ol{\syn{\mu}_{\ul{i}}}\}$.
These are interpreted as sequences aligned
in horizontal direction in a \ac{VDC}.
Note that a mere sequence of terms in a context, for instance, is not written 
with the overline notation.

\begin{lemma}[Substitution lemmas]
    \label{lem:subst}
    The following equations hold for substitution and prosubstitution.
    \begin{enumerate}
        \item $\syn{\alpha}\left[\syn{S}/\syn{\Delta}\smcl\syn{T}/\syn{\Theta}\right]\left[\syn{S}'/\syn{\Delta}'\smcl\syn{T}'/\syn{\Theta}'\right]\equiv
        \syn{\alpha}\left[\syn{S}\left[\syn{S}'/\syn{\Delta}'\right]/\syn{\Delta}\smcl\syn{T}\left[\syn{T}'/\syn{\Theta}'\right]/\syn{\Theta}\right]$. 
        \item $\syn{\mu}\left[\ol{\syn{S}_{\ul{i}}}/\ol{\syn{\Delta}_{\ul{i}}}\right]\left[\ol{\syn{S}'_{\ul{i}}}/\ol{\syn{\Delta}'_{\ul{i}}}\right] 
        \equiv \syn{\mu}\left[\ol{\syn{S}_{\ul{i}}\left[\syn{S}'_{i}/\syn{\Delta}_{i}\right]}/\ol{\syn{\Delta}_{\ul{i}}}\right]$. 
        \item $\syn{\mu}\psb{\ol{\syn{\nu}_{\ul{i}}}/\ol{\syn{b}_{\ul{i}}}}\psb{\ol{\syn{\nu}'_{\ul{i,j}}}/\ol{\syn{b}'_{\ul{i,j}}}}
        \equiv \syn{\mu}\psb{\left.\ol{\syn{\nu}_{\ul{i}}\psb{\ol{\syn{\nu}'_{i,\ul{j}}}/\ol{\syn{b}'_{i,\ul{j}}}}} \right/ \ol{\syn{b}_{\ul{i}}}}$. 
        \item $\syn{\mu}\psb{\ol{\syn{\nu}_{\ul{i}}}/\ol{\syn{b}_{\ul{i}}}}\left[\ol{\syn{S}_{\ul{i,j}}}/\ol{\syn{\Delta}_{\ul{i,j}}}\right]
        \equiv \left(\syn{\mu}\left[\ol{\syn{S}_{\ul{i},n_{\ul{i}}}}/\ol{\syn{\Delta}_{\ul{i},n_{\ul{i}}}}\right]\right)\psb{\left.\ol{\syn{\nu}_{\ul{i}}\left[\ol{\syn{S}_{i,\ul{j}}}/\ol{\syn{\Delta}_{i,\ul{j}}}\right]}\right/\ol{\syn{b}_{\ul{i}}}}$. 
    \end{enumerate}
\end{lemma}
\begin{proof}
    The proof is straightforward by induction on the structure of terms, protypes, and proterms. 
\end{proof}

%% file: newsubfiles/typetheory/semantics.tex
As previously mentioned, the semantics of \ac{FVDblTT} are taken in \acp{CFVDC}.
The elements in the type theory are to be interpreted as the following elements in a \ac{CFVDC} $\dbl{D}$:
    \begin{itemize}
        \item $\syn{I}\ \textsf{type}$ 
        and $\syn{\Gamma}\ \textsf{ctx}$
        are to be interpreted as an object $\sem{\syn{I}}$
        and $\sem{\syn{\Gamma}}$ in $\dbl{D}$, respectively.
        \item $\syn{\Gamma}\vdash \syn{t}:\syn{I}$ and $\syn{\Gamma}\vdash \syn{S}\,/\,\syn{\Delta}$ are to be interpreted as tight arrows 
        $\sem{\syn{t}}\colon\sem{\syn{\Gamma}}\to \sem{\syn{I}}$ and $\sem{\syn{S}}\colon\sem{\syn{\Gamma}}\to \sem{\syn{\Delta}}$ in $\dbl{D}$, respectively.  
        \item $\syn{\Gamma}\smcl\syn{\Delta}\vdash \syn{\alpha}\ \textsf{protype}$ is to be interpreted as a loose arrow
        $\sem{\syn{\alpha}}\colon\sem{\syn{\Gamma}}\sto\sem{\syn{\Delta}}$ in $\dbl{D}$.
        \item $\syn{\Gamma}_0\smcl\dots\smcl\syn{\Gamma}_n\mid \syn{a}_1:\syn{\alpha}_1\smcl\dots\smcl\syn{a}_n:\syn{\alpha}_n\ \textsf{proctx}$ is 
        to be interpreted as a path of loose arrows 
        \[ 
        \sem{\syn{\Gamma}_0}\stonormal["\sem{\syn{\alpha}_1}"]\sem{\syn{\Gamma}_1}\stonormal \dots \stonormal["\sem{\syn{\alpha}_n}"]\sem{\syn{\Gamma}_n}\quad\text{in}\ \dbl{D}.
        \]
        \item $\syn{\ol{\Gamma}}\mid \syn{a}_1:\syn{\alpha}_1\smcl\dots\smcl\syn{a}_n:\syn{\alpha}_n\vdash \syn{\mu}:\syn{\beta}$ is to be interpreted as a globular cell 
        $\sem{\syn{\mu}}\colon\ol{\sem{\syn{\alpha}_{\ul{i}}}}\Rightarrow\sem{\syn{\beta}}$ in $\dbl{D}$.

    \end{itemize}

%% file: newsubfiles/typetheory/inddefsem.tex
The semantics of \ac{FVDblTT} consists of two parts: 
assignment of data in a \ac{CFVDC} to the ingredients of a signature,
and inductive definition of the interpretation.

\begin{definition}
    For a signature $\Sigma$ and a \ac{CFVDC} $\dbl{D}$,
    a \emph{$\Sigma$-structure} $\one{M}$ in $\dbl{D}$ is a morphism of signatures $\Sigma\to\Sigma_{\dbl{D}}$.
    The identity morphism on $\Sigma_{\dbl{D}}$ can be deemed a $\Sigma_{\dbl{D}}$-structure in $\dbl{D}$,
    which we call the \emph{canonical ($\Sigma_\dbl{D}$-)structure} in $\dbl{D}$.
\end{definition}

Instead of writing $\one{M}(\syn{\sigma})$ for the image of a category symbol $\syn{\sigma}$ under $\one{M}$, 
we write $\sem{\syn{\sigma}}_{\one{M}}$, or simply $\sem{\syn{\sigma}}$ when $\one{M}$ is clear from the context.

\begin{definition}
\label{def:semantics}
    Suppose we are given a $\Sigma$-structure $\one{M}$ in a \ac{CFVDC} $\dbl{D}$.
    The interpretation of the terms, protypes, protype isomorphisms, and proterms for $\Sigma$ 
    in $\dbl{D}$ is defined inductively as follows:
\begin{itemize}
    \item The interpretation of the type $\syn{\sigma}$ is the object $\sem{\syn{\sigma}}$ of $\dbl{D}$.
    \item The interpretation of the context $\cdot$ is the terminal object of $\dbl{D}$. 
    \item The interpretation of the context $\syn{\Gamma},\syn{x}:\syn{\sigma}$ is the product
    $\sem{\syn{\Gamma}}\times\sem{\syn{\sigma}}$ of $\sem{\syn{\Gamma}}$ and $\sem{\syn{\sigma}}$. 
    \item The interpretation of the term $\syn{\Gamma},\syn{x}:\syn{\sigma},\syn{\Delta}\vdash\syn{x}:\syn{\sigma}$
    is the projection onto $\sem{\syn{\sigma}}$. 
    \item The interpretation of the term $\syn{f}(\syn{t})$ is the composite $\sem{\syn{f}}\circ\sem{\syn{t}}$
    of the tight arrows $\sem{\syn{f}}\colon\sem{\syn{\sigma}}\to\sem{\syn{\tau}}$ and $\sem{\syn{t}}\colon\sem{\syn{\Gamma}}\to\sem{\syn{\sigma}}$. 
    \item Product types $\times,\syn{1}$ are interpreted as the product and terminal object of $\dbl{D}$, respectively.
    Pairing, projections, and the unit are interpreted in an obvious way. 
    \item The interpretation of the protype $\syn{\rho}(\syn{s}\smcl\syn{t})$ 
    is the restriction of the loose arrow $\sem{\syn{\rho}}:\sem{\syn{\sigma}}\sto\sem{\syn{\tau}}$ along the tight arrows 
    $\sem{\syn{s}}:\sem{\syn{\Gamma}}\to\sem{\syn{\sigma}}$ and $\sem{\syn{t}}:\sem{\syn{\Delta}}\to\sem{\syn{\tau}}$. 
    \[
        \begin{tikzcd}[column sep=12ex,virtual]
            \sem{\syn{\Gamma}}
            \sar[r, "{\sem{\syn{\rho}(\syn{s}\smcl\syn{t})}}"]
            \ar[d, "\sem{\syn{s}}"']
            \ar[dr,phantom, "\restc" description]
            &
            \sem{\syn{\Delta}}
            \ar[d, "\sem{\syn{t}}"]
            \\
            \sem{\syn{\sigma}}
            \sar[r, "\sem{\syn{\rho}}"']
            &
            \sem{\syn{\tau}}
        \end{tikzcd}
    \]
    \item Product protypes $\land,\top$ in context $\syn{\Gamma}\smcl\syn{\Delta}$ 
    are interpreted as the product and terminal loose arrow from $\sem{\syn{\Gamma}}$ to $\sem{\syn{\Delta}}$, respectively.
    Pairing, projections, and the unit are interpreted in an obvious way.
    \item The interpretation of the proterm $\syn{a}:\syn{\alpha}\vdash \syn{a}:\syn{\alpha}$ is the identity cell on $\sem{\syn{\alpha}}$.
    \item To define the interpretation of the proterm $\syn{\kappa}(\ol{\syn{s}_{\ul{i}}})\{\ol{\syn{\mu}_{\ul{i}}}\}$,
    we first define a cell $\sem{\syn{\kappa}(\ol{\syn{s}_{\ul{i}}})}$ 
    as the restriction $\sem{\syn{\kappa}}\left[\ol{\sem{\syn{s}_{\ul{i}}}}\right]\colon\sem{\syn{\rho}_1(\syn{s}_0\smcl\syn{s}_1)}\smcl\dots\smcl\sem{\syn{\rho}_n(\syn{s}_{n-1}\smcl\syn{s}_n)}\Rightarrow\sem{\syn{\omega}(\syn{s}_0\smcl\syn{s}_n)}$
    in the sense of \Cref{def:restrictioncell}.
    Then, the interpretation of the proterm $\syn{\kappa}(\ol{\syn{s}_{\ul{i}}})\{\ol{\syn{\mu}_{\ul{i}}}\}$ is the composite 
    $\sem{\syn{\kappa}(\ol{\syn{s}_{\ul{i}}})}\{\sem{\syn{\mu}_{1}}\smcl\dots\smcl\sem{\syn{\mu}_{n}}\}$ 
    of the cell $\sem{\syn{\kappa}(\ol{\syn{s}_{\ul{i}}})}$ 
    and the cells $\sem{\syn{\mu}_{i}}\colon\sem{\syn{A}_i}\Rightarrow\sem{\syn{\rho}_i(\syn{s}_{i-1}\smcl\syn{s}_i)}$
    for $i=1,\dots,n$.
\end{itemize}
\end{definition}

Taking semantics in the \acp{VDC} listed in \Cref{example:prof,example:rel} justifies how \ac{FVDblTT} expresses
formal category theory and predicate logic.

We have naively used restrictions in the interpretation of protypes,
but they are only defined up to isomorphism \textit{a priori}.
To make the definition precise, we need to consider strict functoriality in the following sense.
\begin{definition}
    A \ac{CFVDC} $\dbl{D}$ is \emph{split} if it comes with
    chosen finite products of its tight category,
    chosen restrictions $(-)[-\smcl-]$,
    and chosen terminals $\top$ and binary products $(-)\land(-)$ in the loose hom-categories.
    that satisfy the following equalities:
    \begin{itemize}
        \item $\alpha[\id_I\smcl\id_J]=\alpha$ for any $\alpha\colon I\sto J$.
        \item $\alpha[s\smcl t][s'\smcl t']=\alpha[s\circ s'\smcl t\circ t']$ for any $\alpha\colon I\sto J$ and $s,t,s',t'$.
        \item $\top[s\smcl t]=\top$ for any $s,t$.
        \item $(\alpha\land\beta)[s\smcl t]=\alpha[s\smcl t]\land\beta[s\smcl t]$ for any $\alpha,\beta\colon I\sto J$ and $s,t$. 
    \end{itemize}
    A \emph{morphism of split} \acp{CFVDC} is a 1-cell in $\FibVDblCart$
    that preserves the chosen tightwise finite products, restrictions, terminals, and binary products on the nose.
    We will denote the category of split \acp{CFVDC} by $\FibVDblCart^{\spl}$.
\end{definition}
Note that in a split \ac{CFVDC}, restrictions of globular cells along tight arrows in \Cref{def:restrictioncell}
are uniquely determined by the chosen restrictions.
\begin{lemma}
    \label{lemma:interpretsubs}
    Let $\dbl{D}$ be a split \ac{CFVDC}, and let $\one{M}$ be a $\Sigma$-structure in $\dbl{D}$.
    Suppose we choose the chosen restrictions in $\dbl{D}$ in the definition of the interpretation.
    \begin{enumerate}
    \item The interpretation of term substitutions is obtained by
    restrictions of loose arrows or globular cells along tight arrows.
    Explicitly, 
    we have 
    $\sem{\syn{\alpha}[\syn{S}/\syn{\Gamma}\smcl\syn{T}/\syn{\Delta}]}=\sem{\syn{\alpha}}[\sem{\syn{S}}\smcl\sem{\syn{T}}]$
    and $\sem{\syn{\mu}[\ol{\syn{S}_{\ul{i}}}/\ol{\syn{\Delta}_{\ul{i}}}]}=\sem{\syn{\mu}}[\ol{\sem{\syn{S}_{\ul{i}}}}]$
    whenever the substitutions are well-typed.
    \item The interpretation of proterm prosubstitutions is obtained by
    composition of globular cells.
    Explicitly, we have
    $\sem{\syn{\mu}\psb{\ol{\syn{\nu}_{\ul{i}}}/\ol{\syn{b}_{\ul{i}}}}}=\sem{\syn{\mu}}\{\ol{\sem{\syn{\nu}_{\ul{i}}}}\}$
    whenever the prosubstitutions are well-typed.
    \end{enumerate}
\end{lemma}
\begin{proof}
    By induction on the structure of term substitutions and prosubstitutions.
\end{proof}

Assuming splitness for a \ac{CFVDC} is too strong for practical purposes,
but we can replace an arbitrary \ac{CFVDC} by an equivalent split one.
\begin{lemma}
    \label{lemma:split}
    For any \ac{CFVDC} $\dbl{D}$, there exists a split \ac{CFVDC} $\dbl{D}^{\spl}$
    that is equivalent to $\dbl{D}$ in the 2-category $\FibVDblCart$.
\end{lemma}
\begin{proof}[Proof sketch]
    The proof is analogous to the proof for splitness of fibrational virtual double categories
    in \cite[Theorem A.1]{arkorNerveTheoremRelative2024}.
    For a \ac{CFVDC} $\dbl{D}$, 
    fix chosen terminals and binary products in each loose hom-category
    and chosen restrictions.
    We define a split \ac{CFVDC} $\dbl{D}^{\spl}$ 
    by taking the same objects and tight arrows as $\dbl{D}$,
    but its loose arrows from $I$ to $J$ are
    finite tuples of triples $(f_i,g_i,\alpha_i)_{i}$ where 
    $f_i\colon I\to K_i$ and $g_i\colon J\to L_i$ are tight arrows and 
    $\alpha_i\colon K_i\sto L_i$ are loose arrows in $\dbl{D}$.
    From a loose arrow $(f_i,g_i,\alpha_i)_{i}$,
    we can define its realization in $\dbl{D}$ by taking $\bigwedge_i\alpha_i[f_i\smcl g_i]$.
    Then, we can define cells in $\dbl{D}^{\spl}$ framed by two tight arrows and loose arrows 
    as those in $\dbl{D}$ framed by the same tight arrows and the realization of the corresponding loose arrows.
    The associativity and unitality of cell composition in $\dbl{D}^{\spl}$ are inherited from those in $\dbl{D}$. 
    There is a virtual double functor $\dbl{D}^{\spl}\to\dbl{D}$ that
    is the identity on the tight part and sends a loose arrow to its realization 
    and a cell to itself.
    This is an equivalence of virtual double categories.
    To verify that $\dbl{D}^{\spl}$ admits the structure of a split \ac{CFVDC},
    we define the chosen restrictions, terminals, and binary products in $\dbl{D}^{\spl}$ 
    as follows:
    \begin{itemize}
        \item The restriction of a loose arrow $(f_i,g_i,\alpha_i)_{i}$ along a pair of tight arrows 
        $(h,k)$ is the tuple $(f_i\circ h,g_i\circ k,\alpha_i)_i$.
        \item The terminal object in the loose hom-category from $I$ to $J$ is
        the empty tuple.
        \item The binary product of two loose arrows $(f_i,g_i,\alpha_i)_{i\in\zero{I}}$ and $(f'_j,g'_j,\alpha'_j)_{j\in\zero{J}}$ is 
        the sum of the two tuples.
    \end{itemize}
    It is straightforward to verify that these chosen structures strictly satisfy the equalities in the definition of split \acp{CFVDC}. 
\end{proof}

%% file: newsubfiles/typetheory/proiso.tex
In category theory, one often proves that two objects, functors, or profunctors are isomorphic
by exhibiting a sequence of those isomorphisms between them that one has already constructed or known to exist.
Protype isomorphisms enable us to do the same in the type theory
without showing proterms in both directions explicitly every time
but still keeping track of the proterms that represent the isomorphisms.
We introduce protype isomorphisms as additional typing judgments 
but they also serve partially as equality judgments for protypes up to isomorphism.
Protype isomorphisms are also considered as codes for the two proterms mutually inverse to each other
so that proterms can track what they actually represent in the type theory.
They are also used to express isomorphisms between functors (terms) as 
we will see in \Cref{sec:examples}.
It should be noted that we do not have equality judgments for protype isomorphisms since one can identify or distinguish them
by the proterms they represent using the equality judgments for proterms.

We call this extension of the type theory with protype isomorphisms \ac{FVDblTT}$^{\ccong}$.
The judgments for protype isomorphisms are presented as
$
    \syn{\Gamma}\smcl\syn{\Delta}\vdash \syn{\Upsilon}:\syn{\alpha}\ccong\syn{\beta}
$ 
where $\syn{\alpha}$ and $\syn{\beta}$ are protypes in the context $\syn{\Gamma}\smcl\syn{\Delta}$.  
The rules for protype isomorphisms are given as follows:

\begin{mathparpagebreakable}
    \inferrule*
    {\syn{\Gamma} \smcl \syn{\Delta} \vdash \syn{\alpha} \ \textsf{protype}}
    {\syn{\Gamma} \smcl \syn{\Delta} \vdash \idt_{\syn{\alpha}}:\syn{\alpha}\ccong\syn{\alpha}}
    \and
    \inferrule*
    {\syn{\Gamma} \smcl \syn{\Delta} \vdash \syn{\Upsilon}:\syn{\alpha}\ccong\syn{\beta}}
    {\syn{\Gamma} \smcl \syn{\Delta} \vdash \syn{\Upsilon}\sinv:\syn{\beta}\ccong\syn{\alpha}}
    \and
    \inferrule*
    {\syn{\Gamma} \smcl \syn{\Delta} \vdash \syn{\Upsilon}:\syn{\alpha}\ccong\syn{\beta} \\ 
    \syn{\Gamma} \smcl \syn{\Delta} \vdash \syn{\Omega}:\syn{\beta}\ccong\syn{\gamma}}
    {\syn{\Gamma} \smcl \syn{\Delta} \vdash \syn{\Omega}\circ\syn{\Upsilon}:\syn{\alpha}\ccong\syn{\gamma}}
    \and
    \inferrule*
    {\syn{\Gamma}\smcl\syn{\Delta} \mid \syn{a}:\syn{\alpha} \vdash \syn{\mu}\{\syn{a}\}:\syn{\beta} \\
    \syn{\Gamma}\smcl\syn{\Delta} \mid \syn{b}:\syn{\beta} \vdash \syn{\nu}\{\syn{b}\}:\syn{\alpha}\\
    \syn{\Gamma}\smcl\syn{\Delta} \mid \syn{b}:\syn{\beta} \vdash \syn{\mu}\{\syn{\nu}\{\syn{b}\}\}\equiv \syn{b}:\syn{\beta}\\
    \syn{\Gamma}\smcl\syn{\Delta} \mid \syn{a}:\syn{\alpha} \vdash \syn{\nu}\{\syn{\mu}\{\syn{a}\}\}\equiv \syn{a}:\syn{\alpha}}
    {\syn{\Gamma}\smcl\syn{\Delta} \vdash \lcp\syn{\mu},\syn{\nu}\rcp:\syn{\alpha}\ccong\syn{\beta}}         
    \and
    \inferrule*
    {\syn{\Gamma}\smcl\syn{\Delta}\vdash \syn{\Upsilon}:\syn{\alpha}\ccong\syn{\beta}}
    {\syn{\Gamma}\smcl\syn{\Delta}\mid \syn{a}:\syn{\alpha} \vdash \syn{\Upsilon}\{\syn{a}\}:\syn{\beta}}
    \and
    \inferrule*
    {\syn{\Gamma}\smcl\syn{\Delta}\vdash \syn{\alpha} \ \textsf{protype}}
    {\syn{\Gamma}\smcl\syn{\Delta}\mid \syn{a}:\syn{\alpha}\vdash \idt_{\syn{\alpha}}\{\syn{a}\}\equiv \syn{a}:\syn{\alpha}}
    \and
    \inferrule*
    {\syn{\Gamma}\smcl\syn{\Delta}\vdash \syn{\Upsilon}:\syn{\alpha}\ccong\syn{\beta}}
    {\syn{\Gamma}\smcl\syn{\Delta}\mid \syn{a}:\syn{\alpha}\vdash \syn{\Upsilon}\sinv\{\syn{\Upsilon}\{\syn{a}\}\}\equiv \syn{a}:\syn{\alpha}}
    \and  
    \inferrule*
    {\syn{\Gamma}\smcl\syn{\Delta}\vdash \syn{\Upsilon}:\syn{\alpha}\ccong\syn{\beta}}
    {\syn{\Gamma}\smcl\syn{\Delta}\mid \syn{b}:\syn{\beta} \vdash \syn{\Upsilon}\{\syn{\Upsilon}\sinv\{\syn{a}\}\}\equiv \syn{a}:\syn{\alpha}}
    \and
    \inferrule*
    {\syn{\Gamma}\smcl\syn{\Delta}\vdash \syn{\Upsilon}:\syn{\alpha}\ccong\syn{\beta}\\
    \syn{\Gamma}\smcl\syn{\Delta}\vdash \syn{\Omega}:\syn{\beta}\ccong\syn{\gamma}}
    {\syn{\Gamma}\smcl\syn{\Delta}\mid \syn{a}:\syn{\alpha}\vdash (\syn{\Omega}\circ\syn{\Upsilon})\{\syn{a}\}\equiv
    \syn{\Omega}\{\syn{\Upsilon}\{\syn{a}\}\}:\syn{\gamma}}
    \and     
    \inferrule*
    {\syn{\Gamma}\smcl\syn{\Delta} \mid \syn{a}:\syn{\alpha} \vdash \syn{\mu}\{\syn{a}\}:\syn{\beta} \\
    \syn{\Gamma}\smcl\syn{\Delta} \mid \syn{b}:\syn{\beta} \vdash \syn{\nu}\{\syn{b}\}:\syn{\alpha}\\
    \syn{\Gamma}\smcl\syn{\Delta} \mid \syn{b}:\syn{\beta} \vdash \syn{\mu}\{\syn{\nu}\{\syn{b}\}\}\equiv \syn{b}:\syn{\beta}\\
    \syn{\Gamma}\smcl\syn{\Delta} \mid \syn{a}:\syn{\alpha} \vdash \syn{\nu}\{\syn{\mu}\{\syn{a}\}\}\equiv \syn{a}:\syn{\alpha}}
    {\syn{\Gamma}\smcl\syn{\Delta} \mid \syn{a}:\syn{\alpha} \vdash \lcp\syn{\mu},\syn{\nu}\rcp\{\syn{a}\}\equiv \syn{\mu}\{\syn{a}\}:\syn{\beta}}
\end{mathparpagebreakable}

If one has a pair of proterms $\syn{\mu}$ and $\syn{\nu}$ that are mutually inverse to each other, 
one can form a protype isomorphism $\lcp\syn{\mu},\syn{\nu}\rcp$.
Conversely, protype isomorphisms are realized as proterms via the rule that introduces
the proterm $\syn{\Upsilon}\{\syn{a}\}$ for a protype isomorphism $\syn{\Upsilon}$.
We have the rule $\lcp\syn{\mu},\syn{\nu}\rcp\{\syn{a}\}\equiv\syn{\mu}\{\syn{a}\}$,
which is sufficient to derive that the inverse of $\lcp\syn{\mu},\syn{\nu}\rcp$ also has the expected behavior: 
    $
        \lcp\syn{\mu},\syn{\nu}\rcp\inv\{\syn{b}\} 
        \equiv \lcp\syn{\mu},\syn{\nu}\rcp\inv\left\{\syn{\mu}\left\{\syn{\nu}\{\syn{b}\}\right\}\right\}
        \equiv \lcp\syn{\mu},\syn{\nu}\rcp\inv\left\{\lcp\syn{\mu},\syn{\nu}\rcp\left\{\syn{\nu}\{\syn{b}\}\right\}\right\}
        \equiv \syn{\nu}\{\syn{b}\}
    $.
The other rules are designed to ensure that protype isomorphisms behave
as a groupoid as a whole.

The semantics of \ac{FVDblTT}$^{\ccong}$ are also given in a \ac{CFVDC}.
A protype isomorphism judgment
$\syn{\Gamma}\smcl\syn{\Delta}\vdash \syn{\Upsilon}:\syn{\alpha}\ccong\syn{\beta}$ is
to be interpreted as an isomorphism of loose arrows
$\sem{\syn{\Upsilon}}\colon\sem{\syn{\alpha}}\Rightarrow\sem{\syn{\beta}}\colon\sem{\syn{\Gamma}}\sto\sem{\syn{\Delta}}$ 
in $\dbl{D}$.
The interpretations of the protype isomorphisms $\idt_{\syn{\alpha}},\syn{\Upsilon}\sinv,\syn{\Upsilon}\circ\syn{\Omega}$ 
are defined as the identity cell, the inverse cell, and the composite cell
of the corresponding cells in $\dbl{D}$,
and the interpretation of the protype isomorphism $\lcp\syn{\mu},\syn{\nu}\rcp$ is the cell $\sem{\syn{\mu}}$.

%% file: newsubfiles/constructor/tableofconstructors.tex
In this section, we will specify the type and protype constructors that can be added to \ac{FVDblTT}.
The virtual double categories of relations and those of profunctors have many structures in common.
We would like to introduce the inductive types and protypes corresponding to the common structures in these kinds of virtual double categories.
We first list the additional types will introduce for the type theory.

\begin{table}[h] 
    \rowcolors{2}{gray!25}{white}
    \begin{tabular}{|c||c|c||c|}
        \hline
        \rowcolor{gray!50}
        \textbf{Structures} & \textbf{Formal category theory} & \textbf{Predicate logic} & \begin{tabular}{@{}c@{}}\textbf{Constructors}\\ \textbf{in FVDblTT}\end{tabular}\\
        \hline
        \hline
        Units \cite{cruttwellUnifiedFrameworkGeneralized2010} & hom-profuntors $\one{C}(-,\bullet)$ & equality $=$ & path $\ides$ \\
        Composition \cite{cruttwellUnifiedFrameworkGeneralized2010} & composition via coends $\int$ & composition via $\exists$ & composition $\odot$ \\
        Extension \cite{riehlElementsCategoryTheory2022} & profunctor extension $\triangleright$ & contraction via $\forall$ & extension $\triangleright$ \\
        Tabulators \cite{grandisLimitsDoubleCategories1999} & \begin{tabular}{@{}c@{}}two-sided\\ Grothendieck construction\end{tabular} & comprehension $\{\mhyphen\}$ & tabulator $\cmpr{\mhyphen}$ \\
        \hline
    \end{tabular} 
    \caption{The common structures and the corresponding constructors}
    \label{table:commonconst}
\end{table}

The constructors we will add to \ac{FVDblTT} are $\ides$, $\odot$, $\triangleright$, $\triangleleft$, and $\cmpr{\mhyphen}$.
Even though we can add the constructors for the loose adjunctions and the companions and conjoints independently of the other constructors,
we would take the approach of defining them in terms of $\ides$ and $\odot$ in this paper.

%% file: newsubfiles/constructor/units.tex
\paragraph*{Path protype $\ides$ for the units.\ ({\Cref{sec:unit-protype,sec:unit-protype-meets-product-type}})}
The path protype is the protype that represents the units in a \ac{VDC}.
In a double category, the units are just the identity loose morphisms,
but in a \ac{VDC}, the units are formulated via a universal property.
The definition of units is due to \cite{cruttwellUnifiedFrameworkGeneralized2010}.

\begin{definition}[Units {\cite[Definition 5.1]{cruttwellUnifiedFrameworkGeneralized2010}}]
    A \emph{unit} of an object $I$ in a \ac{VDC} is a loose arrow $U_I\colon I\sto I$ equipped with a cell 
    $\eta_I\colon \cdot\Rightarrow U_I$
    with the following universal property.
    Given any cell $\nu$ on the left below where $\ol\alpha=\left(I_0\sto["\alpha_1"] \cdots\sto["\alpha_n"] I\right)$ and 
    $\ol\alpha'=\left(I\sto["\alpha'_1"] \cdots\sto["\alpha'_{n'}"] I'_{n'}\right)$
    are arbitrary sequences of loose arrows,
    it uniquely factors through the sequence of the identity cells with $\eta_I$ as on the right below.
    \[
        \begin{tikzcd}[virtual, column sep=small]
            I_0
            \sar[rr, "\ol\alpha",dashed]
            \ar[d, "f"']
            \ar[phantom,rrrrd, "\nu",description]
            &
            &
            I
            \sar[rr, "\ol\alpha'",dashed]
            &&
            {I'_{n'}}
            \ar[d, "f'"]
            \\
            J
            \sar[rrrr, "\beta"']
            &&&&
            J'
        \end{tikzcd}=
        \begin{tikzcd}[virtual, column sep=small]
            I_0
            \sar[rr, "\ol\alpha",dashed]
            \ar[d, equal]
            \ar[dr, phantom, description, "{\rotatebox{90}{=}}",xshift=1ex]
            &
            &
            I
            \ar[ld, equal]
            \ar[rd, equal]
            \ar[d,phantom, "\eta_I"]
            \sar[rr, "\ol\alpha'",dashed]
            &
            &
            I'_{n'}
            \ar[d, equal]
            \ar[dl, phantom, description, "{\rotatebox{90}{=}}",xshift=-1ex]
            \\
            I_0
            \sar[r, "\ol\alpha"',dashed]
            \ar[d, "f"']
            \ar[phantom,rrrrd, "\wt{\nu}",description, yshift=-1ex]
            &
            I
            \sar[rr, "U_I"']
            &\!&
            I
            \sar[r, "\ol\alpha'"',dashed]
            &
            I'_{n'}
            \ar[d, "f'"]
            \\
            J
            \sar[rrrr, "\beta"']
            &&&\!&
            J'
        \end{tikzcd}
    \]
\end{definition}

The formation rule for the \emph{path protype} is on the left below, and it comes equipped with the introduction rule on the right below: 
    
    \begin{mathparpagebreakable}
    \inferrule*[right=$\ides$-Form]
    {\syn{I}\ \textsf{type}\\
    \syn{\Gamma}\vdash \syn{s}:\syn{I}\\
    \syn{\Delta}\vdash \syn{t}:\syn{I}}
    {\syn{\Gamma}\smcl\syn{\Delta}\vdash \syn{s}\ide{\syn{I}}\syn{t}\ \textsf{protype}}
    \and
    \inferrule*
    {\syn{I}\ \textsf{type}}
    {\syn{x}:\syn{I}\mid \quad \vdash \refl_{\syn{I}}(\syn{x}): \syn{x}\ide{\syn{I}}\syn{x}} 
    \end{mathparpagebreakable}
The proterm $\refl$ corresponds to the unit $\eta_I$ in the definition of the units.
To let the path protype encode the units in the \acp{VDC},
we need to add elimination and computation rules as in \Cref{sec:unit-protype}.
The path protypes behave as inductive (pro)types,
and their inductions look very similar to path induction in homotopy type theory,
but with the difference that the path protype is directed.

The semantics of the path protypes $\ides$ are given by the units in any \ac{VDC} with units,
with the proterm constructor $\refl_{\syn{I}}$ interpreted as the cell $\eta_{\sem{\syn{I}}}$.
For instance, in the \acp{VDC} $\Prof$ and $\Rel{}$, the interpretations of the path protypes are given as
the hom profunctors and the equality relations, respectively.
These follow from the fact that the identity loose morphisms in a double category serve as the units when we see it as a \ac{VDC}.

In order to make the path protypes behave well with the product types in \ac{FVDblTT},
we need to add the compatibility rules between the path protypes and the product types
as in \Cref{sec:unit-protype-meets-product-type}.
For instance, when we consider the hom-profunctor on a product category $\one{C}\times\one{D}$,
we expect its components to be isomorphic to the product $\one{C}(C,C')\times\one{D}(D,D')$.
Correspondingly, we would like to add the following rule, which does not follow from other rules \textit{a priori}:

    \begin{mathparpagebreakable}
    \inferrule*
    {\syn{I}\ \textsf{type}\\ \syn{J}\ \textsf{type}}
    {\syn{x}:\syn{I},\syn{y}:\syn{J}\smcl\syn{x}':\syn{I},\syn{y}':\syn{J}
    \vdash \exc_{\ides,\land}:\langle\syn{x},\syn{y}\rangle\ide{\syn{I}\times\syn{J}}\langle\syn{x'},\syn{y'}\rangle
    \ccong \syn{x}\ide{\syn{I}}\syn{x}'\land\syn{y}\ide{\syn{J}}\syn{y'}}.
    \end{mathparpagebreakable}
\Cref{sec:unit-protype-meets-product-type} will give the whole set of rules for the compatibility between the path protypes and the product types. 
The rules we introduce are justified by the fact that 
with them, the syntactic \acp{VDC} we will introduce in \Cref{sec:synsemadj}
become cartesian objects in the 2-category of \acp{CFVDC} with units.
See \Cref{prop:cartesianunital} for a detailed explanation from the 2-categorical perspective.

%% file: newsubfiles/constructor/others.tex
\paragraph*{Composition protype $\odot$ for the composites.\ ({\Cref{sec:composition-protype,sec:compo-protype-meets-product-type}})}
The composition protype is the protype that represents the composition of 
paths of loose arrows just of length 2 in the virtual double categories.
Here we follow the definition of the composition of paths of loose arrows
in \cite{cruttwellUnifiedFrameworkGeneralized2010}.
\begin{definition}[{Composites \cite[Definition 5.2]{cruttwellUnifiedFrameworkGeneralized2010}}]
    A \emph{composite} of a given sequence of loose arrows
\[\ol\alpha=\left(I_0\sto["\alpha_1"] I_{1}\sto\cdots\sto["\alpha_m"] I_m\right)\] in a virtual double category 
    is a loose arrow $\odot\ol\alpha$ from $I_0$ to $I_m$ equipped with a cell
    \[
        \begin{tikzcd}[virtual, column sep=small]
            I_0
            \sar[r, "\alpha_1"]
            \ar[d, equal]
            \ar[rrrd, phantom, "\varrho_{\ol\alpha}"]
            &
            I_1
            \sar[r]
            &
            \cdots
            \sar[r, "\alpha_m"]
            &
            I_m
            \ar[d, equal]
            \\
            I_0
            \sar[rrr, "\odot\ol\alpha"']
            &&&
            I_m
        \end{tikzcd}
    \]
    with the following universal property.
    Given any cell $\nu$ on the left below where $\ol\beta,\ol\beta'$ 
    are arbitrary sequences of loose arrows,
    it uniquely factors through the sequence of the identity cells with $\mu_{\ol\alpha}$ as on the right below.
    \[
        \begin{tikzcd}[virtual, column sep=small]
            J_0
            \sard[r, "\ol\beta"]
            \ar[d, "f"']
            \ar[phantom,rrrrd, "\nu"]
            &
            I_0
            \sard[rr, "\ol\alpha"]
            &&
            I_m
            \sard[r, "\ol\beta'"]
            &
            J'_{n'}
            \ar[d, "f'"]
            \\
            K
            \sar[rrrr, "\gamma"']
            &&&&
            K'
        \end{tikzcd}
=
        \begin{tikzcd}[virtual, column sep=small]
            J_0
            \sard[r, "\ol\beta"]
            \ar[d, equal]
            \ar[dr, phantom, description, "{\rotatebox{90}{=}}"]
            &
            I_0
            \ar[d, equal]
            \ar[drr,phantom, "\varrho_{\ol\alpha}"]
            \sar[rr, "\ol\alpha",dashed]
            &
            &
            I_m 
            \sard[r, "\ol\beta'"]
            \ar[d, equal]
            \ar[dr, phantom, description, "{\rotatebox{90}{=}}"]
            &
            J'_{n'}
            \ar[d, equal]
            \\
            J_0
            \sard[r, "\ol\beta"']
            \ar[d, "f"']
            \ar[phantom,rrrrd, "\wt{\nu}",description, yshift=-1ex]
            &
            I_0
            \sar[rr, "\odot\ol\alpha"']
            &&
            I_m
            \sard[r, "\ol\beta'"']
            &
            J'_{n'}
            \ar[d, "f'"]
            \\
            K
            \sar[rrrr, "\gamma"']
            &&&&
            K'
        \end{tikzcd}
\]
\end{definition}
The units are the special cases of the composition of paths of length 0,
and the composition of paths longer than 2 can be realized by the iterated use of the composition of paths of length 2.
In order to gain access to the composition of paths of positive length in the type theory,
we introduce the \emph{composition protype} $\odot$ to \ac{FVDblTT}.
The formation rule for the composition protype is the following:
\[
    \inferrule*
    {\syn{w}:\syn{I}\smcl\syn{x}:\syn{J}\vdash \syn{\alpha}(\syn{w}\smcl\syn{x})\ \textsf{protype} \\ 
    \syn{x}:\syn{J}\smcl\syn{y}:\syn{K}\vdash \syn{\beta}(\syn{x}\smcl\syn{y})\ \textsf{protype}}
    {\syn{w}:\syn{I}\smcl\syn{y}:\syn{K}\vdash \syn{\alpha}(\syn{w}\smcl\syn{x})\odot_{\syn{x}: \syn{J}}\syn{\beta}(\syn{x}\smcl\syn{y})\ \textsf{protype}}
\]
This comes equipped with the introduction rule:
\[
    \inferrule*
    {\syn{w}:\syn{I}\smcl\syn{x}:\syn{J}\vdash \syn{\alpha}(\syn{w}\smcl\syn{x})\ \textsf{protype} \\
            \syn{x}:\syn{J}\smcl\syn{y}:\syn{K}\vdash \syn{\beta}(\syn{x}\smcl\syn{y})\ \textsf{protype} }
            {\syn{w}:\syn{I}\smcl\syn{x}:\syn{J}\smcl\syn{y}:\syn{K}\mid \syn{a}:\syn{\alpha}(\syn{w}\smcl\syn{x})\smcl\syn{b}:\syn{\beta}(\syn{x}\smcl\syn{y})\vdash
            \syn{a}\odot\syn{b}:\syn{\alpha}(\syn{w}\smcl\syn{x})\odot_{\syn{x}:\syn{J}}\syn{\beta}(\syn{x}\smcl\syn{y})}
\]
For the detailed rules of the composition protype, see \Cref{sec:composition-protype}.
Plus, we need the compatibility rules for the composition protype and the product types
as we did for the path protype, see \Cref{sec:compo-protype-meets-product-type}.

If we load the path protype $\ides$ and the composite protype $\odot$ to \ac{FVDblTT},
procontexts can be equivalently expressed by a single protype.
In this sense, such a type theory can be seen as an internal language of double categories.
This is supported by the fact that a \ac{VDC} is equivalent to one arising from a double category 
if and only if it has composites of all paths of loose arrows, including units
\cite[Theorem 5.2]{cruttwellUnifiedFrameworkGeneralized2010}.

The semantics of the composition protypes $\odot$ is given by the composites in \acp{VDC} if they have ones of sequences
of length 2 in an appropriate way.
For example, in the \ac{VDC} $\Prof$, the composite of paths of length 2 is the composite of profunctors,
given by the coend $\int$.
In the \ac{VDC} $\Rel{}$, the composites of paths of length 2 are the composites of relations,
given by the existential quantification $\exists$.
\[
    \begin{aligned}
        \sem{\syn{\alpha}(\syn{w}\smcl\syn{x})\odot_{\syn{x}:\syn{J}}\syn{\beta}(\syn{x}\smcl\syn{y})}&=
        \int^{X\in\sem{\syn{J}}}\sem{\syn{\alpha}}(-,X)\times\sem{\syn{\beta}}(X,\bullet)\colon\sem{\syn{I}}\sto\sem{\syn{K}}\  \text{in}\ \Prof{}\\
        \sem{\syn{\alpha}(\syn{w}\smcl\syn{x})\odot_{\syn{x}:\syn{J}}\syn{\beta}(\syn{x}\smcl\syn{y})}&=
        \left\{\,(w,y)\mid\exists x\in\sem{\syn{J}}.\sem{\syn{\alpha}}(w,x)\land\sem{\syn{\beta}}(x,y)\,\right\}\colon\sem{\syn{I}}\sto\sem{\syn{K}}\ \text{in}\ \Rel{}
    \end{aligned}
\]

\paragraph*{Filler protype $\triangleright,\triangleleft$ for the closed structure.\ ({\Cref{sec:filler-protype,sec:filler-protype-meets-product-type}})}
Having obtained the ability to express a particular kind of coends in formal category theory,
and existential quantification in predicate logic,
we would like to introduce the protypes for ends and universal quantification in the type theory.
First of all, we recall the definition of the right extension and the right lift \cite{riehlElementsCategoryTheory2022,arkorFormalTheoryRelative2024} in a \ac{VDC},
which are straightforward generalizations of the right extension and the right lift in a bicategory.

\begin{definition}
    A \emph{right extension} of a loose arrow $\beta\colon I\sto K$ along a loose arrow $\alpha\colon I\sto J$ 
    is a loose arrow $\alpha\triangleright\beta\colon J\sto K$ equipped with a cell
    \[
        \begin{tikzcd}[virtual, column sep=small]
            I
            \sar[r, "\alpha"]
            \ar[d, equal]
            \ar[rrd, phantom, "\varpi_{\alpha,\beta}"]
            &
            J
            \sar[r, "\alpha\triangleright\beta"]
            &
            K
            \ar[d, equal]
            \\
            I
            \sar[rr, "\beta"']
            &&
            K
        \end{tikzcd}
    \]
    with the following universal property.
    Given any cell $\nu$ on the left below where $\ol\gamma$ is an arbitrary sequence of loose arrows,
    it uniquely factors through the cell $\varpi_{\alpha,\beta}$ as on the right below.
    \[
        \begin{tikzcd}[virtual, column sep=small]
            I
            \sar[r, "\alpha"]
            \ar[d, equal]
            \ar[phantom,rrrd, "\nu"]
            &
            J
            \sard[rr, "\ol\gamma"]
            &&
            K 
            \ar[d, equal]
            \\
            I
            \sar[rrr, "\beta"']
            &&&
            K
        \end{tikzcd}
        =
        \begin{tikzcd}[virtual, column sep=small]
            I
            \sar[r, "\alpha"]
            \ar[d, equal]
            \ar[dr, phantom, description, "{\rotatebox{90}{=}}"]
            &
            J
            \ar[d, equal]
            \ar[drr,phantom, "\wt{\nu}"]
            \sard[rr, "\ol\gamma"]
            &&
            K
            \ar[d, equal]
            \\
            I
            \sar[r, "\alpha"']
            \ar[d, equal]
            \ar[phantom,rrrd, "\varpi_{\alpha,\beta}", yshift=-1ex]
            &
            J
            \sar[rr, "\alpha\triangleright\beta"']
            &
            &
            K
            \ar[d, equal]
            \\
            I
            \sar[rrr, "\beta"']
            &&&
            K
        \end{tikzcd}
    \]

    A \emph{right lift} of a protype $\beta\colon I\sto K$ along a protype $\alpha\colon J\sto K$
    is a protype $\beta\triangleleft\alpha\colon I\sto J$ equipped with a cell
    \[
        \begin{tikzcd}[virtual, column sep=small]
            I
            \sar[r, "{\beta\triangleleft\alpha}"]
            \ar[d, equal]
            \ar[rrd, phantom, "\varpi'_{\alpha,\beta}"]
            &
            J
            \sar[r, "\alpha"]
            &
            K
            \ar[d, equal]
            \\
            I
            \sar[rr, "\beta"']
            &&
            K
        \end{tikzcd}
    \]
    with the following universal property.
    Given any cell $\nu$ on the left below where $\ol\gamma$ is an arbitrary sequence of loose arrows,
    it uniquely factors through the cell $\varpi'_{\alpha,\beta}$ as on the right below.
    \[
        \begin{tikzcd}[virtual, column sep=small]
            I
            \sard[rr, "\ol\gamma"]
            \ar[d, equal]
            \ar[phantom,rrrd, "\nu"]
            &&
            J
            \sar[r, "\alpha"]
            &
            K 
            \ar[d, equal]
            \\
            I
            \sar[rrr, "\beta"']
            &&&
            K
        \end{tikzcd}
        =
        \begin{tikzcd}[virtual, column sep=small]
            I
            \sard[rr, "\ol\gamma"]
            \ar[d, equal]
            \ar[drr, phantom, description, "\wt\nu"]
            &&
            J
            \ar[d, equal]
            \ar[dr,phantom, "{\rotatebox{90}{=}}"]
            \sar[r, "\alpha"]
            &
            K
            \ar[d, equal]
            \\
            I
            \sar[rr, "{\beta\triangleleft\alpha}"']
            \ar[d, equal]
            \ar[phantom,rrrd, "\varpi'_{\alpha,\beta}", yshift=-1ex]
            &&
            J
            \sar[r, "\alpha"']
            &
            K
            \ar[d, equal]
            \\
            I
            \sar[rrr, "\beta"']
            &&&
            K
        \end{tikzcd}
    \]
\end{definition}
With this notion, one can handle the concept of weighted limits and colimits internally in virtual double categories.
We now introduce the \emph{filler protypes} $\triangleright$ and $\triangleleft$ to \ac{FVDblTT} 
to express the right extension and the right lift in the type theory.
The formation rule for the \emph{right extension protype} is the following:
\[
    \inferrule*
    {\syn{w}:\syn{I}\smcl\syn{x}:\syn{J}\vdash \syn{\alpha}(\syn{w}\smcl\syn{x})\ \textsf{protype} \\
    \syn{w}:\syn{I}\smcl\syn{y}:\syn{K}\vdash \syn{\beta}(\syn{w}\smcl\syn{y})\ \textsf{protype}}
    {\syn{x}:\syn{J}\smcl\syn{y}:\syn{K}\vdash \syn{\alpha}(\syn{w}\smcl\syn{x})\triangleright_{\syn{w}:\syn{I}}\syn{\beta}(\syn{w}\smcl\syn{y})\ \textsf{protype}}
\]
The constructor for the right extension protype is given in the elimination rule 
since the orientation of the universal property of the right extension is 
opposite to that of the composition protype and the path protype.
\[
    \inferrule*
    {\syn{w}:\syn{I}\smcl\syn{x}:\syn{J}\vdash \syn{\alpha}(\syn{w}\smcl\syn{x})\ \textsf{protype} \\
    \syn{w}:\syn{I}\smcl\syn{y}:\syn{K}\vdash \syn{\beta}(\syn{w}\smcl\syn{y})\ \textsf{protype}}
    {\syn{w}:\syn{I}\smcl\syn{x}:\syn{J}\smcl\syn{y}:\syn{K}\mid \syn{a}:\syn{\alpha}(\syn{w}\smcl\syn{x})\smcl \syn{e}:\syn{\alpha}(\syn{w}\smcl\syn{x})\triangleright_{\syn{w}:\syn{I}}\syn{\beta}(\syn{w}\smcl\syn{y})\vdash \syn{a}\rbl\syn{e}:\syn{\beta}(\syn{w}\smcl\syn{y})}
\]

The semantics of the right extension protype $\triangleright$ is given by the right extension in \acp{VDC}.
The constructor $\rbl$ is interpreted using the cell $\varpi_{\sem{\syn{\alpha}},\sem{\syn{\beta}}}$ above.
To illustrate the semantics of the right extension protype,
we give the interpretations of the right extension protype in the \acp{VDC} $\Prof{}$ and $\Rel{}$.
\[
    \begin{aligned}
        \sem{\syn{\alpha}(\syn{w}\smcl\syn{x})\triangleright_{\syn{w}:\syn{I}}\syn{\beta}(\syn{w}\smcl\syn{y})}&=
        \int_{W\in\sem{\syn{I}}}\left[\sem{\syn{\alpha}}(W,-),\sem{\syn{\beta}}(W,\bullet)\right]\colon\sem{\syn{J}}\sto\sem{\syn{K}}\  \text{in}\ \Prof{}\\
        \sem{\syn{\alpha}(\syn{w}\smcl\syn{x})\triangleright_{\syn{w}:\syn{I}}\syn{\beta}(\syn{w}\smcl\syn{y})}&=
        \left\{\,(x,y)\mid\forall w\in\sem{\syn{I}}.\left(\sem{\syn{\alpha}}(w,x)\Rightarrow\sem{\syn{\beta}}(w,y)\right)\,\right\}\colon\sem{\syn{J}}\sto\sem{\syn{K}}\ \text{in}\ \Rel{}
    \end{aligned}
\]
Here, $[X,Y]$ is the function set from $X$ to $Y$.

\paragraph*{Comprehension type $\cmpr{\mhyphen}$ for the tabulators.\ ({\Cref{sec:comprehension-type,sec:comprehension-type-meets-unit-protype}})}
The last one is not a protype but a type constructor.
A relation $R\colon A\sto B$ in a general category can be seen as a subobject of the product $A\times B$,
and its legs to $A$ and $B$ give a triangle cell 
\[
    \begin{tikzcd}[virtual, column sep=small]
        &
        X_R
        \ar[dl, "\pi_1"']
        \ar[dr, "\pi_2"]
        &
        \\
        A
        \sar[rr, "R"']
        \ar[rr, phantom, "\tau_R", yshift=2ex]
        &&
        B
    \end{tikzcd}
    \quad\text{where}\quad
    X_R=\left\{\,(a,b)\in A\times B\mid R(a,b)\,\right\}.
\]
This triangle cell is a universal triangle cell whose base is $R$.
In a general virtual double category, such a universal object is called a \emph{tabulator} of a loose arrow $A\sto B$.

\begin{definition}[Tabulators\cite{grandisLimitsDoubleCategories1999}]
    \label{def:tabulators}
    A \emph{(1-dimensional) tabulator} of a loose arrow $\alpha\colon I\sto J$ is an object $\cmpr{\alpha}$ equipped
    with a pair of tight arrows $\ell_{\alpha}\colon \{\alpha\}\to I$ and $r_{\alpha}\colon \{\alpha\}\to J$
    and a cell
    \[
        \begin{tikzcd}[virtual, column sep=small]
            &
            \{\alpha\}
            \ar[dl, "\ell_{\alpha}"']
            \ar[dr, "r_{\alpha}"]
            &
            \\
            I
            \sar[rr, "\alpha"']
            \ar[rr, phantom, "\tau_\alpha", yshift=2ex]
            &&
            J
        \end{tikzcd}
    \]
    such that,
    for any cell $\nu$ on the left below, there exists a unique tight arrow $t_{\nu}\colon X\to \cmpr{\alpha}$ 
    that makes the following two cells equal.
    \[
        \begin{tikzcd}[virtual, column sep=small]
            &
            X
            \ar[dl, "h"']
            \ar[dr, "k"]
            &
            \\
            I
            \sar[rr, "\alpha"']
            \ar[rr, phantom, "\nu", yshift=2ex]
            &&
            J
        \end{tikzcd}
        =
        \begin{tikzcd}[virtual, column sep=small]
            &
            X 
            \ar[ddl, "h"', bend right=30]
            \ar[ddl, phantom, "\circlearrowleft"]
            \ar[ddr, "k", bend left=30]
            \ar[ddr, phantom, "\circlearrowleft"]
            \ar[d, "t_{\nu}"]
            \\
            &
            \{\alpha\}
            \ar[dl, "\ell_{\alpha}"']
            \ar[dr, "r_{\alpha}"]
            &
            \\
            I
            \sar[rr, "\alpha"']
            \ar[rr, phantom, "\tau_\alpha", yshift=2ex]
            &&
            J
        \end{tikzcd}
    \]
    Henceforth, we call a dataset $(X,h,k,\nu)$ a \emph{cone} over $\alpha$ with the apex $X$.
\end{definition}

Corresponding to the tabulators in the virtual double categories,
we introduce the \emph{comprehension type} $\cmpr{\mhyphen}$ to \ac{FVDblTT}.
The formation rule for the comprehension type is the following:
\[
    \inferrule*
    {\syn{x}:\syn{I}\smcl\syn{y}:\syn{J}\vdash \syn{\alpha}\ \textsf{protype}}
    {\cmpr{\syn{\alpha}} \ \textsf{type}}
\]
This comes equipped with the constructor
\[
    \inferrule*
    {\syn{x}:\syn{I}\smcl\syn{y}:\syn{J}\vdash \syn{\alpha}\ \textsf{protype}}
    {\syn{w}:\cmpr{\syn{\alpha}}\vdash \syn{l}(\syn{w}):\syn{I}\\
    \syn{w}:\cmpr{\syn{\alpha}}\vdash \syn{r}(\syn{w}):\syn{J}\\
    \syn{w}:\cmpr{\syn{\alpha}}\mid \vdash \tabb_{\cmpr{\syn{\alpha}}}(\syn{w}):\syn{\alpha}[\syn{l}(\syn{w})/\syn{x}\smcl\syn{r}(\syn{w})/\syn{y}]}
\]

The comprehension type $\cmpr{\mhyphen}$ is interpreted as the tabulators in the \acp{VDC}.
In the \ac{VDC} $\Prof$,
the tabulator of a profunctor $P\colon\one{C}\sto\one{D}$ is given 
by \emph{two-sided Grothendieck construction},
which results in \emph{a two-sided discrete fibration} from $\one{C}$ to $\one{D}$.
A frequently used example of this construction is
the comma category for a pair of functors $F\colon\one{C}\to\one{E}$ and $G\colon\one{D}\to\one{E}$
as the tabulator of the profunctor $\one{E}(F(-),G(-))$,
see \cite{loregianCategoricalNotionsFibration2020} for more details.
The \ac{VDC} $\Rel{}$ has the tabulators
if we ground the double category to an axiomatic system of set theory with the comprehension axiom,
as the tabulator of a relation $R\colon A\sto B$ is given by the set of all the pairs $(a,b)$ such that $R(a,b)$ holds.

In the presence of the unit protype $\ide{}$, 
we should add some rules concerning the compatibility between the comprehension type and the path protype.
This is because, in many examples of double categories, the tabulators have not only the universal property
as in \Cref{def:tabulators} but also respect the units,
although the original universal property of the tabulators is enough to detect the tabulators in a double category.
This issue is thoroughly discussed in \cite{grandisLimitsDoubleCategories1999}.
Here, we give a slightly generalized version of the tabulators in virtual double categories with units.
\begin{definition}[2-dimensional universal property of tabulators]
    In a virtual double category with units, an \emph{unital tabulator} $\{\alpha\}$ of a loose arrow $\alpha\colon I\sto J$
    is a tabulator of $\alpha$ in the sense of \Cref{def:tabulators}, which also satisfies the following universal property.
    Suppose we are given any pair of cones $(X,h,k,\nu)$ and $(X',h',k',\nu')$ over $\alpha$
    and a pair of cells $\varsigma,\vartheta$ such that the following equality holds.
    \[
        \begin{tikzcd}[virtual]
            X
            \sard[r, "\ol\gamma"]
            \ar[rd, phantom, "\varsigma"]
            \ar[d, "h"']
            &
            X'
            \ar[d, "h'"']
            \ar[dr, phantom, "\nu'", yshift=-1ex, xshift=-2ex]
            \ar[dr,"k'"]
            \\
            I
            \sar[r, "U_I"']
            \ar[d, equal]
            \ar[drr, phantom, "{\rotatebox{90}{$\cong$}}"]
            &
            I
            \sar[r, "\alpha"']
            &
            J
            \ar[d, equal]
            \\
            I 
            \sar[rr, "\alpha"']
            &&
            J
        \end{tikzcd}
        =
        \begin{tikzcd}[virtual]
            &
            X
            \ar[dl, "h"']
            \ar[d, "k"]
            \sard[r, "\ol\gamma"]
            \ar[dr, phantom, "\vartheta"]
            \ar[dl, phantom, "\nu", yshift=-1ex, xshift=1ex]
            &
            X'
            \ar[d, "k'"]
            \\
            I
            \sar[r, "\alpha"']
            \ar[d, equal]
            \ar[rrd, phantom, "{\rotatebox{90}{$\cong$}}"]
            &
            J
            \sar[r, "U_J"']
            &
            J
            \ar[d, equal]
            \\
            I
            \sar[rr, "\alpha"']
            &&
            J
        \end{tikzcd}
    \]
    Then, there exists a unique cell $\varrho$ for which the following equalities hold.
    \[
        \begin{tikzcd}[virtual]
            X
            \sard[r, "\ol\gamma"]
            \ar[d, "t_{\nu}"']
            \ar[rd, phantom, "\varrho"]
            &
            X'
            \ar[d, "t_{\nu'}"]
            \\
            \{\alpha\}
            \sar[r, "U_{\{\alpha\}}"']
            \ar[d, "\ell_{\alpha}"']
            \ar[dr, phantom, "U_{\ell_\alpha}"description, yshift=-1ex]
            &
            \{\alpha\}
            \ar[d, "\ell_{\alpha}"]
            \\
            I
            \sar[r, "U_I"']
            &
            I
        \end{tikzcd}
        =
        \begin{tikzcd}[virtual]
            X
            \sard[r, "\ol\gamma"]
            \ar[d, "h"']
            \ar[rd, phantom, "\varsigma"]
            &
            X'
            \ar[d, "h'"]
            \\
            I
            \sar[r, "U_I"']
            &
            I
        \end{tikzcd}
        \quad
        ,
        \quad
        \begin{tikzcd}[virtual]
            X
            \sard[r, "\ol\gamma"]
            \ar[d, "t_{\nu}"']
            \ar[rd, phantom, "\varrho"]
            &
            X'
            \ar[d, "t_{\nu'}"]
            \\
            \{\alpha\}
            \sar[r, "U_{\{\alpha\}}"']
            \ar[d, "r_{\alpha}"']
            \ar[dr, phantom, "U_{r_\alpha}"description, yshift=-1ex]
            &
            \{\alpha\}
            \ar[d, "r_{\alpha}"]
            \\
            J 
            \sar[r, "U_J"']
            &
            J
        \end{tikzcd}
        =
        \begin{tikzcd}[virtual]
            X
            \sard[r, "\ol\gamma"]
            \ar[d, "k"']
            \ar[rd, phantom, "\vartheta"]
            &
            X'
            \ar[d, "k'"]
            \\
            J
            \sar[r, "U_J"']
            &
            J
        \end{tikzcd}
    \]
\end{definition}
This universal property determines what the unit on the apex of the tabulator should be. 
\Cref{sec:comprehension-type-meets-unit-protype} will present the corresponding rules for the comprehension type $\cmpr{\mhyphen}$ in \ac{FVDblTT} with the unit protype $\ides$.

\begin{remark}[Substitution into the additional constructor]
    \label{rem:substitution-into-additional-constructor}
    There are options how we define substitution for the additional constructors.
    For example, we may define the substitution for the composition protype as follows.
    \[
    (\syn{\alpha}\odot_{\syn{y}:\syn{J}}\syn{\beta})[\syn{s}/\syn{x}\smcl\syn{t}/\syn{z}]
    \coloneqq
    \syn{\alpha}[\syn{s}/\syn{x}\smcl\syn{y}/\syn{y}]\odot_{\syn{y}:\syn{J}}\syn{\beta}[\syn{y}/\syn{y}\smcl\syn{t}/\syn{z}] 
    \]
    This seems reasonable for our use in formal category theory,
    but this equality is not always satisfied in a general PL-composable \ac{FVDC}
    unless it is actually a virtual equipment.
    Instead, we may extend the introduction rule for the composition protype 
    so that the substituted composition protypes are directly introduced.
    \[
    \inferrule*
    {\syn{w}:\syn{I}\smcl\syn{x}:\syn{J}\vdash \syn{\alpha}(\syn{w}\smcl\syn{x})\ \textsf{protype} \\
    \syn{x}:\syn{J}\smcl\syn{y}:\syn{K}\vdash \syn{\beta}(\syn{x}\smcl\syn{y})\ \textsf{protype} \\
    \syn{\Gamma}\vdash \syn{s}:\syn{I} \\
    \syn{\Delta}\vdash \syn{t}:\syn{K}}
    {\syn{\Gamma}\smcl\syn{\Delta}\vdash (\syn{\alpha}\odot_{\syn{x}:\syn{J}}\syn{\beta})[\syn{s}\smcl\syn{t}/\syn{x}\smcl\syn{y}]: \textsf{protype}}
    \]
    Then, the substitution for the composition protype is obvious.
    Indeed we take the latter approach for the path protype.

    Therefore, it depends on the purpose of the type theory how we define the substitution for the additional constructors,
    and we do not specify it in this paper because our main focus is the syntax-semantics duality
    for the very basic type theory.
\end{remark}

%% file: newsubfiles/constructor/predlogic.tex
\paragraph*{Predicate logic.}
\label{sec:arrangement-for-predicate-logic}

When we work with the type theory \ac{FVDblTT} for the purpose of reasoning about predicate logic,
we consider types, terms, protypes, and proterms to represent sets, functions, predicates (or propositions), and proofs, respectively. 
However, the type theory \ac{FVDblTT}, as it is, treats the protypes in a context $\syn{\Gamma}\smcl\syn{\Delta}$ and
those in a context $\syn{\Delta}\smcl\syn{\Gamma}$ as different things.
In this sense, the type theory \ac{FVDblTT} as predicated logic has directionality.
If one wants to develop a logic without a direction,
one can simply add the following rules to the type theory.

    \begin{mathparpagebreakable}
        \inferrule*
        {\syn{\Gamma}\smcl\syn{\Delta}\vdash \syn{\alpha} \ \textsf{protype}}
        {\syn{\Delta}\smcl\syn{\Gamma}\vdash \syn{\alpha}^\circ \ \textsf{protype}}
        \and 
        \inferrule*
        {\syn{\Gamma}_0\smcl\cdots\smcl\syn{\Gamma}_m\mid \syn{a}_1:\syn{\alpha}_1\cdots \syn{a}_n:\syn{\alpha}_n
        \vdash \syn{\mu}:\syn{\beta}}
        {\syn{\Gamma}_m\smcl\cdots\smcl\syn{\Gamma}_0\mid \syn{a}_n:\syn{\alpha}_n^\circ\cdots \syn{a}_1:\syn{\alpha}_1^\circ
        \vdash \syn{\mu}^\circ:\syn{\beta}^\circ}
        \and
        \inferrule*
        {\syn{\Gamma}_0\smcl\cdots\smcl\syn{\Gamma}_m\mid \syn{a}_1:\syn{\alpha}_1\cdots \syn{a}_n:\syn{\alpha}_n
        \vdash \syn{\mu}:\syn{\beta}}
        {\syn{\Gamma}_0\smcl\cdots\smcl\syn{\Gamma}_m\mid \syn{a}_1:\syn{\alpha}_1\cdots \syn{a}_n:\syn{\alpha}_n
        \vdash \syn{\mu}^{\circ\circ}\equiv\syn{\mu}:\syn{\beta}}
    \end{mathparpagebreakable}
These rules are the counterparts of the structure of involution in \acp{VDC}.

If one also wants to make the type theory \ac{FVDblTT} proof irrelevant,
one can reformulate protype isomorphism judgment as equality judgments of protypes 
and add the rule stating that all the proterms are equal.
It is the counterpart of the flatness \cite{grandisLimitsDoubleCategories1999} or local preorderedness
\cite{hoshinoDoubleCategoriesRelations2023} of \acp{VDC}. 

%% file: newsubfiles/constructor/examples.tex
This section exemplifies how one can reason about category theory and logic formally in the type theory \ac{FVDblTT}.

\begin{example}[(co)Yoneda Lemma]
    \label{example:ninja-yoneda}
    One of the most fundamental results in category theory is the Yoneda Lemma,
    and it has a variety of presentations in the literature.
    Here we present one called the Yoneda Lemma \cite[Proposition 2.2.1]{loregianCoEndCalculus2021}:
    \textit{given a category $\one{C}$ and a functor $F\colon\one{C}\op\to\Set$, we have the canonical isomorphism}
    \[
        F\cong\int_{X\in\one{C}}[\one{C}(X,-),FX] .
    \]
    This follows from the categorical fact that $\Prof$ is an \ac{FVDC} with the structures listed above. 
    Indeed, in the type theory \ac{FVDblTT} with the path protype $\ides$ and the filler protype $\triangleright$,
    one can deduce the following:
    \[
        \syn{y}:\syn{I}\smcl\cdot\vdash \textsf{Yoneda}:\left(\syn{x}\ide{\syn{I}}\syn{y}\right)\triangleright_{\syn{x}:\syn{I}}\syn{\alpha}(\syn{x})\ccong\syn{\alpha}(\syn{y})
    \]
    Similarly, we have
    \[
        \syn{y}:\syn{I}\smcl\cdot\vdash \textsf{CoYoneda}:\left(\syn{y}\ide{\syn{I}}\syn{x}\right)\odot_{\syn{x}:\syn{I}}\syn{\alpha}(\syn{x})\ccong\syn{\alpha}(\syn{y})
    \]
    which expresses
    the coYoneda Lemma:
    \[
        \int^{X\in\one{C}}\one{C}(-,X)\times FX\cong F.
    \]
    In short, all the theorems in category theory that can be proven using this type theory 
    fall into corollaries of the theorem that $\Prof$ is a \ac{CFVDC} with the structures corresponding to the constructors.
    Other examples include the unit laws and the associativity of the composition of profunctors or the iteration 
    of extensions and lifts of profunctors.

    Turning to the aspect of predicate logic, we can interpret the protype isomorphisms as the following logical equivalences.
    \[
        \begin{aligned}
            \varphi(y)\quad&\equiv\quad\forall x\in I.\left(x=y\right)\Rightarrow\varphi(x)\\
            \varphi(y)\quad&\equiv\quad\exists x\in I.\left(x=y\right)\land\varphi(x)
        \end{aligned}
    \]
\end{example}

\begin{example}[Isomorphism of functors]
    A natural transformation $\xi\colon F\to G$ between two functors $F,G\colon\one{C}\to\one{D}$ is given by 
    a family of arrows $\xi_X\colon FX\to GX$ satisfying some naturality conditions.
    In the type theory \ac{FVDblTT} with the path protype $\ides$,
    this natural transformation can be represented by a proterm $\syn{x}:\syn{I}\mid\vdash \syn{\xi}(\syn{x}):\syn{f}(\syn{x})\ide{\syn{I}}\syn{g}(\syn{x})$.
    Here, the naturality condition automatically holds because we describe it as a proterm.
    The isomorphism of functors can be expressed using this notion, 
    but an alternative way is to use the protype isomorphism. 

    \begin{lemma}
        \label{lemma:isomorphism-of-functors}
        Given two terms, $\syn{f}(\syn{x})$ and $\syn{g}(\syn{x})$, in the same context,
        the following are equivalent.
        \begin{enumerate}
            \labeleditem There are proterms $\syn{\xi}(\syn{x}):\syn{f}(\syn{x})\ide{\syn{I}}\syn{g}(\syn{x})$ and $\syn{\eta}(\syn{x}):\syn{g}(\syn{x})\ide{\syn{I}}\syn{f}(\syn{x})$
            such that $\syn{\xi}(\syn{x})\boxdot\syn{\eta}(\syn{x})\equiv\refl_{\syn{f}(\syn{x})}$ and $\syn{\eta}(\syn{x})\boxdot\syn{\xi}(\syn{x})\equiv\refl_{\syn{g}(\syn{x})}$.
            \label{lemma:isomorphism-of-functors:1}
            \labeleditem There is a protype isomorphism $\syn{Z}:\syn{y}\ide{\syn{J}}\syn{f}(\syn{x})\ccong\syn{y}\ide{\syn{I}}\syn{g}(\syn{x})$. 
            \label{lemma:isomorphism-of-functors:2}
        \end{enumerate}
        Here, $\boxdot$ is a tailored constructor
        defined as follows.
        \[
            \syn{y}:\syn{J}\smcl\syn{y}':\syn{J}\smcl\syn{y}'':\syn{J}\mid\
            \syn{a}:\syn{y}\ide{\syn{J}}\syn{y}'\smcl\syn{b}:\syn{y}'\ide{\syn{J}}\syn{y}''\vdash
            \syn{a}\boxdot\syn{b}\colequiv \ind_{\ide{\syn{J}}}(\syn{a}):\syn{y}\ide{\syn{J}}\syn{y}''.
        \]
    \end{lemma}
    \begin{proof}
        First, suppose \Cref{lemma:isomorphism-of-functors:1} holds.
        We define a proterm $\syn{\zeta}$ by the following:
        \[
        \small
        \inferrule*
        {
        \syn{x}:\syn{I}\mid\vdash \syn{\xi}: \syn{f}(\syn{x})\ide{\syn{J}}\syn{g}(\syn{x})\\
        \inferrule*
        {
        \syn{y}:\syn{J}\smcl\syn{y'}:\syn{J}\smcl\syn{y''}:\syn{J}\mid
        \syn{a}:\syn{y}\ide{\syn{J}}\syn{y'}\smcl\syn{b}:\syn{y'}\ide{\syn{J}}\syn{y''}
        \vdash \syn{a}\boxdot\syn{b}:\syn{y}\ide{\syn{J}}\syn{y''}
        }
        {
        \syn{y}:\syn{J}\smcl\syn{x}:\syn{I}\smcl\syn{x'}:\syn{I}\mid
        \syn{a}:\syn{y}\ide{\syn{J}}\syn{f(x)}\smcl\syn{b}:\syn{f(x)}\ide{\syn{J}}\syn{g(x)}
        \vdash \syn{a}\boxdot\syn{b}[\syn{y}/\syn{y}\smcl\syn{f(x)}/\syn{y'}\smcl\syn{g(x)}/\syn{y''}]:\syn{y}\ide{\syn{J}}\syn{g(x)}
        }
        }
        {\syn{y}:\syn{J}\smcl\syn{x}:\syn{I}\mid \syn{a}:\syn{y}\ide{\syn{J}}\syn{f}(\syn{x})\vdash \syn{\zeta}(\syn{a}):
        \syn{y}\ide{\syn{J}}\syn{g}(\syn{x})}
        \]
        Therefore, we have
        $\syn{\zeta}(\syn{a})$,
        and in the same way, we can define a proterm $\syn{b}:\syn{y}\ide{\syn{J}}\syn{g}(\syn{x})\vdash\syn{\zeta'}(\syn{b}):
        \syn{y}\ide{\syn{J}}\syn{f}(\syn{x})$,
        which is the inverse of $\syn{\zeta}$ by simple reasoning.
    
        Next, suppose \Cref{lemma:isomorphism-of-functors:2} holds.
        Let $\syn{a}:\syn{y}\ide{\syn{J}}\syn{f}(\syn{x})\vdash \syn{\zeta}(\syn{a}):
        \syn{y}\ide{\syn{J}}\syn{g}(\syn{x})$ be the proterm witnessing the isomorphism.
        By substituting $\syn{f}(\syn{x})$ for $\syn{y}$ and the $\refl$ for $\syn{a}$,
        we obtain a proterm $\syn{\xi}(\syn{x}):\syn{f}(\syn{x})\ide{\syn{J}}\syn{g}(\syn{x})$.
        In the same way, we can define a proterm $\syn{\eta}(\syn{x}):\syn{g}(\syn{x})\ide{\syn{J}}\syn{f}(\syn{x})$,
        for which the two desired equalities hold.
    \end{proof}
    
    We therefore use the equalities $\syn{y}\ide{\syn{J}}\syn{f}(\syn{x})$ and $\syn{y}\ide{\syn{J}}\syn{g}(\syn{x})$ 
    when $\syn{f}$ and $\syn{g}$ are already proven to be isomorphic. 
\end{example}

\begin{example}[Adjunction]
    In a virtual double category, the \emph{companion} and \emph{conjoint} of a tight arrow $f\colon A\to B$ is defined 
    as the loose arrows $f_*\colon A\sto B$ and $f^*\colon B\sto A$ equipped with cells satisfying some equations of cells \cite{grandisAdjointDoubleCategories2004,cruttwellUnifiedFrameworkGeneralized2010}.
    In a virtual equipment,
    it is known that the companion and conjoint of a tight arrow $f\colon A\to B$ are
    the restrictions of the units on $B$ along the pairs of tight arrows $(f,\id_B)$ and $(\id_B,f)$, respectively.
    These notions are the formalization of the representable profunctors in the virtual double categories.
    Therefore, the companions and conjoints of a term $\syn{t}(\syn{x})$ in the type theory \ac{FVDblTT} should be
    defined as $\syn{t}(\syn{x})\ide{\syn{I}}\syn{y}$ and $\syn{y}\ide{\syn{I}}\syn{t}(\syn{x})$, respectively.

    The \emph{adjunction} between two functors is described in terms of representable profunctors,
    which motivates the following definition of the adjunction in the type theory \ac{FVDblTT}.
    Remember a functor $F\colon\one{C}\to\one{D}$ is left adjoint to a functor $G\colon\one{D}\to\one{C}$ 
    if there is a natural isomorphism between the hom-sets
    \[
        \one{D}(F-,\bullet)\cong\one{C}(-,G\bullet).
    \]
    In the type theory \ac{FVDblTT}, a term $\syn{t}(\syn{x})$ is announced to be a left adjoint to a term $\syn{u}(\syn{y})$ 
    if the following equality holds:
    \[
        \syn{x}:\syn{I}\smcl\syn{y}:\syn{J}\vdash \syn{t}(\syn{x})\ide{\syn{J}}\syn{y}\equiv \syn{x}\ide{\syn{I}}\syn{u}(\syn{y}).
    \]
\end{example}

\begin{example}[Kan extension]
    In \cite{kellyBasicConceptsEnriched2005},
    the (pointwise) left Kan extension $\Lan_{G}F$ of a functor $F\colon\one{C}\to\one{D}$ along a functor $G\colon\one{C}\to\one{E}$ is defined as a functor $H\colon\one{D}\to\one{E}$
    equipped with a natural transformation 
    \[
        \begin{tikzcd}[column sep=small]
            \one{C}
            \ar[rr, "F"]
            \ar[dr, "G"']
            &\!&
            \one{D}
            \ar[from=dl, "H"']
            \\
            &
            \one{E}
            \ar[u, phantom, "{\rotatebox{270}{$\Rightarrow$}}"]
            \ar[u, phantom, "\mu"{xshift=-1.5ex,yshift=1ex}]
        \end{tikzcd}
    \]
    with the following canonical natural transformation being an isomorphism:
    \[
        \one{D}(HE,D)\overset{\cong}{\to}\widehat{\one{C}}\left(\one{E}(G-,E),\one{D}(F-,D)\right) \quad \text{naturally in}\ D\in\one{D}, E\in\one{E}.
    \]
    A protype isomorphism corresponding to this isomorphism is given by the following.
    \[
            \syn{z}:\syn{K}\smcl\syn{y}:\syn{J}\vdash \textsf{LeftKan}: \syn{h}(\syn{z})\ide{\syn{J}}\syn{y}\ccong
            \left(\syn{g}(\syn{x})\ide{\syn{K}}\syn{z}\right)\triangleright_{\syn{x}:\syn{I}}
            \left(\syn{f}(\syn{x})\ide{\syn{J}}\syn{y}\right)
    \]
    We will demonstrate how proofs in category theory can be done in the type theory \ac{FVDblTT}.
    \begin{proposition}[{\cite[Theorem 4.47]{kellyBasicConceptsEnriched2005}}]
        $\Lan_{G'}\Lan_GF\cong\Lan_{G'\circ G}F$ hold
        for any functors $F\colon\one{C}\to\one{D}$, $G\colon\one{C}\to\one{E}$, and $G'\colon\one{E}\to\one{F}$
        if the Kan extensions exist.
        \[
            \begin{tikzcd}[column sep=small, row sep=small]
                \one{C}
                \ar[rrr, "F"]
                \ar[dr, "G"']
                &\!&\!&
                \one{D}
                \ar[from=dll, "\Lan_{G}F"description]
                \ar[from=ddl, "\Lan_{G'}\Lan_{G}F\cong\Lan_{G'\circ G}F"']
                \\
                &
                \one{E}
                \ar[u, phantom, "{\rotatebox{270}{$\Rightarrow$}}"]
                \ar[dr, "G'"']
                \\
                &&
                \one{E}'
                \ar[uu, phantom, "{\rotatebox{270}{$\Rightarrow$}}",yshift=-2ex]
            \end{tikzcd}
        \]
    \end{proposition}
    \begin{proof}
        We associate $F,G,G',\Lan_{G}F,\Lan_{G'}\Lan_GF,\Lan_{G'\circ G}F$ 
        with the terms $\syn{f}(\syn{x}),\syn{g}(\syn{x}),\syn{g'}(\syn{z}),\syn{h}(\syn{z}),\syn{h'}(\syn{z}')$, and $\syn{h''}(\syn{z}')$.
        We will have the desired protype isomorphism judgment
        by composing the protype isomorphisms in the following order.
        {\small
        \begin{align*}
            \syn{z'}:\syn{K'}\smcl\syn{y}:\syn{J}&\mid
            \syn{h'}(\syn{z'})\ide{\syn{J}}\syn{y}\\
            &\ccong 
            \left(\syn{g}'(\syn{z})\ide{\syn{K'}}\syn{z'}\right)\triangleright_{\syn{z}:\syn{K}}
            \left(\syn{h}(\syn{z})\ide{\syn{J}}\syn{y}\right)
            &(\textsf{LeftKan})\\
            &\ccong 
            \left(\syn{g}'(\syn{z})\ide{\syn{K'}}\syn{z'}\right)
            \triangleright_{\syn{z}:\syn{K}}
            \left(\left(\syn{g}(\syn{x})\ide{\syn{K}}\syn{z}\right)\triangleright_{\syn{x}:\syn{I}}
            \left(\syn{f}(\syn{x})\ide{\syn{J}}\syn{y}\right)\right)
            &(\,\left(\syn{g}'(\syn{z})\ide{\syn{K'}}\syn{z'}\right)\triangleright_{\syn{z}:\syn{K}}\textsf{LeftKan}\,)\\
            &\ccong
            \left(\left(\syn{g}(\syn{x})\ide{\syn{K}}\syn{z}\right)\odot_{\syn{z}:\syn{K}}
            \left(\syn{g}'(\syn{z})\ide{\syn{K'}}\syn{z'}\right)\right)
            \triangleright_{\syn{x}:\syn{I}}
            \left(\syn{f}(\syn{x})\ide{\syn{J}}\syn{y}\right)
            &(\textsf{Fubini})\\
            &\ccong
            \left(\syn{g'}(\syn{g}(\syn{x}))\ide{\syn{K'}}\syn{z'}\right)\triangleright_{\syn{x}:\syn{I}}
            \left(\syn{f}(\syn{x})\ide{\syn{J}}\syn{y}\right)
            &(\,\textsf{CoYoneda}\triangleright_{\syn{x}:\syn{I}}\left(\syn{f}(\syn{x})\ide{\syn{J}}\syn{y}\right)\,)\\
            &\ccong
            \syn{h''}(\syn{z'})\ide{\syn{K'}}\syn{y} &(\textsf{LeftKan}\inv)
        \end{align*}
        }
        Here, the protype isomorphism $\textsf{Fubini}$ is given as $\lcp\textsf{Fubini}_1,\textsf{Fubini}_2\rcp$,
        where $\textsf{Fubini}_1$ and $\textsf{Fubini}_2$ are the proterms derived as follows.
            \begin{mathparpagebreakable}
            \small
            \inferrule*{
            \inferrule*{
            \syn{x}_0:\syn{I}_0\smcl\syn{x}_1:\syn{I}_1\smcl\syn{x}_2:\syn{I}_2\smcl\syn{x}_3:\syn{I}_3\mid 
            \syn{a}:\syn{\alpha}\smcl\syn{b}:\syn{\beta}\smcl\syn{c}:\syn{\beta}\triangleright_{\syn{x}_1:\syn{I}_1}\left(\syn{\alpha}\triangleright_{\syn{x}_0:\syn{I}_0}\syn{\gamma}\right)
            \vdash \syn{a}\rbl(\syn{b}\rbl\syn{c}): \syn{\gamma}
            }
            {
            \syn{x}_0:\syn{I}_0\smcl\syn{x}_2:\syn{I}_2\smcl\syn{x}_3:\syn{I}_3\mid
            \syn{d}:\syn{\alpha}\odot_{\syn{x}_1:\syn{I}_1}\syn{\beta}\smcl\syn{c}:\syn{\beta}\triangleright_{\syn{x}_1:\syn{I}_1}\left(\syn{\alpha}\triangleright_{\syn{x}_0:\syn{I}_0}\syn{\gamma}\right)
            \vdash \_: \syn{\gamma}
            }
            }
            {
                \syn{x}_2:\syn{I}_2\smcl\syn{x}_3:\syn{I}_3\mid
                \syn{c}:\syn{\beta}\triangleright_{\syn{x}_1:\syn{I}_1}\left(\syn{\alpha}\triangleright_{\syn{x}_0:\syn{I}_0}\syn{\gamma}\right)
                \vdash \textsf{Fubini}_1: \left(\syn{\alpha}\odot_{\syn{x}_1:\syn{I}_1}\syn{\beta}\right)\triangleright_{\syn{x}_0:\syn{I}_0}\syn{\gamma}
            }
            \end{mathparpagebreakable}
            \begin{mathparpagebreakable}
            \small
            \inferrule*{
            \inferrule*{
            \inferrule*{
                \syn{x}_0:\syn{I}_0\smcl\syn{x}_1:\syn{I}_1\smcl\syn{x}_2:\syn{I}_2 
                \mid \syn{a}:\syn{\alpha}\smcl\syn{b}:\syn{\beta}\vdash \syn{a}\odot\syn{b}:\syn{\alpha}\odot_{\syn{x}_1:\syn{I}_1}\syn{\beta}\\
                \syn{x}_0:\syn{I}_0\smcl\syn{x}_2:\syn{I}_2\smcl\syn{x}_3:\syn{I}_3\mid 
                \syn{d}:\syn{\alpha}\odot_{\syn{x}_1:\syn{I}_1}:\syn{\beta}\smcl\syn{e}:\left(\syn{\alpha}\odot_{\syn{x}_1:\syn{I}_1}\syn{\beta}\right)\triangleright_{\syn{x}_0:\syn{I}_0}\syn{\gamma}
                \vdash \syn{d}\rbl\syn{e}: \syn{\gamma}
            }
            {
                \syn{x}_0:\syn{I}_0\smcl\syn{x}_1:\syn{I}_1\smcl\syn{x}_2:\syn{I}_2\smcl\syn{x}_3:\syn{I}_3\mid 
                \syn{a}:\syn{\alpha}\smcl\syn{b}:\syn{\beta}\smcl\syn{e}:\left(\syn{\alpha}\odot_{\syn{x}_1:\syn{I}_1}\syn{\beta}\right)\triangleright_{\syn{x}_0:\syn{I}_0}\syn{\gamma}
                \vdash \_ : \syn{\gamma}
            }
            }
            {
            \syn{x}_1:\syn{I}_1\smcl\syn{x}_2:\syn{I}_2\smcl\syn{x}_3:\syn{I}_3\mid
            \syn{b}:\syn{\beta}\smcl\syn{e}:\left(\syn{\alpha}\odot_{\syn{x}_1:\syn{I}_1}\syn{\beta}\right)\triangleright_{\syn{x}_0:\syn{I}_0}\syn{\gamma}
            \vdash \_: \syn{\alpha}\triangleright_{\syn{x}_0:\syn{I}_0}\syn{\gamma}
            }}
            {
                \syn{x}_2:\syn{I}_2\smcl\syn{x}_3:\syn{I}_3\mid
                \syn{e}: \left(\syn{\alpha}\odot_{\syn{x}_1:\syn{I}_1}\syn{\beta}\right)\triangleright_{\syn{x}_0:\syn{I}_0}\syn{\gamma}
                \vdash \textsf{Fubini}_2: \syn{\beta}\triangleright_{\syn{x}_1:\syn{I}_1}\left(\syn{\alpha}\triangleright_{\syn{x}_0:\syn{I}_0}\syn{\gamma}\right)
            }
            \end{mathparpagebreakable}
    \end{proof}
\end{example}

%% file: newsubfiles/adjunction/synsem.tex
Stating that a type theory is the internal language of a categorical structure
always comes with the notion of a syntax-semantics adjunction.
We set out to construct the term model of \ac{FVDblTT} by following the standard procedure of categorical logic.

%% file: newsubfiles/adjunction/syntacticpres.tex
Now, we turn to the definition of a specification for a signature in the type theory.

\begin{definition}
    Let $\Phi\colon\Sigma\to\Sigma'$ be a morphism of signatures,
    and $J$ be a judgment in the type theory based on $\Sigma$.
    We write $J^{\Phi}$ for the judgment in $(\Sigma,\zero{E})$ defined 
    by replacing each symbol in $J$ with its image under $\Phi$.
    $J^{\Phi}$ is called the \emph{translation} of $J$ via $\Phi$.
\end{definition}

\begin{definition}
    \label{def:specification}
    A \emph{specification} $\zero{E}$ for a signature $\Sigma$ is a pair $(\zero{E}^{\tm},\zero{E}^{\ptm})$ where
    \begin{itemize}
        \item $\zero{E}^{\tm}$ is a class of pair of terms of the same type that are well-formed in $\Sigma$, 
        \item $\zero{E}^{\ptm}$ is a class of proterm equality judgments that are well-formed in $\Sigma$ and $\zero{E}^{\tm}$.
    \end{itemize}
    When we say $(\Sigma,\zero{E})$ is a specification, we mean that $\Sigma$ is a signature and $\zero{E}$ is a specification for $\Sigma$.

    A \emph{morphism of specifications} $\Phi\colon (\Sigma,\zero{E})\to(\Sigma',\zero{E}')$ is a morphism of signatures $\Phi\colon\Sigma\to\Sigma'$
    by which every judgment in $\zero{E}$ is translated to a judgment that is derivable from $\zero{E}'$.
\end{definition}

\begin{definition}[Validity of equality judgments]
    We define the validity of equality judgments in a \ac{CFVDC} as follows.
    \begin{itemize}
        \item A term equality judgment $\syn{t}\equiv\syn{t'}$ is \emph{valid} in a $\Sigma$-structure $\one{M}$ in a \ac{CFVDC} $\dbl{D}$ if
        $\sem{\syn{t}}_{\one{M}}$ and $\sem{\syn{t'}}_{\one{M}}$ are equal as tight arrows in $\dbl{D}$.
        \item A proterm equality judgment $\syn{\mu}\equiv\syn{\mu'}$ is \emph{valid} in
        a $\Sigma$-structure $\one{M}$ in a \ac{CFVDC} $\dbl{D}$ if
        $\sem{\syn{\mu}}_{\one{M}}$ and $\sem{\syn{\mu'}}_{\one{M}}$ are equal as cells in $\dbl{D}$.
    \end{itemize}
\end{definition}

With the definition of validity, one can canonically associate a specification $\zero{E}_{\dbl{D}}$ to a \ac{CFVDC} $\dbl{D}$,
which exhaustively contains the information of $\dbl{D}$.

\begin{definition}
    \label{def:associatedspec}
    The \emph{associated specification} $\Sp(\dbl{D})$ of a \ac{CFVDC} $\dbl{D}$ is
    the specification $(\Sigma_{\dbl{D}},\zero{E}_{\dbl{D}})$ with $\Sigma_{\dbl{D}}$ as above,
    $\zero{E}_{\dbl{D}}^{\tm}$ (resp. $\zero{E}_{\dbl{D}}^{\ptm}$)
    the set of all the valid equality judgments for terms (resp. proterms) in the canonical structure in $\dbl{D}$.
\end{definition}

%% file: newsubfiles/adjunction/newproof.tex
We will construct a biadjunction between the 2-category of virtual double categories and the 2-category of specifications in \ac{FVDblTT}. 

The first goal is to construct a 1-adjunction between the category of specifications and the category of
split \acp{CFVDC} and morphisms between them.
\begin{definition}
    For a specification $(\Sigma,\zero{E})$,
        the \emph{syntactic virtual double category} (or classifying virtual double category) $\nS(\Sigma,\zero{E})$ is the virtual double category whose 
        \begin{itemize}
            \item objects are contexts $\syn{\Gamma}\ \textsf{ctx}$ in $\Sigma$,
            \item tight arrows $\syn{\Gamma}\to\syn{\Delta}=(\syn{y}_1:\syn{J}_1,\dots,\syn{y}_n:\syn{J}_n)$ are equivalence classes of
            sequences of terms (or, term substitutions) $\syn{\Gamma}\vdash \syn{s}_1:\syn{J}_1,\dots,\syn{s}_n:\syn{J}_n$ (or substitutions)
            modulo equality judgments derivable from $(\Sigma,\zero{E})$,
            \item loose arrows $\syn{\Gamma}\sto\syn{\Delta}$ are 
            protypes $\syn{\Gamma}\smcl\syn{\Delta}\vdash \syn{\alpha}\ \textsf{protype}$ in $\Sigma$
            modulo equality judgments derivable from $(\Sigma,\zero{E})$,
            \item cells of form 
            \begin{equation}
                \label{eq:cell}
                \begin{tikzcd}[column sep=8ex,virtual]
                    {\syn{\Gamma}_0}
                    \sar[r, "{\syn{\alpha}_1}"]
                    \ar[d, "{\syn{S}_0}"']
                    \ar[phantom,rrrd, "{\syn{\mu}}" description]
                    &
                    \cdots
                    &
                    \cdots
                    \sar[r, "{\syn{\alpha}_n}"]
                    &
                    {\syn{\Gamma}_n}
                    \ar[d,"\syn{S}_1"]
                    \\
                    {\syn{\Delta}_0}
                    \sar[rrr, "{\syn{\beta}}"'] 
                    &
                    &&
                    {\syn{\Delta}_1}
                \end{tikzcd}
            \end{equation}
            are equivalence classes of proterms
            \[
            \syn{\ol{\Gamma}}\mid \syn{a}_1:\syn{\alpha}_1\smcl\dots\smcl\syn{a}_n:\syn{\alpha}_n
            \vdash \syn{\mu}:\syn{\beta}[\syn{S}_0/\syn{\Delta}_0\smcl\syn{S}_n/\syn{\Delta}_n]
            \]
            modulo equality judgments derivable from $(\Sigma,\zero{E})$.
            It makes no difference which representatives we choose for the equivalence classes of terms $\syn{S}_i$'s
            and protypes $\syn{\alpha}_i$'s because of the replacement axioms,
            and the congruence problem does not arise
            because the equality judgments for protypes are limited to those coming from the equality judgments for terms
            by the replacement axiom.
        \end{itemize}
\end{definition}

\begin{proposition}
    \label{prop:crudevdc}
    The syntactic \ac{VDC} $\nS(\Sigma,\zero{E})$ for a specification $(\Sigma,\zero{E})$
    has a structure of a split \ac{CFVDC}.
\end{proposition}
\begin{proof}
    The tight structure is given as usual in algebraic theories.
    The composite of the following cells 
    \[
        \begin{tikzcd}[column sep=4em,virtual]
            \syn{\Gamma}_{1,0}
            \ar[d, "\syn{S}_0"']
            \sard[r, "\syn{\ol\alpha}_1"]
            \ar[dr, phantom, "\syn{\mu}_1"]
            & \syn{\Gamma}_{1,n_1}
            \ar[d, "\syn{S}_1"']
            \sard[r]
            & \cdots
            \sard[r, "\ol{\syn{\alpha}}_{n}"]
            \ar[dr, phantom, "\syn{\mu}_n"]
            & \syn{\Gamma}_{n,m_n}
            \ar[d, "\syn{S}_n"] \\
            \syn{\Delta}_{0}
            \ar[d, "\syn{T}_0"']
            \sar[r, "\syn{\beta}_1"'] 
            \ar[drrr, phantom, "\syn{\nu}"]
            & \syn{\Delta}_{1}
            \sar[r]
            & \cdots
            \sar[r, "\syn{\beta}_n"']  
            & \syn{\Delta}_{n}
            \ar[d, "\syn{T}_1"] \\
            \syn{\Theta}_{0}
            \sar[rrr, "\syn{\gamma}"']
            & & & \syn{\Theta}_{1}
        \end{tikzcd}
    \]
    is given as 
    \[
    \ol{\syn{\Gamma}}\mid\ol{\syn{\alpha}_1}\smcl\dots\smcl\ol{\syn{\alpha}_n}
    \vdash
    \syn{\nu}[\ol{\syn{S}_{\ul{i}}}/\ol{\syn{\Delta}_{\ul{i}}}]
    \psb{\syn{\mu}_1\smcl\dots\smcl\syn{\mu}_n}
    : \syn{\gamma}[\syn{T}_0/\syn{\Theta}_0\smcl\syn{T}_1/\syn{\Theta}_1][\syn{S}_0/\syn{\Delta}_0\smcl\syn{S}_n/\syn{\Delta}_n].
    \]
    The associativity and unit laws follow from \Cref{lem:subst}.

    The chosen restrictions are given by the term substitution into protypes.
    It is straightforward to check that
    the canonical cell 
    \[
        \begin{tikzcd}[virtual,column sep=12ex]
            \syn{\Gamma}_0
            \sar[r, "{\syn{\alpha}[\syn{S}_0/\syn{\Delta}_0\smcl\syn{S}_1/\syn{\Delta}_1]}"]
            \ar[d,"\syn{S}_0"']
            \ar[phantom,rd, "\restc" description]
            &
            \syn{\Gamma}_1
            \ar[d,"\syn{S}_1"]
            \\
            \syn{\Delta}_0
            \sar[r, "{\syn{\alpha}}"']
            &
            \syn{\Delta}_1
        \end{tikzcd}
        \quad
        \text{given by }
        \quad
        \syn{\Gamma}_0\smcl\syn{\Gamma}_1\mid \syn{a}:\syn{\alpha}[\syn{S}_0/\syn{\Delta}_0\smcl\syn{S}_1/\syn{\Delta}_1]
        \vdash \syn{a}:\syn{\alpha}[\syn{S}_0/\syn{\Delta}_0\smcl\syn{S}_1/\syn{\Delta}_1]
    \]
    exhibits $\syn{\alpha}[\syn{S}_0/\syn{\Delta}_0\smcl\syn{S}_1/\syn{\Delta}_1]$ 
    as a restriction of a loose arrow $\syn{\alpha}$ along $\syn{S}_0$ and $\syn{S}_1$ as tight arrows.
    The chosen terminals and binary products are given by the constructors $\top$ and $\land$,
    whose universal properties can be confirmed by the computation rules for them.
    By \Cref{lem:subst}, the choice gives a split \ac{CFVDC}.
\end{proof}

The functoriality is easy to check.

\begin{lemma}
    For any morphism of specifications $\Phi\colon(\Sigma,\zero{E})\to(\Sigma',\zero{E}')$,
    the translation $(-)^\Phi$ by $\Phi$
    defines a morphism $\nS(\Phi)\colon \nS(\Sigma,\zero{E})\to \nS(\Sigma',\zero{E}')$.
    This defines a (1-)functor $\nS\colon\Speci\to\FibVDblCart^\spl$.

\end{lemma}

\begin{theorem}
    \label{prop:crudeadj}
    The assignment that sends a \ac{CFVDC} $\dbl{D}$ to the associated specification $(\Sigma_{\dbl{D}},\zero{E}_{\dbl{D}})$
    extends to a functor $\Sp\colon\FibVDblCart^\spl\to\Speci$ which is a right adjoint to $\nS$.
    The counit components of the adjunction $\varepsilon_{\dbl{D}}\colon\nS(\zero{Sp}(\dbl{D}))\to\dbl{D}$ are
    an equivalence as a 1-cell in $\FibVDblCart$.
\end{theorem}
\begin{proof}
    We construct a virtual double functor $\varepsilon_\dbl{D}\colon\nS(\zero{Sp}(\dbl{D}))\to\dbl{D}$.
    We have the canonical $\Sigma_{\dbl{D}}$-structure in $\dbl{D}$.
    In the way we showed in \cref{sec:semantics},
    we can interpret all the items in $\zero{Sp}(\dbl{D})$ in $\dbl{D}$.
    Now, we show that this defines a virtual double functor from $\nS(\Sigma_{\dbl{D}},\zero{E}_{\dbl{D}})$ to $\dbl{D}$.
    The actions on the objects, tight arrows, and loose arrows are straightforward using \Cref{def:semantics}.
    A cell of $\nS(\Sigma_{\dbl{D}},\zero{E}_{\dbl{D}})$ of the form \cref{eq:cell} is
    interpreted as the composite of the cartesian cell on the left and the cell $\sem{\mu}$ on the right,
    which is inductively defined in \Cref{def:semantics}.
    \[
        \begin{tikzcd}[virtual, column sep=8ex]
            \sem{\syn{\Gamma}_0}
            \sar[r, "{\sem{\syn{\beta}[\syn{S}_0/\syn{\Delta}_0\smcl\syn{S}_1/\syn{\Delta}_1]}}"{yshift=1ex}]
            \ar[d, "{\sem{\syn{S}_0}}"']
            \ar[dr,phantom, "\restc" description]
            &
            \sem{\syn{\Gamma}_1}
            \ar[d, "{\sem{\syn{S}_1}}"]
            \\
            \sem{\syn{\Delta}_0}
            \sar[r, "{\sem{\syn{\beta}}}"']
            &
            \sem{\syn{\Delta}_1}
        \end{tikzcd}
        \hspace{1em}
        ,
        \hspace{1em}
        \begin{tikzcd}[column sep=8ex,virtual]
            {\sem{\syn{\Gamma}_0}}
            \sar[r, "{\sem{\syn{\alpha}_1}}"]
            \ar[d, equal]
            \ar[phantom,rrrd, "{\sem{\syn{\mu}}}" description]
            &
            \cdots
            &
            \cdots
            \sar[r, "{\sem{\syn{\alpha}_n}}"]
            &
            {\sem{\syn{\Gamma}_n}}
            \ar[d,equal]
            \\
            {\sem{\syn{\Gamma}_0}}
            \sar[rrr, "{\sem{\syn{\beta}[\syn{S}_0/\syn{\Delta}_0\smcl\syn{S}_1/\syn{\Delta}_1]}}"'] 
            &
            &&
            {\sem{\syn{\Gamma}_n}}
        \end{tikzcd}
    \]
    These assignments are independent of the choice of terms and proterms 
    since in $\nS(\Sigma_{\dbl{D}},\zero{E}_{\dbl{D}})$, we take equivalence classes with respect to the equality judgments belonging to $\zero{E}_{\dbl{D}}$.
    Proving that this defines a morphism in $\FibVDblCart^\spl$ is a routine verification.
    For instance, it sends a chosen restriction $\syn{\alpha}[\syn{S}_0/\syn{\Delta}_0\smcl\syn{S}_1/\syn{\Delta}_1]$
    of $\syn{\alpha}$ along $\syn{S}_0$ and $\syn{S}_1$ to 
    $\sem{\syn{\alpha}[\syn{S}_0/\syn{\Delta}_0\smcl\syn{S}_1/\syn{\Delta}_1]}$,
    which is the same as $\sem{\syn{\alpha}}[\sem{\syn{S}_0}\smcl\sem{\syn{S}_1}]$ by \Cref{lemma:interpretsubs}.

    We show that $\varepsilon_{\dbl{D}}$ is an equivalence as a virtual double functor.
    The surjectiveness part directly follows from the construction.
    The proofs of the fully-faithfulness on tight arrows and cells are parallel:
    if two terms or proterms in $\Sp(\dbl{D})$ are interpreted as the same term or proterm in $\dbl{D}$,
    then this equality is reflected in the equality judgments in $\zero{E}_{\dbl{D}}$,
    and hence the terms or proterms are already derivably equal in $\Sp(\dbl{D})$. 

    Now, we show that $\varepsilon_{\dbl{D}}$ is a terminal object in the comma category $\nS\,\downarrow\,\dbl{D}$.
    Suppose we are given a morphism $F\colon\nS(\Sigma,\zero{E})\to\dbl{D}$.
    If $\wh{F}\colon(\Sigma,\zero{E})\to\Sp(\dbl{D})$ satisfies $\varepsilon_{\dbl{D}}\circ\nS(\wh{F})=F$,
    then it satisfies the following:
    \begin{itemize}
        \item $\syn{x}:\wh{F}(\syn{\sigma})$ is interpreted as $F(\syn{x}:\syn{\sigma})$ in $\dbl{D}$ for each category symbol $\syn{\sigma}$, 
        \item $(\wh{F}(\syn{f}))(\syn{x})$ is interpreted as $F(\syn{f}(\syn{x}))$ in $\dbl{D}$ for each function symbol $\syn{f}$, 
        \item $(\wh{F}(\syn{\rho}))(\syn{x}\smcl\syn{y})$ is interpreted as $F(\syn{\rho}(\syn{x}\smcl\syn{y}))$ in $\dbl{D}$ for each profunctor symbol $\syn{\rho}$, and 
        \item $(\wh{F}(\syn{\kappa}))(\ol{\syn{x}_i})\{\ol{\syn{a}_i}\}$ is interpreted as $F(\syn{\kappa}(\ol{\syn{x}_i})\{\ol{\syn{a}_i}\})$ in $\dbl{D}$ for each proterm symbol $\syn{\kappa}$. 
    \end{itemize}
    However, $\varepsilon_{\dbl{D}}$ is injective on primitive contexts and procontexts,
    and also is injective on the terms and proterms by the fully-faithfulness.
    Hence, $\wh{F}$ is uniquely determined for $F$ by the above conditions:
    \[
    \wh{F}(\syn{\sigma})=\is{F(\syn{x}:\syn{\sigma})},\quad
    \wh{F}(\syn{f})=\is{F(\syn{f}(\syn{x}))},\quad
    \wh{F}(\syn{\rho})=\is{F(\syn{\rho}(\syn{x},\syn{y}))},\quad
    \wh{F}(\syn{\kappa})=\is{F(\syn{\kappa}(\ol{\syn{x}_i})\{\ol{\syn{a}_i}\})}.
    \]
    Conversely, the assignment $\wh{F}$ defined by the above gives a morphism 
    $\wh{F}\colon(\Sigma,\zero{E})\to\Sp(\dbl{D})$.
    The well-definedness of $\wh{F}$ depends on the fact that a equality judgment in $\zero{E}$ 
    induces an equality in $\nS(\Sigma,\zero{E})$,
    which is sent to an equality in $\dbl{D}$ by $F$.
    It also satisfies the equation $\varepsilon_{\dbl{D}}\circ\nS(\wh{F})=F$,
    which is confirmed by induction on the structure of the judgments in
    $(\Sigma,\zero{E})$.
    Therefore, $\varepsilon_{\dbl{D}}$ has the desired universal property.
\end{proof}

\begin{remark}
    Owing to the splitness lemma \Cref{lemma:split},
    this adjunction achieves the desired syntax-semantics duality without loss of generality.
    It would be more precise to say that this 1-adjunction combines with 
    the biequivalence between the 2-category of split \acp{CFVDC} and
    the 2-category of (cloven) \acp{CFVDC} to form a biadjunction.
\end{remark}

\subsection{Specifications with protype isomorphisms}
We can extend the biadjunction to the type theory with protype isomorphisms.
First, we introduce a notion of specification with protype isomorphisms.
We use the term ``\emph{multi-class}'' to mean a class $\zero{X}$ with multiplicities $(\zero{M}_x)_{x\in\zero{X}}$,
where $\zero{M}_x$ is a class.
One can think of a multi-class as a (class-large) family of classes.

\begin{definition}
    By a \emph{multi-class} $(\zero{M})_{x}$, we mean a class $\zero{X}$ with multiplicities $(\zero{M}_x)_{x\in\zero{X}}$,
    where $\zero{M}_x$ is a class.
    A \emph{multi-class of isomorphism symbols} for a signature $\Sigma$ is a multi-class
    $\zero{PI}_{\syn{\rho},\syn{\omega}}$ indexed by pairs of profunctor symbols $(\syn{\rho},\syn{\omega})$ of the same two-sided arity
    in $\Sigma$.
    We call the elements of $\zero{PI}_{\syn{\rho},\syn{\omega}}$ \emph{isomorphism symbols}.
\end{definition}

\begin{definition}
    A \emph{specification with protype isomorphisms} $(\Sigma,\zero{PI},\zero{E})$ consists of
    \begin{itemize}
        \item a signature $\Sigma$,
        \item $\zero{PI}$, a multi-class of isomorphism symbols for $\Sigma$, and
        \item a pair $(\zero{E}^{\tm},\zero{E}^{\ptm})$ as in \Cref{def:specification},
        but the derivation of proterms can refer to the following rule.
        {\small
        \[
            \inferrule*
                {\syn{m}\in \zero{PI}_{\syn{\rho},\syn{\omega}}}
                {\syn{x}:\syn{\sigma}\smcl\syn{y}:\syn{\tau}\vdash \syn{\Lambda}_\syn{m} : \syn{\rho}(\syn{x}\smcl\syn{y})\ccong\syn{\omega}(\syn{x}\smcl\syn{y})}
        \]
        }
    \end{itemize}

    A \emph{morphism of specifications with protype isomorphisms} 
    $\Phi\colon(\Sigma,\zero{PI},\zero{E})\to(\Sigma',\zero{PI}',\zero{E}')$ 
    consists of a morphism of signatures $\Phi\colon\Sigma\to\Sigma'$ and
    a multi-class function $\breve{\Phi}\colon\zero{PI}_{\syn{\rho},\syn{\omega}}\to\zero{PI}'_{\Phi(\syn{\rho}),\Phi(\syn{\omega})}$
    compatible with the index function of $\zero{PI}$ defined by $\Phi$
    such that every judgment in $\zero{E}$ is translated to a judgment that is derivable from $\zero{E}'$ by $(\Phi,\breve{\Phi})$. 

    We write $\Speci^{\ccong}$ for the 2-category of specifications with protype isomorphisms and morphisms between them. 
\end{definition}

We will construct a functor $\Ufd\colon\Speci^{\ccong}\to\Speci$
which has a partial right adjoint. 
Since the right adjoint is defined on the image of $\Sp$,
we will obtain an adjunction between the category of specifications with protype isomorphisms
and the category of split \acp{CFVDC} in the end.

\begin{definition}
    We define a specification (without protype isomorphisms) $\Ufd(\Sigma,\zero{PI},\zero{E})$ for a specification 
    with protype isomorphisms $(\Sigma,\zero{PI},\zero{E})$ as follows.
    \begin{itemize}
        \item the signature consists of data in $\Sigma$ plus additional transformation symbols $\syn{\varphi}_{m}\colon \syn{\rho}\Rightarrow\syn{\omega}$ 
        and $\syn{\psi}_{m}\colon \syn{\omega}\Rightarrow\syn{\rho}$ for each element $\syn{m}\in\zero{PI}_{\syn{\rho},\syn{\omega}}$,
        \item the equality judgments consist of the original equality judgments in $\zero{E}$ 
        with all occurrences of protype isomorphisms inductively replaced by the corresponding proterms 
        as shown in \Cref{fig:crudeencoding},
        plus the following additional equality judgments:
        \begin{equation}
        \label{eq:crudeencoding}
        \syn{x}:\syn{\sigma}\smcl\syn{y}:\syn{\tau}\mid \syn{a}:\syn{\rho}\vdash \syn{\psi}_{m}\{\syn{\varphi}_{m}\{\syn{a}\}\}\equiv\syn{a} : \syn{\rho}\quad \text{and}\quad
        \syn{x}:\syn{\sigma}\smcl\syn{y}:\syn{\tau}\mid \syn{b}:\syn{\omega}\vdash \syn{\varphi}_{m}\{\syn{\psi}_{m}\{\syn{b}\}\}\equiv\syn{b} : \syn{\omega}
        \end{equation}
        for each $m\in\zero{PI}_{\syn{\rho},\syn{\omega}}$.
    \end{itemize}
\end{definition} 

\begin{figure}[h]
    {\small
    \begin{align*}
        \idt_{\syn{\alpha}}\{\syn{a}\} &\rightsquigarrow \syn{a} & \lcp\syn{\mu},\syn{\nu}\rcp\{\syn{a}\} &\rightsquigarrow \syn{\mu}\{\syn{a}\}\\
        \idt_{\syn{\alpha}}\sinv\{\syn{a}\} &\rightsquigarrow \syn{a} & \lcp\syn{\mu},\syn{\nu}\rcp\sinv\{\syn{a}\} &\rightsquigarrow \syn{\nu}\{\syn{a}\}\\
        (\syn{\Omega}\circ\syn{\Upsilon})\{\syn{a}\} &\rightsquigarrow \syn{\Omega}\{\syn{\Upsilon}\{\syn{a}\}\} & \syn{\Lambda}_m\{\syn{a}\} &\rightsquigarrow \syn{\varphi}_m\{\syn{a}\}\\
        (\syn{\Omega}\circ\syn{\Upsilon})\sinv\{\syn{a}\} &\rightsquigarrow \syn{\Upsilon}\sinv\{\syn{\Omega}\sinv\{\syn{a}\}\} & \syn{\Lambda}_m\sinv\{\syn{a}\} &\rightsquigarrow \syn{\psi}_m\{\syn{a}\}
    \end{align*}   
    }
    \caption{Translation of protype isomorphisms}
    \label{fig:crudeencoding}
\end{figure}

\begin{lemma}
    \label{lemma:crudeencoding}
    The assignment $(\Sigma,\zero{PI},\zero{E})\mapsto \Ufd(\Sigma,\zero{PI},\zero{E})$ induces a functor 
    $\Ufd\colon\Speci^{\ccong}\to\Speci$.
\end{lemma}
\begin{proof}[Proof sketch]
    For a morphism of specifications $\Phi\colon(\Sigma,\zero{E})\to(\Sigma',\zero{E}')$,
    the assignment $\Ufd(\Phi)$ sends the transformation symbols $\syn{\varphi}_{m}$ and $\syn{\psi}_{m}$ to
    $\syn{\varphi}_{\Phi(m)}$ and $\syn{\psi}_{\Phi(m)}$.
    The equality judgments \Cref{eq:crudeencoding} are translated into the equality judgments of the same form
    and hence derivable from $\Ufd(\Sigma',\zero{E}')$.
\end{proof}

The functor does not have a right adjoint globally but a partial one.
\begin{definition}
    A specification $(\Sigma,\zero{E})$ is \emph{unary-cell-saturated} if, for 
    any proterm judgment $\syn{x}:\syn{\sigma}\smcl\syn{y}:\syn{\tau}\mid \syn{a}:\syn{\rho}\vdash \syn{\vartheta}:\syn{\omega}$ derivable from $\zero{E}$
    where $\syn{\sigma},\syn{\tau},\syn{\rho},\syn{\omega}$ belongs to the signature $\Sigma$,
    there uniquely exists a transformation symbol $\syn{\kappa}_{\syn{\vartheta}}\colon \syn{\rho}\Rightarrow\syn{\omega}$ in $\Sigma$
    such that the equality judgment
    \[
        \syn{x}:\syn{\sigma}\smcl\syn{y}:\syn{\tau}\mid \syn{a}:\syn{\rho}\vdash \syn{\kappa}_{\syn{\vartheta}}(\syn{x}\smcl\syn{y})\{\syn{a}\}\equiv\syn{\vartheta}:\syn{\omega}
    \]
    is derivable from $\zero{E}$.
    Let $\Speci_{\essat}$ be the full subcategory of $\Speci$
    whose objects are unary-cell-saturated crude specifications.
\end{definition}

It is easy to see that
the associated specification $(\Sigma_{\dbl{D}},\zero{E}_{\dbl{D}})$ of a \ac{CFVDC} $\dbl{D}$
is unary-cell-saturated.
A specification being saturated means that 
the symbols in the signature constitute a virtual double category that 
is equivalent to the syntactic \ac{VDC} of the specification.

\begin{proposition}
    \label{prop:partialbiadj}
    The functor $\Ufd\colon\Speci^{\ccong}\to\Speci$ has a relative right coadjoint $\Fd$
    over the inclusion $J\colon\Speci_{\essat}\hto\Speci$.
    \[
        \begin{tikzcd}[ampersand replacement=\&]
            \Speci^{\ccong}
            \ar[r, shift left, "\Ufd"]
            \ar[dr,phantom, "{\overset{\upsilon}{\Rightarrow}}",xshift=2ex, yshift=2ex]
            \&
            \Speci
            \\
            \&
            \Speci_{\essat} 
            \ar[u, hook, "J"']
            \ar[ul,"\Fd"]
        \end{tikzcd}.
    \]
    The components of the counit $\upsilon_{(P,\zero{D})}\colon\Ufd(\Fd(P,\zero{D}))\to(P,\zero{D})$
    are sent to the equivalence by $\nS$.
\end{proposition}

Here, the relative right coadjoint
means that there is a natural isomorphism
\[
    \Speci(\Ufd(-),J(*))\cong\Speci^{\ccong}(-,\Fd(*))
\]
induced by the $\upsilon$.

\begin{proof}
    For a unary-cell-saturated crude specification $(P,\zero{D})$,
    a specification $\Fd(P,\zero{D})$ consists of the same signature $P$,
    the multi-class $\zero{D}^{\cong}$ defined from $\zero{D}$ by setting
    $\zero{D}^{\cong}_{\syn{\rho},\syn{\omega}}$ to be the class of the pairs $(\syn{\vartheta},\syn{\varsigma})$ of transformation symbols in $\zero{D}$
    \[ 
        \syn{\vartheta}\colon\syn{\rho}\Rightarrow\syn{\omega}\quad\text{and}\quad\syn{\varsigma}\colon\syn{\omega}\Rightarrow\syn{\rho}
    \]
    for which $\zero{D}$ 
    derives the equality judgments
    that express the two cells are inverses of each other,
    and the classes of term and proterm equality judgments
    in $\zero{D}$
    plus the equality judgments
    \begin{align*}
        \label{eq:ruleforDcd} 
        \syn{x}:\syn{\sigma}\smcl\syn{y}:\syn{\tau}\mid \syn{a}:\syn{\rho}(\syn{x}\smcl\syn{y})&\vdash 
        \syn{\Lambda}_{(\syn{\vartheta},\syn{\varsigma})}\{\syn{a}\}\equiv\syn{\vartheta}(\syn{x}\smcl\syn{y})\{\syn{a}\}:\syn{\omega}(\syn{x}\smcl\syn{y})
        \\
        \syn{x}:\syn{\sigma}\smcl\syn{y}:\syn{\tau}\mid \syn{b}:\syn{\omega}(\syn{x}\smcl\syn{y})&\vdash
        \syn{\Lambda}_{(\syn{\vartheta},\syn{\varsigma})}\inv\{\syn{b}\}\equiv\syn{\varsigma}(\syn{x}\smcl\syn{y})\{\syn{b}\}:\syn{\rho}(\syn{x}\smcl\syn{y})
    \end{align*}
    for each isomorphism symbol $(\syn{\vartheta},\syn{\varsigma})$ in $\zero{D}^{\cong}_{\syn{\rho},\syn{\omega}}$.
    Then we will have a morphism of specifications $\upsilon_{(P,\zero{D})}$
    that sends the new transformation symbols $\syn{\varphi}_{(\syn{\vartheta},\syn{\varsigma})}$ 
    and $\syn{\psi}_{(\syn{\vartheta},\syn{\varsigma})}$ to the transformation symbols $\syn{\vartheta}$ and $\syn{\varsigma}$.
    It follows that $\upsilon_{(P,\zero{D})}$ defines a morphism of specifications
    since the equality judgments in $\Ufd(\Fd(P,\zero{D}))$
    are either in $\zero{D}$ or 
    those of the form \Cref{eq:crudeencoding} for the pairs in $\zero{D}^{\cong}$,
    which are translated to equality judgments
    derivable from $\zero{D}$.

    We prove that this $\upsilon_{(P,\zero{D})}$ satisfies the universal property
    for the relative right coadjoint of $\Ufd$.
    That is, for a morphism of specifications $\Phi\colon\Ufd(\Sigma,\zero{PI},\zero{E})\to(P,\zero{D})$,
    there uniquely exists a morphism of specifications with protype isomorphisms
    $\wh{\Phi}\colon(\Sigma,\zero{PI},\zero{E})\to\Fd(P,\zero{D})$ 
    such that the following diagram commutes
    \[
        \begin{tikzcd}
            \Ufd(\Sigma,\zero{PI},\zero{E})\ar[rd,"\Phi"]
            \ar[d,"\Ufd(\wh{\Phi})"']
            &
            \!
            \ar[dl,phantom, "{\rotatebox{45}{$=$}}", xshift=-2.5ex, yshift=-2ex, description] 
            \\
            \Ufd(\Fd(P,\zero{D}))\ar[r,"\upsilon_{(P,\zero{D})}"']
            &
            (P,\zero{D})
        \end{tikzcd}    
        \quad
        \text{in}\ \Speci.
    \]
    To make this diagram commute,
    the signature part of $\wh{\Phi}$ must be the same as $\Phi$.
    Suppose we have a morphism $\wh{\Phi}$ and we determine how it should act on the
    isomorphism symbols in $\zero{PI}$.
    Let $(\syn{\chi}_\syn{m},\syn{\lambda}_\syn{m})$ be the image of $\syn{m}$ under $\wh{\Phi}$.
    Then, the symbol $\syn{\chi}_\syn{m}$ equals to $\upsilon_{(P,\zero{D})}(\syn{\varphi}_{(\syn{\chi}_\syn{m},\syn{\lambda}_\syn{m})})
    =\upsilon_{(P,\zero{D})}(\syn{\varphi}_{\wh{\Phi}(\syn{m})})$,
    which is the image of $\syn{m}$ under $\Phi$.
    Similarly, we must have $\syn{\lambda}_\syn{m}=\Phi(\syn{\psi}_{\syn{m}})$.
    Therefore, the morphism $\wh{\Phi}$ must send $\syn{m}$ to the pair $(\Phi(\syn{\varphi}_{\syn{m}}),\Phi(\syn{\psi}_{\syn{m}}))$. 
    This assignment $\wh{\Phi}$ is a morphism of specifications with protype isomorphisms
    since the equality judgments in $\zero{E}$ 
    with the isomorphism symbols suitably replaced
    are translated by $\Phi$ to 
    the equality judgments provable from $\zero{D}$.
    Note that the proterm $\syn{\Lambda}_{\syn{m}}\{\syn{a}\}$ is sent to
    $\syn{\Lambda}_{\wh{\Phi}(\syn{m})}\{\syn{a}\}$,
    which behaves the same as $\Phi(\syn{\varphi}_{\syn{m}})(\syn{x}\smcl\syn{y})\{\syn{a}\}$
    up to derivable equality in $\zero{D}$.
    
    To see that $\nS(\upsilon_{(P,\zero{D})})$ is an equivalence, we confer \Cref{lemma:fibvdblequiv}.
    The equivalence on the tight part is apparent since $\upsilon_{(P,\zero{D})}$ does not change anything on types and terms.
    Next, for each loose arrow in $\nS(\Ufd(\Fd(P,\zero{D})))$,
    we can find a corresponding loose arrow in $\nS(\Ufd(\Fd(P,\zero{D})))$ by taking the protype with precisely the same presentation.
    Finally, when fixing a frame, the function on globular cells defined by $\upsilon_{(P,\zero{D})}$ sends proterm judgments with the additional transformation symbols $\syn{\varphi}_{(\syn{\vartheta},\syn{\varsigma})}$ and $\syn{\psi}_{(\syn{\vartheta},\syn{\varsigma})}$
    to the proterm judgments without them by replacing those transformation symbols with $\syn{\vartheta}$ and $\syn{\varsigma}$.
    The surjectiveness is checked similarly to the above argument.
    We can also see the injectiveness up to derivable equality by induction on the construction of the proterms.
    For instance,
    the equalities $\syn{\varphi}_{(\syn{\vartheta},\syn{\varsigma})}(\syn{x}\smcl\syn{y})\{\syn{a}\}\equiv\syn{\vartheta}(\syn{x}\smcl\syn{y})\{\syn{a}\}$ 
    and $\syn{\psi}_{(\syn{\vartheta},\syn{\varsigma})}(\syn{x}\smcl\syn{y})\{\syn{a}\}\equiv\syn{\varsigma}(\syn{x}\smcl\syn{y})\{\syn{a}\}$ are already derivable from $\Ufd(\Fd(P,\zero{D}))$. 
\end{proof}

    \begin{corollary}
    The composite $\nS\circ\Ufd\colon\Speci^{\ccong}\to\FibVDblCart^\spl$ 
    has a right adjoint $\Fd\circ\Sp$:
        \[
            \begin{tikzcd}[ampersand replacement=\&]
                \Speci^{\ccong} 
                \ar[r, shift left=2, "\nS\circ\Ufd"]
                \ar[r, phantom, "\rotatebox{90}{$\vdash$}"]
                \&
                \FibVDblCart^\spl
                \ar[l, shift left=2]
            \end{tikzcd},
            \quad 
            \text{given by}\quad
            \quad
            \begin{tikzcd}[ampersand replacement=\&, row sep=2ex]
                \Speci^{\ccong}
                \ar[r, shift left, "\Ufd"]
                \ar[rr, shift right=3ex, phantom, "\rotatebox{270}{$\dashv$}"]
                \&
                \Speci
                \ar[r, shift left, "{\dbl{S}}"]
                \&
                \FibVDblCart^\spl \ar[ld, shift left, "{\Sp}"]
                \\
                \&
                \Speci_{\essat}
                \ar[lu, shift left, "{\Fd}"]
            \end{tikzcd}.
        \]
    Moreover, the counit component of the adjunction is pointwise an equivalence
    as a virtual double functor.
    \end{corollary}
    \begin{proof}
        Through \Cref{prop:crudeadj,prop:partialbiadj}, the expected adjunction follows from the general theory of relative coadjunctions.
        Explicitly,
        for a specification $\zero{S}$ and a \ac{CFVDC} $\dbl{D}$,
        \begin{align*}
            \FibVDblCart^\spl\left(\nS(\Ufd(\zero{S})),\dbl{D}\right)
            &\cong \Speci\left(\Ufd(\zero{S}),\Sp(\dbl{D})\right)& (\text{by\ \Cref{prop:crudeadj}})\\
            &\cong \Speci^{\ccong}\left(\zero{S},\Fd(\Sp(\dbl{D}))\right)& (\text{by\ \Cref{prop:partialbiadj}})\\
        \end{align*}
    The counit component of the adjunction is an equivalence by the construction of the adjunctions.
    \end{proof}
    
\begin{remark}
    The specification $\Fd(\Sp(\dbl{D}))$ is not the same as the associated specification $(\Sigma_{\dbl{D}},\zero{E}_{\dbl{D}})$
    equipped with the isomorphism symbols, but the two give the equivalent virtual double categories.
\end{remark}

\begin{remark}
    For extensions of \ac{FVDblTT} with additional constructors as in \Cref{sec:additional},
    we can also obtain a syntax-semantics biadjunction analogously
    once one determines the treatment of substitutions as explained \Cref{rem:substitution-into-additional-constructor}.
    The procedure goes as follows:
    (i) Prove the splitness lemma for \acp{CFVDC} with the additional structure of interest,
    where the splitness is defined in reflection of the treatment of substitutions;
    (ii) Construct the syntactic \acp{VDC} for the extended type theory and verify
    that they have the structures in question;
    (iii) Prove the adjunction between the category of split \acp{CFVDC} with the additional structures
    and the category of specifications with the additional constructors in the same way as in \Cref{prop:crudeadj}.
    The biadjunction is again obtained by combining this adjunction with the biequivalence between the 2-categories of split and cloven \acp{CFVDC} with the structures. 
\end{remark}

%% file: newsubfiles/discussion/conclusion.tex
There are several directions for future work.
First, we would like to extend the type theory \ac{FVDblTT}
to include more advanced structures studied in formal category theory using virtual double categories.  
In particular, we are interested in the extension of the type theory \ac{FVDblTT} to
\textit{augmented} virtual double categories \cite{koudenburgAugmentedVirtualDouble2020,koudenburgFormalCategoryTheory2024}.
The latter paper conceptualizes the notion of a Kan extension and a Yoneda embedding inside this framework
and develops formal category theory more flexibly than the original virtual double categories.
Second, the dependent version of the type theory \ac{FVDblTT} should be developed 
from the perspective of directed type theory.
There are several studies on directed type theory \cite{licata2DimensionalDirectedType2011,northDirectedHomotopyType2019,ahrensBicategoricalTypeTheory2023},
and those are all based on dependent types.
One of the primary objectives of those studies is to obtain a substantial type theory for higher categories 
as Martin-L\"of type theory is for higher groupoids.
The dependent version of the type theory \ac{FVDblTT} might offer another candidate for this purpose
using the unit protypes and the comprehension types.
Finally, we are interested in the relationship between the type theory \ac{FVDblTT} and 
other type theories or calculi for relations.
In particular, we are interested in the connection to diagrammatic calculi for relations 
such as the one in \cite{bonchiFunctorialSemanticsRelational2017,bonchiDiagrammaticAlgebraFirst2024}
or, more directly, the string diagrams for double categories \cite{myers2018stringdiagramsdoublecategories}.
They may be understood as a string-diagrammatic presentation of the type theory \ac{FVDblTT}.
We hope to explore these connections in future work.

%% file: newsubfiles/appendix/cartesianstr.tex
We provide rationale for the rules in \cref{sec:unit-protype-meets-product-type,sec:compo-protype-meets-product-type,sec:filler-protype-meets-product-type},
by unpacking the cartesianness of virtual double categories with structures.
\begin{lemma}
    \label{lem:cartobjinlffsub}
    Let $\bi{B},\bi{B}'$ be 2-categories with finite products $(\monunit,\otimes)$,
    and $\abs{-}\colon\bi{B}'\to\bi{B}$ be a 2-functor
    preserving finite products and locally full-inclusion, \textit{i.e.}, injective on 1-cells and bijective on 2-cells.
    For an object $x$ of $\bi{B}'$ to be cartesian,
    it is necessary and sufficient that $\abs{x}$ is cartesian in $\bi{B}$ and that
    the 1-cells $1\colon\monunit\to \abs{x}$ and $\times\colon \abs{x}\otimes \abs{x}\to \abs{x}$ 
    right adjoint to the canonical 1-cells are essentially in the image of $\abs{-}$.

    Moreover, for a 1-cell $f\colon x\to y$ of $\bi{B}'$ where $x$ and $y$ are cartesian in $\bi{B}'$,
    $f$ is cartesian in $\bi{B}'$ if and only if $\abs{f}$ is cartesian in $\bi{B}$.
\end{lemma}
\begin{proof}
    The necessity of the first condition follows from the fact that any 2-functor preserves adjunctions,
    that right adjoints are unique up to isomorphism, and that $\abs{-}$ preserves finite products.
    Since $\abs{-}$ is locally fully faithful, it also reflects units, counits, and the triangle identities with respect to the adjunctions,
    and hence the sufficiency of the first condition follows.

    The necessity of the second condition is again immediate from the fact that $\abs{-}$ preserves finite products.
    The sufficiency follws from the fact that $\abs{-}$ is locally fully faithful, in particular, reflects isomorphisms.
\end{proof}

\begin{proposition}
    \label{prop:cartesianunital}
    Let $\FibUVDbl$ be the locally-full sub-2-category of $\FibVDbl$ spanned by the \acp{FVDC} with units
    and functors preserving units.
    Then, a \ac{FVDC} $\dbl{D}$ with units is cartesian in $\FibUVDbl$ 
    if and only if
    \begin{enumerate}
        \item $\dbl{D}$ is a cartesian \ac{FVDC},
        \item $U_1\cong\top_{1,1}$ in $\dbl{D}(1,1)$ canonically, and
        \item for any $I,J\in\dbl{D}$, $U_{I,J}\cong U_I\times U_J$ canonically in $\dbl{D}(I\times J,I\times J)$.
    \end{enumerate}
\end{proposition}
\begin{proof}
By \cref{lem:cartobjinlffsub}, $\dbl{D}$ is cartesian as a unital \ac{FVDC} if and only if 
it is cartesian as a \ac{FVDC} and the 1-cells $1\colon\dbl{1}\to \dbl{D}$ and $\times\colon \dbl{D}\times \dbl{D}\to \dbl{D}$ are in $\FibUVDbl$.
The first condition is equivalent to \textit{(ii)} since it sends the only loose arrow in $\dbl{1}$, which is the unit loose arrow, to $\top_{1,1}$.
The second condition is equivalent to \textit{(iii)} since the unit loose arrow of $(I,J)$ in $\dbl{D}(I\times J,I\times J)$ is $(\delta_I,\delta_J)$,
which is sent to $\delta_I\times\delta_J$ in $\dbl{D}(I\times J,I\times J)$.
\end{proof}

The key idea is that in the virtual double categories $\dbl{D}\times\dbl{D}$ and $\dbl{1}$,
the unit loose arrows are given pointwise by the unit loose arrows of $\dbl{D}$.
We can discuss the cartesianness of some classes of \acp{FVDC} in parallel with the above proposition.

\begin{proposition}
    \label{prop:cartesiancompose}
    Let $\FibCVDbl$ be the locally-full sub-2-category of $\FibVDbl$ spanned by the \acp{FVDC} with composites 
    of sequences of loose arrows of positive length and functors preserving those composites.
    A \ac{VDC} $\dbl{D}$ in $\FibCVDbl$ is cartesian in this 2-category if and only if
    \begin{enumerate}
        \item $\dbl{D}$ is a cartesian \ac{FVDC},
        \item $\top_{1,1}\odot\dots\odot\top_{1,1}\cong\top_{1,1}$ canonically in $\dbl{D}(1,1)$, and
        \item for any paths of positive length 
        \[
        \begin{tikzcd}
            I_0\sar["\alpha_1",r] & I_1\sar[r] & \dots\sar[r,"\alpha_n"] & I_n
        \end{tikzcd}  
        \quad
        \text{and}
        \quad
        \begin{tikzcd}
            J_0\sar["\beta_1",r] & J_1\sar[r] & \dots\sar[r,"\beta_n"] & J_n
        \end{tikzcd}
        \]
        in $\dbl{D}$,
        we have
        \[
        (\alpha_1\odot\dots\odot\alpha_n)\times(\beta_1\odot\dots\odot\beta_n)\cong(\alpha_1\times\beta_1)\odot\dots\odot(\alpha_n\times\beta_n)
        \]
        canonically in $ \dbl{D}(I_0\times J_0,I_n\times J_n)$.
    \end{enumerate}
\end{proposition}

\begin{proposition}
    \label{prop:cartesianextension}
    Let $\FibEVDbl$ be the locally-full sub-2-category of $\FibVDbl$ spanned by the \acp{FVDC} with right extensions
    and functors preserving right extensions.
    A \ac{VDC} $\dbl{D}$ in $\FibEVDbl$ is cartesian in this 2-category if and only if
    \begin{enumerate}
        \item $\dbl{D}$ is a cartesian \ac{FVDC},
        \item $\top_{1,1}\triangleright\top_{1,1}\cong\top_{1,1}$ canonically in $\dbl{D}(1,1)$, and
        \item for any quadruples of loose arrows
        \[
        \begin{tikzcd}
            I_0\sar["\alpha_1",r]
            \sar[rr,"\alpha_2"',bend right=20]
            & I_1 &
            I_2
        \end{tikzcd}
        \quad
        \text{and}
        \quad
        \begin{tikzcd}
            J_0\sar["\beta_1",r]
            \sar[rr,"\beta_2"',bend right=20]
            & J_1 &
            J_2
        \end{tikzcd}
        \]
        in $\dbl{D}$,
        we have
        \[
        (\alpha_1\triangleright\alpha_2)\times(\beta_1\triangleright\beta_2)
        \cong
        (\alpha_1\times\beta_1)\triangleright(\alpha_2\times\beta_2)
        \]
        canonically in $\dbl{D}(I_1\times J_1,I_2\times J_2)$.
    \end{enumerate}
\end{proposition}

%% file: newsubfiles/appendix/cartesiansyn.tex
In \Cref{sec:additional}, we explain some additional constructors of \ac{FVDblTT} 
that are meaningful both in the contexts of formal category theory and predicate logic.
In this section, we provide all the derivation rules of the constructs.

\myparagraph{Unit protype.}\ 
\label{sec:unit-protype}
\begin{mathparpagebreakable}
    \goodbreak
    \inferrule*[right=$\ides$-Form]
    {\syn{I}\ \textsf{type}\\
    \syn{\Gamma}\vdash \syn{s}:\syn{I}\\
    \syn{\Delta}\vdash \syn{t}:\syn{I}}
    {\syn{\Gamma}\smcl\syn{\Delta}\vdash \syn{s}\ide{\syn{I}}\syn{t}\ \textsf{protype}}
    \and
    \inferrule*[right=$\ides$-Intro]
    {\syn{I}\ \textsf{type}}
    {\syn{x}:\syn{I}\mid \quad \vdash \refl_{\syn{I}}(\syn{x}): \syn{x}\ide{\syn{I}}\syn{x}}
    \and
    \inferrule*[right=$\ides$-Elim]
    {\syn{w}_0:\syn{J}_0\smcl\syn{z}_m:\syn{K}_m\vdash \syn{\gamma}(\syn{w}_0\smcl\syn{z}_m)\  \textsf{protype}\\
    \ol{\syn{w}}:\ol{\syn{J}}\smcl\syn{x}:\syn{I}\smcl\ol{\syn{z}}:\ol{\syn{K}}\mid \ol{\syn{A}}(\ol{\syn{w}}\smcl\syn{x})\smcl\ol{\syn{B}}(\syn{x}\smcl\ol{\syn{z}})\vdash \syn{\mu}:\syn{\gamma}(\syn{w}_0\smcl\syn{z}_m)}
    {\ol{\syn{w}}:\ol{\syn{J}}\smcl\syn{x}:\syn{I}\smcl\syn{y}:\syn{I}\smcl\ol{\syn{z}}:\ol{\syn{K}}\mid \ol{\syn{A}}(\ol{\syn{w}}\smcl\syn{x})\smcl\syn{p}:\syn{x}\ide{\syn{I}}\syn{y}\smcl\ol{\syn{B}}(\syn{y}\smcl\ol{\syn{z}})\vdash
    \ideind{\syn{I}}\{\syn{\mu}\}:\syn{\gamma}(\syn{w}_0\smcl\syn{z}_m)}
    \and 
    \inferrule*[right=$\ides$-Comp$\beta$]
    {\ol{\syn{w}}:\ol{\syn{J}}\smcl\syn{x}:\syn{I}\smcl\ol{\syn{z}}:\ol{\syn{K}}\mid \ol{\syn{A}}(\ol{\syn{w}}\smcl\syn{x})\smcl\ol{\syn{B}}(\syn{x}\smcl\ol{\syn{z}})\vdash \syn{\mu}:\syn{\gamma}(w_0\smcl\syn{z}_m)}
    {
    \ol{\syn{w}}:\ol{\syn{J}}\smcl\syn{x}:\syn{I}\smcl\ol{\syn{z}}:\ol{\syn{K}}\mid
    \ol{\syn{A}}(\ol{\syn{w}}\smcl\syn{x})\smcl\ol{\syn{B}}(\syn{x}\smcl\ol{\syn{z}})\vdash
    \left(\ideind{\syn{I}}\{\syn{\mu}\}\right)[\syn{x}/\syn{y}]\psb{\refl_{\syn{I}}(\syn{x})/\syn{p}}\equiv \syn{\mu}:
    \syn{\gamma}(\syn{w}_0\smcl\syn{z}_m)
    }
    \and
    \inferrule*[right=$\ides$-Comp$\eta$]
    {\ol{\syn{w}}:\ol{\syn{J}}\smcl\syn{x}:\syn{I}\smcl\syn{y}:\syn{I}\smcl\ol{\syn{z}}:\ol{\syn{K}}\mid
    \ol{\syn{A}}(\ol{\syn{w}}\smcl\syn{x})\smcl\syn{p}:\syn{x}\ide{\syn{I}}\syn{y}\smcl
    \ol{\syn{B}}(\syn{y}\smcl\ol{\syn{z}})\vdash
    \syn{\nu}:\syn{\gamma}(\syn{w}_0\smcl\syn{z}_m)}
    {\ol{\syn{w}}:\ol{\syn{J}}\smcl\syn{x}:\syn{I}\smcl\syn{y}:\syn{I}\smcl\ol{\syn{z}}:\ol{\syn{K}}\mid
    \ol{\syn{A}}(\ol{\syn{w}}\smcl\syn{x})\smcl\syn{p}:\syn{x}\ide{\syn{I}}\syn{y}\smcl
    \ol{\syn{B}}(\syn{y}\smcl\ol{\syn{z}})\vdash
    \ideind{\syn{I}}\{\syn{\nu}[\syn{x}/\syn{y}]\psb{\refl_{\syn{I}}(\syn{x})/\syn{p}}\}\equiv \syn{\nu}:\syn{\gamma}(\syn{w}_0\smcl\syn{z}_m)
    }
\end{mathparpagebreakable}
\myparagraph{Unit protype meets product type.}\ 
\label{sec:unit-protype-meets-product-type}
\begin{mathparpagebreakable}
    \goodbreak
    \inferrule*[right=$\ides$-$\top$]
    {\ }
    {\cdot\smcl\cdot\vdash \exc_{\ides,\top} : \langle \rangle\ide{\syn{1}}\langle \rangle\ccong \top}
    \and
    \inferrule*[right=$\ides$-$\land$]
    {\syn{I}\ \textsf{type}\\ \syn{J}\ \textsf{type}}
    {\syn{x}:\syn{I},\syn{y}:\syn{J}\smcl\syn{x}':\syn{I},\syn{y}':\syn{J}
    \vdash \exc_{\ides,\land}:\langle\syn{x},\syn{y}\rangle\ide{\syn{I}\times\syn{J}}\langle\syn{x'},\syn{y'}\rangle
    \ccong \syn{x}\ide{\syn{I}}\syn{x}'\land\syn{y}\ide{\syn{J}}\syn{y'}}
    \and
    \inferrule*
    {\syn{I}\ \textsf{type}\\ \syn{J}\ \textsf{type}}
    {\syn{x}:\syn{I},\syn{y}:\syn{J}\smcl\syn{x}':\syn{I},\syn{y}':\syn{J}
    \mid \syn{a}:\langle\syn{x},\syn{y}\rangle\ide{\syn{I}\times\syn{J}}\langle\syn{x}',\syn{y'}\rangle
    \vdash \exc_{\ides,\land}\{\syn{a}\}\equiv \ind_{\ide{\syn{I}\times\syn{J}}}\{\left\langle\refl_{\syn{I}}(\syn{x}),
    \refl_{\syn{J}}(\syn{y})\right\rangle\}:\syn{x}\ide{\syn{I}}\syn{x'}\land\syn{y}\ide{\syn{J}}\syn{y'}}
    \and
    \text{where}\quad
    \inferrule*
    {\syn{x}:\syn{I},\syn{y}:\syn{J}
    \mid \left\langle\refl_{\syn{I}}(\syn{x}),
    \refl_{\syn{J}}(\syn{y})\right\rangle:\syn{x}\ide{\syn{I}}\syn{x'}\land\syn{y}\ide{\syn{J}}\syn{y'}}
    {\syn{x}:\syn{I},\syn{y}:\syn{J}\smcl\syn{x}':\syn{I},\syn{y}':\syn{J}
    \mid \syn{a}:\langle\syn{x},\syn{y}\rangle\ide{\syn{I}\times\syn{J}}\langle\syn{x}',\syn{y'}\rangle
    \vdash \ind_{\ide{\syn{I}\times\syn{J}}}\{\left\langle\refl_{\syn{I}}(\syn{x}),
    \refl_{\syn{J}}(\syn{y})\right\rangle\}:\syn{x}\ide{\syn{I}}\syn{x'}\land\syn{y}\ide{\syn{J}}\syn{y'}
    }
\end{mathparpagebreakable}
\myparagraph{Composition protype.}\ 
\label{sec:composition-protype}
\begin{mathparpagebreakable}
    \goodbreak
    \inferrule*[right=$\odot$-Form]
    {\syn{w}:\syn{I}\smcl\syn{x}:\syn{J}\vdash \syn{\alpha}(\syn{w}\smcl\syn{x})\ \textsf{protype} \\ 
    \syn{x}:\syn{J}\smcl\syn{y}:\syn{K}\vdash \syn{\beta}(\syn{x}\smcl\syn{y})\ \textsf{protype}}
    {\syn{w}:\syn{I}\smcl\syn{y}:\syn{K}\vdash \syn{\alpha}(\syn{w}\smcl\syn{x})\odot_{\syn{x}: \syn{J}}\syn{\beta}(\syn{x}\smcl\syn{y})\ \textsf{protype}}
    \and
    \inferrule*[right=$\odot$-Intro]
    {\syn{w}:\syn{I}\smcl\syn{x}:\syn{J}\vdash \syn{\alpha}(\syn{w}\smcl\syn{x})\ \textsf{protype} \\
    \syn{x}:\syn{J}\smcl\syn{y}:\syn{K}\vdash \syn{\beta}(\syn{x}\smcl\syn{y})\ \textsf{protype} }
    {\syn{w}:\syn{I}\smcl\syn{x}:\syn{J}\smcl\syn{y}:\syn{K}\mid \syn{a}:\syn{\alpha}(\syn{w}\smcl\syn{x})\smcl\syn{b}:\syn{\beta}(\syn{x}\smcl\syn{y})\vdash
    \syn{a}\odot\syn{b}:\syn{\alpha}(\syn{w}\smcl\syn{x})\odot_{\syn{x}:\syn{J}}\syn{\beta}(\syn{x}\smcl\syn{y})}
    \and
    \inferrule*[right=$\odot$-Elim]
    {
    \ol{\syn{v}}:\ol{\syn{H}}\smcl\syn{w}:\syn{I}\smcl\syn{x}:\syn{J}\smcl\syn{y}:\syn{K}\smcl\ol{\syn{z}}:\ol{\syn{L}}
    \mid \ol{\syn{C}}(\ol{\syn{v}}\smcl\syn{w})\smcl\syn{a}:\syn{\alpha}(\syn{w}\smcl\syn{x})\smcl\syn{b}:\syn{\beta}(\syn{x}\smcl\syn{y})\smcl
    \ol{\syn{D}}(\syn{y}\smcl\ol{\syn{z}})\vdash \syn{\mu}:\syn{\gamma}(\syn{v}_0\smcl\syn{z}_m)}
    {\ol{\syn{v}}:\ol{\syn{H}}\smcl\syn{w}:\syn{I}\smcl\syn{y}:\syn{K}\smcl\ol{\syn{z}}:\ol{\syn{L}}\mid
    \ol{\syn{C}}(\ol{\syn{v}}\smcl\syn{w})\smcl\syn{p}:\syn{\alpha}(\syn{w}\smcl\syn{x})\odot_{\syn{x}:\syn{J}}\syn{\beta}(\syn{x}\smcl\syn{y})\smcl
    \ol{\syn{D}}(\syn{y}\smcl\ol{\syn{z}})\vdash
    \compind{\syn{\alpha}}{\syn{\beta}}\{\syn{\mu}\}:\syn{\gamma}(\syn{v}_0\smcl\syn{z}_m)}
    \and
    \inferrule*[right=$\odot$-Comp$\beta$]
    {
    \ol{\syn{v}}:\ol{\syn{H}}\smcl\syn{w}:\syn{I}\smcl\syn{x}:\syn{J}\smcl\syn{y}:\syn{K}\smcl\ol{\syn{z}}:\ol{\syn{L}}
    \mid \ol{\syn{C}}(\ol{\syn{v}}\smcl\syn{w})\smcl\syn{\alpha}(\syn{w}\smcl\syn{x})\smcl\syn{\beta}(\syn{x}\smcl\syn{y})\smcl
    \ol{\syn{D}}(\syn{y}\smcl\ol{\syn{z}})\vdash \syn{\mu}:\syn{\gamma}(\syn{v}_0\smcl\syn{z}_m)}
    {{\begin{array}{r}{\ol{\syn{v}}:\ol{\syn{H}}\smcl\syn{w}:\syn{I}\smcl\syn{x}:\syn{J}\smcl\syn{y}:\syn{K}\smcl\ol{\syn{z}}:\ol{\syn{L}}\mid
    \ol{\syn{C}}(\ol{\syn{v}}\smcl\syn{w})\smcl\syn{a}:\syn{\alpha}(\syn{w}\smcl\syn{x})\smcl
    \syn{b}:\syn{\beta}(\syn{x}\smcl\syn{y})\smcl\ol{\syn{D}}(\syn{y}\smcl\ol{\syn{z}})\quad}\\{\vdash
    \left(\compind{\syn{\alpha}}{\syn{\beta}}\{\syn{\mu}\}\right)\psb{\syn{a}\odot\syn{b}/\syn{p}}\equiv \syn{\mu}:
    \syn{\gamma}(\syn{v}_0\smcl\syn{z}_m)}
    \end{array}}
    }
    \and
    \inferrule*[right=$\odot$-Comp$\eta$]
    {\ol{\syn{v}}:\ol{\syn{H}}\smcl\syn{w}:\syn{I}\smcl\syn{y}:\syn{K}\smcl\ol{\syn{z}}:\ol{\syn{L}}
    \mid \ol{\syn{C}}(\ol{\syn{v}}\smcl\syn{w})\smcl\syn{p}:\syn{\alpha}(\syn{w}\smcl\syn{x})\odot_{\syn{x}:\syn{J}}\syn{\beta}(\syn{x}\smcl\syn{y})\smcl
    \ol{\syn{D}}(\syn{y}\smcl\ol{\syn{z}})\vdash 
    \syn{\nu}:\syn{\gamma}(\syn{v}_0\smcl\syn{z}_m)}
    {\ol{\syn{v}}:\ol{\syn{H}}\smcl\syn{w}:\syn{I}\smcl\syn{y}:\syn{K}\smcl\ol{\syn{z}}:\ol{\syn{L}}\mid
    \ol{\syn{C}}(\ol{\syn{v}}\smcl\syn{w})\smcl\syn{p}:\syn{\alpha}(\syn{w}\smcl\syn{x})\odot_{\syn{x}:\syn{J}}\syn{\beta}(\syn{x}\smcl\syn{y})\smcl
    \ol{\syn{D}}(\syn{y}\smcl\ol{\syn{z}})\vdash
    \compind{\syn{\alpha}}{\syn{\beta}}\{\left(\syn{\nu}\psb{\syn{a}\odot\syn{b}/\syn{p}}\right)\}\equiv \syn{\nu}:\syn{\gamma}(\syn{v}_0\smcl\syn{z}_m)
    }
\end{mathparpagebreakable}
\myparagraph{Composition protype meets product type.}\ 
\label{sec:compo-protype-meets-product-type}
\begin{mathparpagebreakable}
    \goodbreak
    \inferrule*[right=$\odot$-$\top$]
    {\ }
    {\cdot\smcl\cdot\vdash \exc_{\odot,\top}:\top\odot_{\langle \rangle:\cdot}\top\ccong \top}
    \and
    \inferrule*[right=$\odot$-$\land$]
    {\syn{x}:\syn{I}\smcl\syn{y}:\syn{J}\vdash \syn{\alpha}(\syn{x}\smcl\syn{y})\ \textsf{protype}\\
    \syn{y}:\syn{J}\smcl\syn{z}:\syn{K}\vdash \syn{\beta}(\syn{y}\smcl\syn{z})\ \textsf{protype}\\
    \syn{u}:\syn{L}\smcl\syn{v}:\syn{M}\vdash \syn{\gamma}(\syn{u}\smcl\syn{v})\ \textsf{protype}\\
    \syn{v}:\syn{M}\smcl\syn{w}:\syn{N}\vdash \syn{\delta}(\syn{v}\smcl\syn{w})\ \textsf{protype}}
    {
    {\begin{array}{r}
            {\syn{x}:\syn{I},\syn{u}:\syn{L}\smcl\syn{z}:\syn{K},\syn{w}:\syn{N}\vdash
    \exc_{\odot,\land}:
    \left(\syn{\alpha}(\syn{x}\smcl\syn{y})\land\syn{\gamma}(\syn{u}\smcl\syn{v})\right)
    \odot_{\langle\syn{y},\syn{v}\rangle:\syn{J}\times\syn{M}}
    \left(\syn{\beta}(\syn{y}\smcl\syn{z})\land\syn{\delta}(\syn{v}\smcl\syn{w})\right)}
    \quad
    \\
    {\ccong
    \left(\syn{\alpha}(\syn{x}\smcl\syn{y})\odot_{\syn{y}:\syn{J}}\syn{\beta}(\syn{y}\smcl\syn{z})\right)
    \land
    \left(\syn{\gamma}(\syn{u}\smcl\syn{v})\odot_{\syn{v}:\syn{M}}\syn{\delta}(\syn{v}\smcl\syn{w})\right)}
    \end{array}
    }
    }
    \and
    \inferrule*
    {\syn{x}:\syn{I}\smcl\syn{y}:\syn{J}\vdash \syn{\alpha}(\syn{x}\smcl\syn{y})\ \textsf{protype}\\
    \syn{y}:\syn{J}\smcl\syn{z}:\syn{K}\vdash \syn{\beta}(\syn{y}\smcl\syn{z})\ \textsf{protype}\\
    \syn{u}:\syn{L}\smcl\syn{v}:\syn{M}\vdash \syn{\gamma}(\syn{u}\smcl\syn{v})\ \textsf{protype}\\
    \syn{v}:\syn{M}\smcl\syn{w}:\syn{N}\vdash \syn{\delta}(\syn{v}\smcl\syn{w})\ \textsf{protype}}
    {
    {\begin{array}{l}    
    {\syn{x}:\syn{I},\syn{u}:\syn{L}\smcl\syn{z}:\syn{K},\syn{w}:\syn{N}\mid 
        \syn{e}:
        \left(\syn{\alpha}(\syn{x}\smcl\syn{y})\land\syn{\gamma}(\syn{u}\smcl\syn{v})\right)
        \odot_{\langle\syn{y},\syn{v}\rangle:\syn{J}\times\syn{M}}
        \left(\syn{\beta}(\syn{y}\smcl\syn{z})\land\syn{\delta}(\syn{v}\smcl\syn{w})\right)}\\
        {\qquad
        \vdash
        \exc_{\odot,\land}\{\syn{e}\}
        \equiv
        \ind_{\odot_{\syn{\alpha}\land\syn{\gamma},\syn{\beta}\land\syn{\delta}}}\left\{
            \left\langle\syn{\pi}_0\{\syn{a}\}\odot\syn{\pi}_0\{\syn{b}\},\syn{\pi}_1\{\syn{a}\}\odot\syn{\pi}_1\{\syn{b}\}\right\rangle
            \right\}:\left(\syn{\alpha}(\syn{x}\smcl\syn{y})\odot_{\syn{y}:\syn{J}}\syn{\beta}(\syn{y}\smcl\syn{z})\right)
            \land
            \left(\syn{\gamma}(\syn{u}\smcl\syn{v})\odot_{\syn{v}:\syn{M}}\syn{\delta}(\syn{v}\smcl\syn{w})\right)}
    \end{array}}
    }
    \and
    \text{where}\quad
    \inferrule*
    {\syn{x}:\syn{I}\smcl\syn{u}:\syn{L}\smcl\syn{y}:\syn{J}\smcl\syn{v}:\syn{M}\smcl
    \syn{z}:\syn{K}\smcl\syn{w}:\syn{N}\mid \syn{a}:\syn{\alpha}(\syn{x}\smcl\syn{y})\land
    \syn{\gamma}(\syn{u}\smcl\syn{v})\smcl\syn{b}:\syn{\beta}(\syn{y}\smcl\syn{z})\land
    \syn{\delta}(\syn{v}\smcl\syn{w})\\
    \qquad
    \vdash 
    \left\langle\syn{\pi}_0\{\syn{a}\}\odot\syn{\pi}_0\{\syn{b}\},\syn{\pi}_1\{\syn{a}\}\odot\syn{\pi}_1\{\syn{b}\}\right\rangle
    :\left(\syn{\alpha}(\syn{x}\smcl\syn{y})\odot_{\syn{y}:\syn{J}}\syn{\beta}(\syn{y}\smcl\syn{z})\right)
    \land
    \left(\syn{\gamma}(\syn{u}\smcl\syn{v})\odot_{\syn{v}:\syn{M}}\syn{\delta}(\syn{v}\smcl\syn{w})\right)
    }
    {
    {\begin{array}{l}    
    \syn{x}:\syn{I}\smcl\syn{u}:\syn{L}\smcl\syn{z}:\syn{K}\smcl\syn{w}:\syn{N}\mid
    \syn{e}:\left(\syn{\alpha}(\syn{x}\smcl\syn{y})\land\syn{\gamma}(\syn{u}\smcl\syn{v})\right)
    \odot_{\langle\syn{y},\syn{v}\rangle:\syn{J}\times\syn{M}}
    \left(\syn{\beta}(\syn{y}\smcl\syn{z})\land\syn{\delta}(\syn{v}\smcl\syn{w})\right)\\
    \qquad 
    \vdash
    \ind_{\odot_{\syn{\alpha}\land\syn{\gamma},\syn{\beta}\land\syn{\delta}}}\left\{
    \left\langle\syn{\pi}_0\{\syn{a}\}\odot\syn{\pi}_0\{\syn{b}\},\syn{\pi}_1\{\syn{a}\}\odot\syn{\pi}_1\{\syn{b}\}\right\rangle
    \right\}:\left(\syn{\alpha}(\syn{x}\smcl\syn{y})\odot_{\syn{y}:\syn{J}}\syn{\beta}(\syn{y}\smcl\syn{z})\right)
    \land
    \left(\syn{\gamma}(\syn{u}\smcl\syn{v})\odot_{\syn{v}:\syn{M}}\syn{\delta}(\syn{v}\smcl\syn{w})\right)
    \end{array}}
    }
\end{mathparpagebreakable}
\myparagraph{Filler protype.}\ 
\label{sec:filler-protype}
\begin{mathparpagebreakable}
    \goodbreak
    \inferrule*[right=$\triangleright$-Form]
    {\syn{w}:\syn{I}\smcl\syn{x}:\syn{J}\vdash \syn{\alpha}(\syn{w}\smcl\syn{x})\ \textsf{protype} \\
    \syn{w}:\syn{I}\smcl\syn{y}:\syn{K}\vdash \syn{\beta}(\syn{w}\smcl\syn{y})\ \textsf{protype}}
    {\syn{x}:\syn{J}\smcl\syn{y}:\syn{K}\vdash \syn{\alpha}(\syn{w}\smcl\syn{x})\triangleright_{\syn{w}:\syn{I}}\syn{\beta}(\syn{w}\smcl\syn{y})\ \textsf{protype}}
    \and
    \inferrule*[right=$\triangleright$-Intro]
    {\syn{w}:\syn{I}\smcl\syn{x}:\syn{J}\smcl\ol{\syn{y}}:\ol{\syn{L}}\mid \syn{a}:\syn{\alpha}(\syn{w}\smcl\syn{x})\smcl\ol{\syn{C}}(\syn{x}\smcl\ol{\syn{y}})\vdash \syn{\mu}:\syn{\beta}(\syn{w}\smcl\syn{y}_m)}
    {\syn{x}:\syn{J}\smcl\ol{\syn{y}}:\ol{\syn{L}}\mid \ol{\syn{C}}(\syn{x}\smcl\ol{\syn{y}})\vdash \ind_{\triangleright_{\syn{\alpha},\syn{\beta}}}\{\syn{\mu}\}:\syn{\alpha}(\syn{w}\smcl\syn{x})\triangleright_{\syn{w}:\syn{I}}\syn{\beta}(\syn{w}\smcl\syn{y}_m)}
    \and
    \inferrule*[right=$\triangleright$-Elim]
    {\syn{w}:\syn{I}\smcl\syn{x}:\syn{J}\vdash \syn{\alpha}(\syn{w}\smcl\syn{x})\ \textsf{protype} \\
    \syn{w}:\syn{I}\smcl\syn{y}:\syn{K}\vdash \syn{\beta}(\syn{w}\smcl\syn{y})\ \textsf{protype}}
    {\syn{w}:\syn{I}\smcl\syn{x}:\syn{J}\smcl\syn{y}:\syn{K}\mid \syn{a}:\syn{\alpha}(\syn{w}\smcl\syn{x})\smcl \syn{e}:\syn{\alpha}(\syn{w}\smcl\syn{x})\triangleright_{\syn{w}:\syn{I}}\syn{\beta}(\syn{w}\smcl\syn{y})\vdash \syn{a}\rbl\syn{e}:\syn{\beta}(\syn{w}\smcl\syn{y})}
    \and
    \inferrule*[right=$\triangleright$-Comp$\beta$]
    {\syn{w}:\syn{I}\smcl\syn{x}:\syn{J}\smcl\ol{\syn{y}}:\ol{\syn{L}}\mid \syn{a}:\syn{\alpha}(\syn{w}\smcl\syn{x})\smcl\ol{\syn{C}}(\syn{x}\smcl\ol{\syn{y}})\vdash 
    \syn{\mu}:\syn{\beta}(\syn{w}\smcl\syn{y}_m)}
    {\syn{w}:\syn{I}\smcl\syn{x}:\syn{J}\smcl\ol{\syn{y}}:\ol{\syn{L}}\mid 
    \syn{a}:\syn{\alpha}(\syn{w}\smcl\syn{x})\smcl\ol{\syn{C}}(\syn{x}\smcl\ol{\syn{y}})\vdash
    \syn{a}\rbl\left(\ind_{\triangleright_{\syn{\alpha},\syn{\beta}}}\{\syn{\mu}\}\right)\equiv \syn{\mu}:\syn{\beta}(\syn{w}\smcl\syn{y}_m)}
    \and
    \inferrule*[right=$\triangleright$-Comp$\eta$]
    {\syn{x}:\syn{J}\smcl\ol{\syn{y}}:\ol{\syn{L}}\mid\ol{\syn{C}}(\syn{x}\smcl\ol{\syn{y}})\vdash
    \syn{\nu}:\syn{\alpha}(\syn{w}\smcl\syn{x})\triangleright_{\syn{w}:\syn{I}}\syn{\beta}(\syn{w}\smcl\syn{y}_m)}
    {\syn{x}:\syn{J}\smcl\ol{\syn{y}}:\ol{\syn{L}}\mid \ol{\syn{C}}(\syn{x}\smcl\ol{\syn{y}})\vdash
    \ind_{\triangleright_{\syn{\alpha},\syn{\beta}}}\left\{\syn{a}\rbl\syn{\nu}\right\}\equiv \syn{\nu}:\syn{\beta}(\syn{w}\smcl\syn{y}_m)}
\end{mathparpagebreakable}
\begin{mathparpagebreakable}
    \goodbreak
    \inferrule*[right=$\triangleleft$-Form]
    {\syn{y}:\syn{J}\smcl\syn{z}:\syn{K}\vdash \syn{\alpha}(\syn{y}\smcl\syn{z})\ \textsf{protype} \\
    \syn{x}:\syn{I}\smcl\syn{z}:\syn{K}\vdash \syn{\beta}(\syn{x}\smcl\syn{z})\ \textsf{protype}}
    {\syn{x}:\syn{I}\smcl\syn{y}:\syn{J}\vdash \syn{\beta}(\syn{x}\smcl\syn{z})\triangleleft_{\syn{z}:\syn{K}}
    \syn{\alpha}(\syn{y}\smcl\syn{z})\ \textsf{protype}}
    \and
    \inferrule*[right=$\triangleleft$-Intro]
    {\ol{\syn{x}}:\ol{\syn{J}}\smcl\syn{y}:\syn{J}\smcl\syn{z}:\syn{K}\mid \ol{\syn{C}}(\ol{\syn{x}}\smcl\syn{y})\smcl
    \syn{a}:\syn{\alpha}(\syn{y}\smcl\syn{z})\vdash \syn{\mu}:\syn{\beta}(\syn{x}\smcl\syn{z})}
    {\ol{\syn{x}}:\ol{\syn{J}}\smcl\syn{y}:\syn{J}\mid \ol{\syn{C}}(\ol{\syn{x}}\smcl\syn{y})\vdash
    \ind_{\triangleleft_{\syn{\alpha},\syn{\beta}}}\{\syn{\mu}\}:\syn{\beta}(\syn{x}\smcl\syn{z})\triangleleft_{\syn{z}:\syn{K}}\syn{\alpha}(\syn{y}\smcl\syn{z})}
    \and
    \inferrule*[right=$\triangleleft$-Elim]
    {\syn{x}:\syn{I}\smcl\syn{y}:\syn{J}\vdash \syn{\beta}(\syn{x}\smcl\syn{z})\ \textsf{protype} \\
    \syn{y}:\syn{J}\smcl\syn{z}:\syn{K}\vdash \syn{\alpha}(\syn{y}\smcl\syn{z})\ \textsf{protype}}
    {\syn{x}:\syn{I}\smcl\syn{y}:\syn{J}\smcl\syn{z}:\syn{K}\mid \syn{a}:\syn{\beta}(\syn{x}\smcl\syn{z})\smcl \syn{e}:\syn{\beta}(\syn{x}\smcl\syn{z})\triangleleft_{\syn{z}:\syn{K}}\syn{\alpha}(\syn{y}\smcl\syn{z})\vdash \syn{a}\lbl\syn{e}:\syn{\alpha}(\syn{y}\smcl\syn{z})}
    \and
    \inferrule*[right=$\triangleleft$-Comp$\beta$]
    {\syn{x}:\syn{I}\smcl\syn{y}:\syn{J}\smcl\ol{\syn{z}}:\ol{\syn{L}}\mid \syn{a}:\syn{\beta}(\syn{x}\smcl\syn{z})\smcl\ol{\syn{C}}(\syn{x}\smcl\ol{\syn{z}})\vdash
    \syn{\mu}:\syn{\alpha}(\syn{y}\smcl\syn{z}_m)}
    {\syn{x}:\syn{I}\smcl\syn{y}:\syn{J}\smcl\ol{\syn{z}}:\ol{\syn{L}}\mid \syn{a}:\syn{\beta}(\syn{x}\smcl\syn{z})\smcl\ol{\syn{C}}(\syn{x}\smcl\ol{\syn{z}})\vdash
    \syn{a}\lbl\left(\ind_{\triangleleft_{\syn{\alpha},\syn{\beta}}}\{\syn{\mu}\}\right)\equiv \syn{\mu}:\syn{\alpha}(\syn{y}\smcl\syn{z}_m)}
    \and
    \inferrule*[right=$\triangleleft$-Comp$\eta$]
    {\syn{y}:\syn{J}\smcl\ol{\syn{z}}:\ol{\syn{L}}\mid\ol{\syn{C}}(\syn{y}\smcl\ol{\syn{z}})\vdash
    \syn{\nu}:\syn{\beta}(\syn{x}\smcl\syn{z})\triangleleft_{\syn{z}:\syn{K}}\syn{\alpha}(\syn{y}\smcl\syn{z})}
    {\syn{y}:\syn{J}\smcl\ol{\syn{z}}:\ol{\syn{L}}\mid\ol{\syn{C}}(\syn{y}\smcl\ol{\syn{z}})\vdash
    \ind_{\triangleleft_{\syn{\alpha},\syn{\beta}}}\left\{\syn{a}\lbl\syn{\nu}\right\}\equiv \syn{\nu}:\syn{\alpha}(\syn{y}\smcl\syn{z}_m)}
\end{mathparpagebreakable}
\myparagraph{Filler protype meets product type.}\ 
\label{sec:filler-protype-meets-product-type}
\begin{mathparpagebreakable}
    \goodbreak
    \inferrule*[right=$\triangleright$-$\top$]
    {\ }
    {\cdot\smcl\cdot\mid \exc_{\triangleright,\top}:\top\triangleright_{\cdot}\top\ccong \top}
    \and
    \inferrule*[right=$\triangleright$-$\land$]
    {\syn{x}:\syn{I}\smcl\syn{y}:\syn{J}\vdash \syn{\alpha}(\syn{x}\smcl\syn{y})\ \textsf{protype}\\
    \syn{x}:\syn{I}\smcl\syn{z}:\syn{K}\vdash \syn{\beta}(\syn{x}\smcl\syn{z})\ \textsf{protype}\\
    \syn{u}:\syn{L}\smcl\syn{v}:\syn{M}\vdash \syn{\gamma}(\syn{u}\smcl\syn{v})\ \textsf{protype}\\
    \syn{u}:\syn{L}\smcl\syn{w}:\syn{N}\vdash \syn{\delta}(\syn{v}\smcl\syn{w})\ \textsf{protype}}
    {
    {\begin{array}{l}    
    \syn{y}:\syn{J},\syn{v}:\syn{M}\smcl\syn{z}:\syn{K},\syn{w}:\syn{N}\vdash
    \exc_{\triangleright,\land}:
    \left(\syn{\alpha}(\syn{x}\smcl\syn{y})\triangleright_{\syn{x}:\syn{I}}\syn{\beta}(\syn{x}\smcl\syn{z})\right)
    \land
    \left(\syn{\gamma}(\syn{u}\smcl\syn{v})\triangleright_{\syn{u}:\syn{L}}\syn{\delta}(\syn{v}\smcl\syn{w})\right)\\
    \ccong
    \left(\syn{\alpha}(\syn{x}\smcl\syn{y})\land\syn{\gamma}(\syn{u}\smcl\syn{v})\right)
    \triangleright_{\syn{x}:\syn{I},\syn{u}:\syn{L}}
    \left(\syn{\beta}(\syn{x}\smcl\syn{z})\land\syn{\delta}(\syn{v}\smcl\syn{w})\right)
    \end{array}}
    }
    \and
    \inferrule*[right=$\triangleright$-$\land$-canon]
    {\syn{x}:\syn{I}\smcl\syn{y}:\syn{J}\vdash \syn{\alpha}(\syn{x}\smcl\syn{y})\ \textsf{protype}\\
    \syn{x}:\syn{I}\smcl\syn{z}:\syn{K}\vdash \syn{\beta}(\syn{x}\smcl\syn{z})\ \textsf{protype}\\
    \syn{u}:\syn{L}\smcl\syn{v}:\syn{M}\vdash \syn{\gamma}(\syn{u}\smcl\syn{v})\ \textsf{protype}\\
    \syn{u}:\syn{L}\smcl\syn{w}:\syn{N}\vdash \syn{\delta}(\syn{v}\smcl\syn{w})\ \textsf{protype}}
    {
        {\begin{array}{l}
        \syn{y}:\syn{J},\syn{v}:\syn{M}\smcl\syn{z}:\syn{K},\syn{w}:\syn{N}\mid 
        \syn{e}:
        \left(\syn{\alpha}(\syn{x}\smcl\syn{y})\triangleright_{\syn{x}:\syn{I}}\syn{\beta}(\syn{x}\smcl\syn{z})\right)
        \land
        \left(\syn{\gamma}(\syn{u}\smcl\syn{v})\triangleright_{\syn{u}:\syn{L}}\syn{\delta}(\syn{v}\smcl\syn{w})\right)\\
        \qquad
        \vdash
        \exc_{\triangleright,\land}\{\syn{e}\}
        \equiv
        \ind_{\triangleright_{\syn{\alpha\land\gamma},\syn{\beta\land\delta}}}
        \left\{
        \left\langle\syn{\pi}_0\{\syn{a}\}\rbl\syn{\pi}_0(\syn{e}),\syn{\pi}_1\{\syn{a}\}\rbl\syn{\pi}_1(\syn{e})\right\rangle
        \right\}:\left(\syn{\alpha}(\syn{x}\smcl\syn{y})\land\syn{\gamma}(\syn{u}\smcl\syn{v})\right)
        \triangleright_{\syn{x}:\syn{I},\syn{u}:\syn{L}}
        \left(\syn{\beta}(\syn{x}\smcl\syn{z})\land\syn{\delta}(\syn{v}\smcl\syn{w})\right)
        \end{array}}
    }
    \and
    \text{where}\quad
    \inferrule*
    {
    {\begin{array}{l}
    \syn{x}:\syn{I},\syn{u}:\syn{L},\syn{y}:\syn{J},\syn{v}:\syn{M},\syn{z}:\syn{K},\syn{w}:\syn{N}\mid
    \syn{a}:\left(\syn{\alpha}(\syn{x}\smcl\syn{y})\land\syn{\gamma}(\syn{u}\smcl\syn{v})\right)\smcl
    \syn{e}:\left(\syn{\alpha}(\syn{x}\smcl\syn{y})\triangleright_{\syn{x}:\syn{I}}\syn{\beta}(\syn{x}\smcl\syn{z})\right)
    \land
    \left(\syn{\gamma}(\syn{u}\smcl\syn{v})\triangleright_{\syn{u}:\syn{L}}\syn{\delta}(\syn{v}\smcl\syn{w})\right)\\
    \qquad
    \vdash
    \left\langle\syn{\pi}_0\{\syn{a}\}\rbl\syn{\pi}_0(\syn{e}),\syn{\pi}_1\{\syn{a}\}\rbl\syn{\pi}_1(\syn{e})\right\rangle:
    \left(\syn{\beta}(\syn{x}\smcl\syn{z})\land\syn{\delta}(\syn{v}\smcl\syn{w})\right)
    \end{array}}
    }
    {
    {\begin{array}{l}
    \syn{y}:\syn{J},\syn{v}:\syn{M}\smcl\syn{z}:\syn{K},\syn{w}:\syn{N}\mid
    \syn{e}:\left(\syn{\alpha}(\syn{x}\smcl\syn{y})\triangleright_{\syn{x}:\syn{I}}\syn{\beta}(\syn{x}\smcl\syn{z})\right)
    \land
    \left(\syn{\gamma}(\syn{u}\smcl\syn{v})\triangleright_{\syn{u}:\syn{L}}\syn{\delta}(\syn{v}\smcl\syn{w})\right)\\
    \qquad
    \vdash
    \ind_{\triangleright_{\syn{\alpha\land\gamma},\syn{\beta\land\delta}}}
    \left\{
    \left\langle\syn{\pi}_0\{\syn{a}\}\rbl\syn{\pi}_0(\syn{e}),\syn{\pi}_1\{\syn{a}\}\rbl\syn{\pi}_1(\syn{e})\right\rangle
    \right\}
    :\left(\syn{\alpha}(\syn{x}\smcl\syn{y})\land\syn{\gamma}(\syn{u}\smcl\syn{v})\right)
    \triangleright_{\syn{x}:\syn{I},\syn{u}:\syn{L}}
    \left(\syn{\beta}(\syn{x}\smcl\syn{z})\land\syn{\delta}(\syn{v}\smcl\syn{w})\right)
    \end{array}}
    }
\end{mathparpagebreakable}
\begin{mathparpagebreakable}
    \goodbreak
    \inferrule*[right=$\triangleleft$-$\top$]
    {\ }
    {\cdot\smcl\cdot\vdash 
    \exc_{\triangleleft,\top}:\top\triangleleft_{\cdot}\top\equiv \top}
    \and
    \inferrule*[right=$\triangleleft$-$\land$]
    {\syn{x}:\syn{I}\smcl\syn{z}:\syn{K}\vdash \syn{\alpha}(\syn{x}\smcl\syn{z})\ \textsf{protype}\\
    \syn{y}:\syn{J}\smcl\syn{z}:\syn{K}\vdash \syn{\beta}(\syn{y}\smcl\syn{z})\ \textsf{protype}\\
    \syn{u}:\syn{L}\smcl\syn{w}:\syn{N}\vdash \syn{\gamma}(\syn{u}\smcl\syn{w})\ \textsf{protype}\\
    \syn{v}:\syn{M}\smcl\syn{w}:\syn{N}\vdash \syn{\delta}(\syn{v}\smcl\syn{w})\ \textsf{protype}}
    {    
    {\begin{array}{l}
    \syn{x}:\syn{I},\syn{u}:\syn{L}\smcl\syn{y}:\syn{J},\syn{v}:\syn{M}\vdash
    \exc_{\triangleleft,\land}:
    \left(\syn{\alpha}(\syn{x}\smcl\syn{z})\triangleleft_{\syn{z}:\syn{K}}\syn{\beta}(\syn{y}\smcl\syn{z})\right)
    \land
    \left(\syn{\gamma}(\syn{u}\smcl\syn{w})\triangleleft_{\syn{w}:\syn{N}}\syn{\delta}(\syn{v}\smcl\syn{w})\right)\\
    \ccong
    \left(\syn{\alpha}(\syn{x}\smcl\syn{z})\land\syn{\gamma}(\syn{u}\smcl\syn{w})\right)
    \triangleleft_{\syn{z}:\syn{K},\syn{w}:\syn{N}}
    \left(\syn{\beta}(\syn{y}\smcl\syn{z})\land\syn{\delta}(\syn{v}\smcl\syn{w})\right)
    \end{array}}
    }
    \and
    \inferrule*[right=$\triangleleft$-$\land$-canon]
    {\syn{x}:\syn{I}\smcl\syn{z}:\syn{K}\vdash \syn{\alpha}(\syn{x}\smcl\syn{z})\ \textsf{protype}\\
    \syn{y}:\syn{J}\smcl\syn{z}:\syn{K}\vdash \syn{\beta}(\syn{y}\smcl\syn{z})\ \textsf{protype}\\
    \syn{u}:\syn{L}\smcl\syn{w}:\syn{N}\vdash \syn{\gamma}(\syn{u}\smcl\syn{w})\ \textsf{protype}\\
    \syn{v}:\syn{M}\smcl\syn{w}:\syn{N}\vdash \syn{\delta}(\syn{v}\smcl\syn{w})\ \textsf{protype}}
    {
    {\begin{array}{l}
        \syn{x}:\syn{I},\syn{u}:\syn{L}\smcl\syn{y}:\syn{J},\syn{v}:\syn{M}\mid 
        \syn{e}:
        \left(\syn{\alpha}(\syn{x}\smcl\syn{z})\triangleleft_{\syn{z}:\syn{K}}\syn{\beta}(\syn{y}\smcl\syn{z})\right)
        \land
        \left(\syn{\gamma}(\syn{u}\smcl\syn{w})\triangleleft_{\syn{w}:\syn{N}}\syn{\delta}(\syn{v}\smcl\syn{w})\right)\\
        \qquad
        \vdash
        \exc_{\triangleleft,\land}\{\syn{e}\}
        \equiv
        \ind_{\triangleleft_{\syn{\alpha\land\gamma},\syn{\beta\land\delta}}}
        \left\{
        \left\langle\syn{\pi}_0\{\syn{a}\}\lbl\syn{\pi}_0(\syn{e}),\syn{\pi}_1\{\syn{a}\}\lbl\syn{\pi}_1(\syn{e})\right\rangle
        \right\}:\left(\syn{\alpha}(\syn{x}\smcl\syn{z})\land\syn{\gamma}(\syn{u}\smcl\syn{w})\right)
        \triangleleft_{\syn{z}:\syn{K},\syn{w}:\syn{N}}
        \left(\syn{\beta}(\syn{y}\smcl\syn{z})\land\syn{\delta}(\syn{v}\smcl\syn{w})\right)
    \end{array}}
    }
    \and
    \text{where}\quad
    \inferrule*
    {\
    {\begin{array}{l}
    \syn{x}:\syn{I},\syn{u}:\syn{L}\smcl\syn{y}:\syn{J},\syn{v}:\syn{M},\syn{z}:\syn{K},\syn{w}:\syn{N}\mid
    \syn{a}:\left(\syn{\alpha}(\syn{x}\smcl\syn{z})\land\syn{\gamma}(\syn{u}\smcl\syn{w})\right)\smcl
    \syn{e}:\left(\syn{\alpha}(\syn{x}\smcl\syn{z})\triangleleft_{\syn{z}:\syn{K}}\syn{\beta}(\syn{y}\smcl\syn{z})\right)
    \land
    \left(\syn{\gamma}(\syn{u}\smcl\syn{w})\triangleleft_{\syn{w}:\syn{N}}\syn{\delta}(\syn{v}\smcl\syn{w})\right)\\
    \qquad
    \vdash
    \left\langle\syn{\pi}_0\{\syn{a}\}\lbl\syn{\pi}_0(\syn{e}),\syn{\pi}_1\{\syn{a}\}\lbl\syn{\pi}_1(\syn{e})\right\rangle:
    \left(\syn{\beta}(\syn{y}\smcl\syn{z})\land\syn{\delta}(\syn{v}\smcl\syn{w})\right)
    \end{array}}
    }
    {
    {\begin{array}{l}
    \syn{x}:\syn{I},\syn{u}:\syn{L}\smcl\syn{y}:\syn{J},\syn{v}:\syn{M}\mid
    \syn{e}:\left(\syn{\alpha}(\syn{x}\smcl\syn{z})\triangleleft_{\syn{z}:\syn{K}}\syn{\beta}(\syn{y}\smcl\syn{z})\right)
    \land
    \left(\syn{\gamma}(\syn{u}\smcl\syn{w})\triangleleft_{\syn{w}:\syn{N}}\syn{\delta}(\syn{v}\smcl\syn{w})\right)\\
    \qquad
    \vdash
    \ind_{\triangleleft_{\syn{\alpha\land\gamma},\syn{\beta\land\delta}}}
    \left\{
    \left\langle\syn{\pi}_0\{\syn{a}\}\lbl\syn{\pi}_0(\syn{e}),\syn{\pi}_1\{\syn{a}\}\lbl\syn{\pi}_1(\syn{e})\right\rangle
    \right\}
    :\left(\syn{\alpha}(\syn{x}\smcl\syn{z})\land\syn{\gamma}(\syn{u}\smcl\syn{w})\right)
    \triangleleft_{\syn{z}:\syn{K},\syn{w}:\syn{N}}
    \left(\syn{\beta}(\syn{y}\smcl\syn{z})\land\syn{\delta}(\syn{v}\smcl\syn{w})\right)
    \end{array}}
    }
\end{mathparpagebreakable}

\myparagraph{Comprehension type.}\ 
\label{sec:comprehension-type}
\begin{mathparpagebreakable}
    \goodbreak
        \inferrule*[right=$\cmpr{}$-Form]
        {\syn{x}:\syn{I}\smcl\syn{y}:\syn{J}\vdash \syn{\alpha}\ \textsf{protype}}
        {\cmpr{\syn{\alpha}} \ \textsf{type}}
        \and
        \inferrule*[right=$\cmpr{}$-Elim-$\ell$]
        {\syn{x}:\syn{I}\smcl\syn{y}:\syn{J}\vdash \syn{\alpha}\ \textsf{protype}}
        {\syn{w}:\cmpr{\syn{\alpha}}\vdash \syn{l}(\syn{w}):\syn{I}}
        \and
        \inferrule*[right=$\cmpr{}$-Elim-$r$]
        {\syn{x}:\syn{I}\smcl\syn{y}:\syn{J}\vdash \syn{\alpha}\ \textsf{protype}}
        {\syn{w}:\cmpr{\syn{\alpha}}\vdash \syn{r}(\syn{w}):\syn{J}}
        \and 
        \inferrule*[right=$\cmpr{}$-Elim-cell]
        {\syn{x}:\syn{I}\smcl\syn{y}:\syn{J}\vdash \syn{\alpha}\ \textsf{protype}}
        {\syn{w}:\cmpr{\syn{\alpha}}\mid \vdash \tabb_{\cmpr{\syn{\alpha}}}\{\syn{w}\}:\syn{\alpha}[\syn{l}(\syn{w})/\syn{x}\smcl\syn{r}(\syn{w})/\syn{y}]}
        \and
        \inferrule*[right=$\cmpr{}$-Intro]
        {\syn{x}:\syn{I}\smcl\syn{y}:\syn{J}\vdash \syn{\alpha}\ \textsf{protype} \\
        \syn{\Gamma}\vdash \syn{s}:\syn{I} \\
        \syn{\Gamma}\vdash \syn{t}:\syn{J} \\
        \syn{\Gamma}\mid \vdash \syn{\nu}:\syn{\alpha}[\syn{s}/\syn{x}\smcl\syn{t}/\syn{y}]}
        {\syn{\Gamma}\vdash \ind_{\cmpr{}}(\syn{s},\syn{t},\syn{\nu}):\cmpr{\syn{\alpha}}}
        \and
        \inferrule*[right=$\cmpr{}$-Comp-$\ell$]
        {\syn{\Gamma}\vdash \syn{s}:\syn{I} \\
        \syn{\Gamma}\vdash \syn{t}:\syn{J} \\
        \syn{\Gamma}\mid \vdash \syn{\nu}:\syn{\alpha}[\syn{s}/\syn{x}\smcl\syn{t}/\syn{y}]}
        {\syn{\Gamma}\vdash \syn{l}(\ind_{\cmpr{}}(\syn{s},\syn{t},\syn{\nu}))\equiv \syn{s}:\syn{I}}
        \and
        \inferrule*[right=$\cmpr{}$-Comp-$r$]
        {\syn{\Gamma}\vdash \syn{s}:\syn{I} \\
        \syn{\Gamma}\vdash \syn{t}:\syn{J} \\
        \syn{\Gamma}\mid \vdash \syn{\nu}:\syn{\alpha}[\syn{s}/\syn{x}\smcl\syn{t}/\syn{y}]}
        {\syn{\Gamma}\vdash \syn{r}(\ind_{\cmpr{}}(\syn{s},\syn{t},\syn{\nu}))\equiv \syn{t}:\syn{J}}
        \and 
        \inferrule*[right=$\cmpr{}$-Comp-$\beta$]
        {\syn{x}:\syn{I}\smcl\syn{y}:\syn{J}\vdash \syn{\alpha}\ \textsf{protype} \\
        \syn{\Gamma}\vdash \syn{s}:\syn{I} \\
        \syn{\Gamma}\vdash \syn{t}:\syn{J} \\
        \syn{\Gamma}\mid \vdash \syn{\nu}:\syn{\alpha}[\syn{s}/\syn{x}\smcl\syn{t}/\syn{y}]}
        {\syn{\Gamma}\vdash \tabb_{\cmpr{\syn{\alpha}}}\{\ind_{\cmpr{}}(\syn{s},\syn{t},\syn{\nu})\}\equiv \syn{\nu}:
        \syn{\alpha}[\syn{s}/\syn{x}\smcl\syn{t}/\syn{y}]}
        \and 
        \inferrule*[right=$\cmpr{}$-Comp-$\eta$]
        {\syn{x}:\syn{I}\smcl\syn{y}:\syn{J}\vdash \syn{\alpha}\ \textsf{protype}}
        {\syn{w}:\cmpr{\syn{\alpha}}\vdash \ind_{\cmpr{}}(\syn{l}(\syn{w}),\syn{r}(\syn{w}),\tabb_{\cmpr{\syn{\alpha}}}\{\syn{w}\}\
        \equiv \syn{w}:\cmpr{\syn{\alpha}}}
\end{mathparpagebreakable}
\myparagraph{Comprehension type meets unit protype.}\ 
\label{sec:comprehension-type-meets-unit-protype}
\begin{mathparpagebreakable}
    \goodbreak
    \inferrule*[right=$\cmpr{}$-Elim]
    {\syn{\Gamma}_0\vdash \syn{s}_0:\syn{I} \\
    \syn{\Gamma}_m\vdash \syn{s}_1:\syn{I} \\
    \syn{\Gamma}_0\vdash \syn{t}_0:\syn{J} \\
    \syn{\Gamma}_m\vdash \syn{t}_1:\syn{J} \\
    \syn{x}:\syn{I},\syn{y}:\syn{J}\vdash \syn{\alpha}(\syn{x},\syn{y})\ \textsf{protype} \\
    \syn{\Gamma}_0\mid\vdash \syn{\mu}_0:\syn{\alpha}(\syn{s}_0\smcl\syn{t}_0) \\
    \syn{\Gamma}_m\mid\vdash \syn{\mu}_1:\syn{\alpha}(\syn{s}_1\smcl\syn{t}_1)\\
    \ol{\syn{\Gamma}}\mid \syn{B} \vdash \syn{i}:\syn{s}_0\ide{\syn{I}}\syn{s}_1 \\
    \ol{\syn{\Gamma}}\mid \syn{B} \vdash \syn{j}:\syn{t}_0\ide{\syn{J}}\syn{t}_1\\
    \ol{\syn{\Gamma}}\mid \syn{B} \vdash \syn{i}\boxdot\syn{\mu_1}\equiv \syn{\mu}_0\boxdot\syn{j}}
    {\ol{\syn{\Gamma}}\mid \syn{B} \vdash \ind_{\cmpr{}}(\syn{i},\syn{j},\syn{\mu}_0,\syn{\mu}_1):\ind_{\cmpr{}}(\syn{s}_0,\syn{t}_0,\syn{\mu}_0)\ide{\cmpr{\syn{\alpha}}}\ind_{\cmpr{}}(\syn{s}_1,\syn{t}_1,\syn{\mu}_1)}
    \and
    \text{where}\ 
    \inferrule*
    {
    \inferrule*
    {\inferrule*
    {\syn{x}:\syn{I}\smcl\syn{y}:\syn{J}\mid \syn{a}:\syn{\alpha}(\syn{x}\smcl\syn{y})\vdash \syn{a}: \syn{\alpha}(\syn{x'}\smcl\syn{y})}
    {\syn{x}:\syn{I}\smcl\syn{x'}:\syn{I}\smcl\syn{y}:\syn{J}\mid \syn{p}:\syn{x}\ide{\syn{I}}\syn{x'}\smcl\syn{a}:\syn{\alpha}(\syn{x}\smcl\syn{y})\vdash \ind_{\ide{}}\{\syn{a}\}:\syn{\alpha}(\syn{x}\smcl\syn{y})}\\
    \syn{\Gamma}_0\vdash \syn{s}_0:\syn{I} \\
    \syn{\Gamma}_m\vdash \syn{s}_1:\syn{I} \\
    \syn{\Gamma}_m\vdash \syn{t}_1:\syn{J} 
    }
    {
    \ol{\syn{\Gamma}}\mid \syn{p}:\syn{s}_0\ide{\syn{I}}\syn{s}_1\smcl\syn{a}:\syn{\alpha}(\syn{s}_1\smcl\syn{t}_1)\vdash \ind_{\ide{}}\{\syn{a}\}[\syn{s}_1/\syn{x'}\smcl\syn{t}_1/\syn{y}]
    :\syn{\alpha}(\syn{s}_0\smcl\syn{t}_1)
    }\\
    \ol{\syn{\Gamma}}\mid \syn{B} \vdash \syn{i}:\syn{s}_0\ide{\syn{I}}\syn{s}_1\\
    \syn{\Gamma}_m\mid\vdash \syn{\mu}_1:\syn{\alpha}(\syn{s}_1\smcl\syn{t}_1)\\
    }
    {
    \ol{\syn{\Gamma}}\mid \syn{B} \vdash \syn{i}\boxdot\syn{\mu}_1\colequiv \ind_{\ide{}}\{\syn{a}\}[\syn{s}_1/\syn{x'}\smcl\syn{t}_1/\syn{y}]\psb{ \syn{i}/\syn{p}:\syn{s}_0\ide{\syn{I}}\syn{s}_1\smcl\syn{\mu_1}/\syn{a}:\syn{\alpha}(\syn{s}_1\smcl\syn{t}_1)}:
    \syn{\alpha}(\syn{s}_0\smcl\syn{t}_1)
    }
    \and
    \text{and similarly for $\syn{\mu}_0\boxdot\syn{j}$.}
    \and
    \inferrule*[right=$\cmpr{}$-Comp]
    {\syn{\Gamma}_0\vdash \syn{s}_0:\syn{I} \\
    \syn{\Gamma}_m\vdash \syn{s}_1:\syn{I} \\
    \syn{\Gamma}_0\vdash \syn{t}_0:\syn{J} \\
    \syn{\Gamma}_m\vdash \syn{t}_1:\syn{J} \\
    \syn{x}:\syn{I},\syn{y}:\syn{J}\vdash \syn{\alpha}(\syn{x},\syn{y})\ \textsf{protype} \\
    \syn{\Gamma}_0\mid\vdash \syn{\mu}_0:\syn{\alpha}[\syn{s}_0/\syn{x}\smcl\syn{t}_0/\syn{y}] \\
    \syn{\Gamma}_m\mid\vdash \syn{\mu}_1:\syn{\alpha}[\syn{s}_1/\syn{x}\smcl\syn{t}_1/\syn{y}]\\
    \ol{\syn{\Gamma}}\mid \syn{B} \vdash \syn{i}:\syn{s}_0\ide{\syn{I}}\syn{s}_1 \\
    \ol{\syn{\Gamma}}\mid \syn{B} \vdash \syn{j}:\syn{t}_0\ide{\syn{J}}\syn{t}_1\\
    \ol{\syn{\Gamma}}\mid \syn{B} \vdash \syn{i}\boxdot\syn{\mu_1}\equiv \syn{\mu}_0\boxdot\syn{j}
    }
    {\ol{\syn{\Gamma}}\mid\syn{B}\vdash \app_{\syn{l}}(\ind_{\cmpr{}}(\syn{i},\syn{j},\syn{\mu}_0,\syn{\mu}_1))\equiv \syn{i}:\syn{s}_0\ide{\syn{I}}\syn{s}_1\\
    }
    \and
    \inferrule*[right=$\cmpr{}$-Comp]
    {\syn{\Gamma}_0\vdash \syn{s}_0:\syn{I} \\
    \syn{\Gamma}_m\vdash \syn{s}_1:\syn{I} \\
    \syn{\Gamma}_0\vdash \syn{t}_0:\syn{J} \\
    \syn{\Gamma}_m\vdash \syn{t}_1:\syn{J} \\
    \syn{x}:\syn{I},\syn{y}:\syn{J}\vdash \syn{\alpha}(\syn{x},\syn{y})\ \textsf{protype} \\
    \syn{\Gamma}_0\mid\vdash \syn{\mu}_0:\syn{\alpha}[\syn{s}_0/\syn{x}\smcl\syn{t}_0/\syn{y}] \\
    \syn{\Gamma}_m\mid\vdash \syn{\mu}_1:\syn{\alpha}[\syn{s}_1/\syn{x}\smcl\syn{t}_1/\syn{y}]\\
    \ol{\syn{\Gamma}}\mid \syn{B} \vdash \syn{i}:\syn{s}_0\ide{\syn{I}}\syn{s}_1 \\
    \ol{\syn{\Gamma}}\mid \syn{B} \vdash \syn{j}:\syn{t}_0\ide{\syn{J}}\syn{t}_1\\
    \ol{\syn{\Gamma}}\mid \syn{B} \vdash \syn{i}\boxdot\syn{\mu_1}\equiv \syn{\mu}_0\boxdot\syn{j}
    }
    {
        \ol{\syn{\Gamma}}\mid\syn{B}\vdash \app_{\syn{r}}(\ind_{\cmpr{}}(\syn{i},\syn{j},\syn{\mu}_0,\syn{\mu}_1))\equiv \syn{j}:\syn{t}_0\ide{\syn{J}}\syn{t}_1
    }
\end{mathparpagebreakable}